\RequirePackage{fix-cm}
\documentclass[smallextended]{svjour3}       

\makeatletter
\let\cl@chapter\undefined
\makeatletter

\usepackage{amsmath,amsfonts}
\usepackage{algorithmic}
\usepackage{algorithm}
\usepackage{array}
\usepackage[caption=false,font=normalsize,labelfont=sf,textfont=sf]{subfig}
\usepackage{textcomp}
\usepackage{stfloats}
\usepackage{booktabs}
\usepackage{multirow}
\usepackage{hyperref}

\usepackage{cleveref}
\usepackage{url}
\usepackage{tikz}
\usepackage{verbatim}
\usepackage{graphicx}
\usepackage{pgfplots}
\usetikzlibrary{shapes.geometric}
\usepgfplotslibrary{groupplots}
\usepackage{pgfplotstable}
\pgfplotsset{compat=newest}
\usetikzlibrary{arrows}
\usetikzlibrary{decorations}
\usetikzlibrary{decorations.markings}
\usetikzlibrary{shapes}
\usetikzlibrary{intersections}
\usetikzlibrary{patterns}

\usetikzlibrary{external}
\Crefname{figure}{Fig.}{Figs.}

\pgfplotsset{
	compat=1.11,
	legend image code/.code={
		\draw[mark repeat=2,mark phase=2]
		plot coordinates {
			(0cm,0cm)
			(0cm,0cm)        
			(0.4cm,0cm)         
		};
	}
}

\usetikzlibrary{arrows.meta, shapes.geometric}
\usepackage{cite}

\emergencystretch=\maxdimen
\newcommand{\glccluster}{\texttt{glcCluster}*}
\newcommand{\tomlab}{\texttt{TOMLAB}}
\newcommand{\dgo}{\texttt{DIRECTGO}}
\newcommand{\mic}{\texttt{MIC}}
\newcommand{\matlab}{\texttt{MATLAB}}
\newcommand{\ypea}{\texttt{YPEA}}
\newcommand{\direct}{\texttt{DIRECT}}
\newcommand{\directrev}{\texttt{DIRECT-rev}*}
\newcommand{\adc}{\texttt{ADC}}
\newcommand{\directrest}{\texttt{DIRECT-restart}}
\newcommand{\birect}{\texttt{BIRECT}}
\newcommand{\birectgb}{\texttt{Gb-BIRECT}}
\newcommand{\birmin}{\texttt{BIRMIN}*}
\newcommand{\dirmin}{\texttt{DIRMIN}*}
\newcommand{\directaggress}{\texttt{Ag. DIRECT}}
\newcommand{\directa}{\texttt{DIRECT-a}}
\newcommand{\directm}{\texttt{DIRECT-m}}
\newcommand{\directl}{\texttt{DIRECT-l}}
\newcommand{\directgl}{\texttt{DIRECT-GL}}
\newcommand{\directg}{\texttt{DIRECT-G}}
\newcommand{\plor}{\texttt{PLOR}}
\newcommand{\glbsolve}{\texttt{glbSolve}}
\newcommand{\disimplc}{\texttt{DISIMPL-C}}
\newcommand{\disimplv}{\texttt{DISIMPL-V}}
\newcommand{\directgolib}{\texttt{DIRECTGOLib v1.2}}
\newcommand{\gkls}{\texttt{GKLS}}
\newcommand{\directmro}{\texttt{MrDIRECT$_{\text{075}}$}}
\newcommand{\dbdvo}{\texttt{1-DBDV-IO}}
\newcommand{\dbdvg}{\texttt{1-DBDV-GL}}
\newcommand{\dbdva}{\texttt{1-DBDV-IA}}
\newcommand{\dtco}{\texttt{1-DTC-IO}}
\newcommand{\dtcg}{\texttt{1-DTC-GL}}
\newcommand{\dtca}{\texttt{1-DTC-IA}}

\newcommand{\surogate}{\texttt{surrogateopt}}
\newcommand{\sfla}{\texttt{SFLA}}
\newcommand{\pa}{\texttt{PSO}}
\newcommand{\ica}{\texttt{ICA}}
\newcommand{\sa}{\texttt{simulannealbnd}}
\newcommand{\be}{\texttt{BE}}
\newcommand{\avoa}{\texttt{AVOA}}
\newcommand{\amc}{\texttt{AMC}}
\newcommand{\multistart}{\texttt{MultiStart}*}
\newcommand{\aos}{\texttt{AOS}}
\newcommand{\ga}{\texttt{ga}}
\newcommand{\fa}{\texttt{FA}}
\newcommand{\iwo}{\texttt{IWO}}
\newcommand{\hs}{\texttt{HS}}

\newcommand{\bbo}{\texttt{BBO}}
\newcommand{\oqnlp}{\texttt{OQNLP}*}
\newcommand{\mqnlp}{\texttt{MSNLP}*}
\newcommand{\mcs}{\texttt{MCS}*}
\newcommand{\de}{\texttt{DE}}
\newcommand{\beso}{\texttt{BESO}}
\newcommand{\tlbo}{\texttt{TLBO}}
\newcommand{\lgo}{\texttt{LGO}*}
\newcommand{\aco}{\texttt{ACO}}
\newcommand{\bat}{\texttt{BAT}}
\newcommand{\cnna}{\texttt{CNNA}}
\newcommand{\cs}{\texttt{CS}}
\newcommand{\abc}{\texttt{ABC}}
\newcommand{\ca}{\texttt{CA}}
\newcommand{\scso}{\texttt{SCSO}}
\newcommand{\sta}{\texttt{STA}}
\newcommand{\mc}{\texttt{MC}}
\newcommand{\multimin}{\texttt{MULTIMIN}*}
\newcommand{\globalsearch}{\texttt{GlobalSearch}*}
\newcommand{\fstde}{\texttt{fstde}}
\newcommand{\nna}{\texttt{NNA}}
\newcommand{\bads}{\texttt{BADS}}
\newcommand{\StochasticRBF}{\texttt{StochasticRBF}}
\newcommand{\cmaes}{\texttt{CMA-ES}}

\newcommand{\PreserveBackslash}[1]{\let\temp=\\#1\let\\=\temp}
\newcolumntype{C}[1]{>{\PreserveBackslash\centering}p{#1}}
\newcolumntype{R}[1]{>{\PreserveBackslash\raggedleft}p{#1}}
\newcolumntype{L}[1]{>{\PreserveBackslash\raggedright}p{#1}}

\definecolor{onyx}{rgb}{0.06, 0.06, 0.06}
\definecolor{sandstorm}{rgb}{0.93, 0.84, 0.25}
\definecolor{princetonorange}{rgb}{1.0, 0.56, 0.0}
\definecolor{sienna}{rgb}{0.53, 0.18, 0.09}
\definecolor{psychedelicpurple}{rgb}{0.87, 0.0, 1.0}
\definecolor{bg}{rgb}{0.95,0.95,0.95}
\definecolor{ao}{rgb}{0.0, 0.5, 0.0}
\definecolor{arsenic}{rgb}{0.23, 0.27, 0.29}
\definecolor{armygreen}{rgb}{0.29, 0.33, 0.13}
\definecolor{antiquebrass}{rgb}{0.8, 0.58, 0.46}
\definecolor{DarkRed}{rgb}{0.55, 0.0, 0.0}
\definecolor{darkblue}{rgb}{0.0, 0.0, 0.55}
\definecolor{blueryb}{rgb}{0.01, 0.28, 1.0}
\definecolor{bluebell}{rgb}{0.64, 0.64, 0.82}
\definecolor{red}{rgb}{1.0, 0.0, 0.0}
\definecolor{redwood}{rgb}{0.67, 0.31, 0.32}
\definecolor{rose}{rgb}{1.0, 0.0, 0.5}
\definecolor{rosybrown}{rgb}{0.74, 0.56, 0.56}
\definecolor{rosewood}{rgb}{0.4, 0.0, 0.04}
\definecolor{saffron}{rgb}{0.96, 0.77, 0.19}
\definecolor{schoolbusyellow}{rgb}{1.0, 0.85, 0.0}
\definecolor{skyblue}{rgb}{0.53, 0.81, 0.92}
\definecolor{unmellowyellow}{rgb}{1.0, 1.0, 0.4}
\definecolor{wheat}{rgb}{0.96, 0.87, 0.7}
\definecolor{aureolin}{rgb}{0.99, 0.93, 0.0}
\definecolor{persianblue}{rgb}{0.11, 0.22, 0.73}
\definecolor{browna}{rgb}{0.59, 0.29, 0.0}
\definecolor{ufogreen}{rgb}{0.24, 0.82, 0.44}
\definecolor{forestgreen}{rgb}{0.13, 0.55, 0.13}
\definecolor{radicalred}{rgb}{1.0, 0.21, 0.37}
\definecolor{LightGreen}{rgb}{0.56, 0.93, 0.56}
\definecolor{LightCoral}{rgb}{0.94, 0.5, 0.5}
\definecolor{LightBlue}{rgb}{0.68, 0.85, 0.9}
\definecolor{DarkGreen}{rgb}{0.0, 0.2, 0.13}
\definecolor{LimeGreen}{rgb}{0.2, 0.8, 0.2}
\definecolor{DarkRed}{rgb}{0.55, 0.0, 0.0}
\definecolor{Tomato}{rgb}{1.0, 0.39, 0.28}
\definecolor{DarkBlue}{rgb}{0.0, 0.0, 0.55}
\definecolor{forestgreen}{rgb}{0.13, 0.55, 0.13}
\definecolor{zaffre}{rgb}{0.0, 0.08, 0.66}
\definecolor{wildstrawberry}{rgb}{1.0, 0.26, 0.64}
\definecolor{venetianred}{rgb}{0.78, 0.03, 0.08}
\definecolor{selectiveyellow}{rgb}{1.0, 0.73, 0.0}
\definecolor{yaleblue}{rgb}{0.06, 0.3, 0.57}

\pgfdeclareplotmark{aaa}{\node[star,star point ratio=2.25,minimum size=6pt,inner sep=0pt,draw=aureolin,solid, xshift=-0.5ex] {};}
\pgfdeclareplotmark{bbb}{\node[star,star points = 9,minimum size=6pt,inner sep=0pt,draw=DarkBlue,solid, xshift=-0.6ex] {};}
\pgfdeclareplotmark{ccc}{\node[isosceles triangle,draw,rotate=180,inner sep=0pt,draw=violet,solid,minimum size =5pt, xshift=-0.8ex] {};}
\pgfdeclareplotmark{ddd}{\node[isosceles triangle,draw,rotate=360,inner sep=0pt,draw=rosybrown,solid,minimum size =5pt, xshift=-0.5ex] {};}
\pgfdeclareplotmark{eee}{\node[isosceles triangle,draw,rotate=270,inner sep=0pt,draw=selectiveyellow,solid,minimum size =5pt, xshift=-1ex] {};}
\pgfdeclareplotmark{fff}{\node[star,star points = 7,minimum size=6pt,inner sep=0pt,draw=ufogreen,solid, xshift=-0.5ex] {};}
\pgfdeclareplotmark{ggg}{\node[regular polygon,regular polygon sides=5,minimum size=6pt,inner sep=0pt,draw=DarkRed,solid, xshift=-0.6ex] {};}
\pgfdeclareplotmark{hhh}{\node[regular polygon,regular polygon sides=7,minimum size=6pt,inner sep=0pt,draw=forestgreen,solid, xshift=-0.6ex] {};}
\pgfdeclareplotmark{iii}{\node[regular polygon,regular polygon sides=3,rotate=45,inner sep=0pt,draw=LightBlue,solid,minimum size =5pt, xshift=-0.4ex] {};}

\makeatletter
\tikzset{
	nomorepostactions/.code={\let\tikz@postactions=\pgfutil@empty},
	mymark/.style 2 args={decoration={markings,
			mark= between positions 0 and 1 step (1/11)*\pgfdecoratedpathlength with{%
				\tikzset{#2,every mark}\tikz@options
				\pgfuseplotmark{#1}%
			},
		},
		postaction={decorate},
		/pgfplots/legend image post style={
			mark=#1,mark options={#2},every path/.append style={nomorepostactions}
		},
	},
}
\makeatother
\smartqed  

\begin{document}

\title{An extensive numerical benchmark study of deterministic vs. stochastic derivative-free global optimization algorithms}

\titlerunning{Numerical benchmarking of deterministic vs. stochastic derivative-free algorithms}


\author{Linas Stripinis \and Remigijus Paulavi\v{c}ius}

\institute{L. Stripinis \and R. Paulavi\v{c}ius \at
	Institute of Data Science and Digital Technologies, Vilnius University, Akademijos 4, LT-08663, Vilnius, Lithuania \\
	\email{linas.stripinis@mif.vu.lt}   \\
	R. Paulavi\v{c}ius \\
	\email{remigijus.paulavicius@mif.vu.lt}
}




\date{Received: date / Accepted: date}

\maketitle

\begin{abstract}
Research in derivative-free global optimization is under active development, and many solution techniques are available today.
Therefore, the experimental comparison of previous and emerging algorithms must be kept up to date.
This paper considers the solution to the bound-constrained, possibly black-box global optimization problem.
It compares $64$ derivative-free deterministic algorithms against classic and state-of-the-art stochastic solvers.
Among deterministic ones, a particular emphasis is on \direct-type, where, in recent years, significant progress has been made.
A set of $800$ test problems generated by the well-known \gkls{} generator and $397$ traditional test problems from \directgolib{} collection are utilized in a computational study.
More than $239,400$ solver runs were carried out, requiring more than $531$ days of single CPU time to complete them.
It has been found that deterministic algorithms perform excellently on \texttt{GKLS}-type and low-dimensional problems, while stochastic algorithms have shown to be more efficient in higher dimensions.
\end{abstract}

\section{Introduction}

Optimization methods are continuously used today to improve and optimize business processes in various areas of human activity and industry.
In this paper, we consider a general single-objective optimization problem, which can be formally stated as:
\begin{equation}
	\label{eq:opt-problem}
	\begin{aligned}
		& \min_{\mathbf{x}\in D} && f(\mathbf{x})
	\end{aligned}
\end{equation}
where $f:\mathbb{R}^n \rightarrow \mathbb{R}$ is a Lipschitz-continuous objective function (with an unknown Lipschitz constant), and $\mathbf{x} \in \mathbb{R}^n$ is the input vector of control variables.
We assume that $f$ can be computed at any point of the feasible region, which is an $n$-dimensional hyper-rectangle
\[
D = [ \mathbf{a},  \mathbf{b}] = \{ \mathbf{x} \in \mathbb{R}^n: a_j \leq x_j \leq b_j, j = 1, \dots, n\}.
\]
Moreover, $f$ can be non-linear, multi-modal, non-convex, and non-differentiable.

It is common practice in many fields of science and engineering to optimize a function derived from an experiment or a highly complex computer simulation.
In addition, as the scale and complexity of applications increase, objective function evaluations become more expensive.
Such situations make derivatives either impossible or impractical to calculate. 
Thus, this paper also assumes that the objective function is potentially ``black-box.''
Therefore, any analytical information about the objective function is unknown, and optimization solvers can only use the objective function values.

Numerical methods are commonly applied to solve ``black-box'' problems~\cite{Zhigljavsky2008:book,William2007}.
The development of derivative-free optimization algorithms has a long history in the field of optimization, dating back to the work of Hooke and Jeeves \cite{Hooke1961} with their famous ``direct search'' approach.
After several years of work, derivative-free techniques have lost popularity in the mathematical optimization community, as some of them converge  slowly or even can not guarantee global convergence~\cite{Xi2020}.
However, in recent decades, derivative-free optimization has regained interest due to the growing number of applications and has attracted much attention from the optimization community \cite{Ezugwu2021,Jones2021}.
Our collected data from the Web of Science (WoS) shows more than $1800$ publications related to derivative-free optimization (see \Cref{tab:document}).

\begin{table}[ht]
	\centering
	\caption{The number of publications related to derivative-free optimization (in WoS database).}
	\begin{tabular*}{\textwidth}{@{\extracolsep{\fill}}lrr}
		\toprule
		Document Type	&	Record Count	&	Percentage	\\
		\midrule
		Articles		&	$1,064$	&	$57.64\%$	\\
		Proceedings Paper	&	$725$	&	$39.27\%$	\\
		Review Articles	&	$27$	&	$1.46\%$	\\
		Book Chapters	&	$29$	&	$1.57\%$	\\
		Book				&	$1$		&	$0.05\%$	\\
		\bottomrule
	\end{tabular*}
	\label{tab:document}
\end{table}

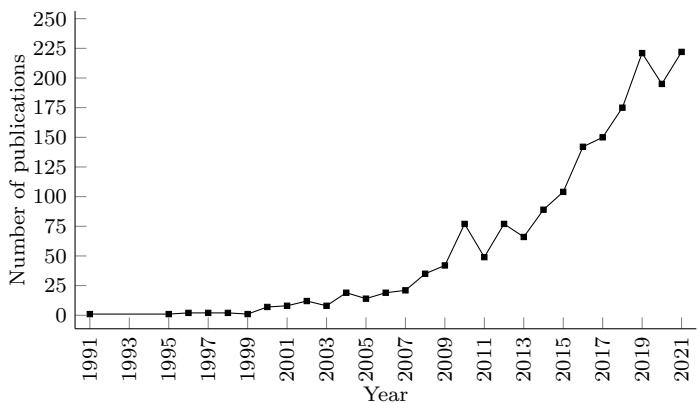
\begin{figure}[ht]
	\centering
	\begin{tikzpicture}
		\begin{axis}[
			xlabel = {Year},
			axis x line*=bottom,
			axis y line*=left,
			xmin=1991,ymax=2021,
			xtick distance=2,
			ymin=-5,ymax=250,
			ytick distance=25,
			xtick = {1991,1993,1995,1997,1999,2001,2003,2005,2007,2009,2011,2013,2015,2017,2019,2021},
			xticklabels = {1991,1993,1995,1997,1999,2001,2003,2005,2007,2009,2011,2013,2015,2017,2019,2021},
			ylabel = {Number of publications},
			xticklabel style = {rotate=90,anchor=east},
			ylabel style={yshift=-0.5em},
			xlabel style={yshift=0.5em},
			legend style={font=\footnotesize},
			legend pos= north west,
			height=0.475\textwidth,width=0.8\textwidth,
			enlargelimits=0.025,
			]
			\addplot[black,mark=square*,mark size=1pt] table[x=F,y=D] {data/direct.txt};
			
		\end{axis}
	\end{tikzpicture}
	\caption{The number of publications related to derivative-free optimization from $1990$ to $2021$ (in WoS database).}
	\label{fig:citations}
\end{figure}

The cumulative number of publications from $1991$ to $2021$ is plotted in \Cref{fig:citations}.
The number of ``derivative-free'' related publications grows each year consistently.
Naturally, the number will be significantly higher in some other databases such as Google Scholar.

Choosing the proper optimization algorithm (solver) for a solution to global optimization problems is complicated.
First, the user needs to know how different solvers perform compared to other algorithms on an extensive set of different test and possibly (related) practical problems.
The choice is further complicated because global optimization techniques are different.
Stochastic nature-inspired global optimization algorithms incorporate probabilistic (stochastic) elements.
The convergence is guaranteed only in a probabilistic sense, i.e., in infinite time, a global optimum will be found with probability one \cite{Zhigljavsky2008:book}.
Deterministic algorithms can guarantee that after a finite time a global optimum will be found within small prescribed tolerances~\cite{Floudas1999book,Horst1996:book}.

Identification of the most effective derivative-free global optimization algorithms is a challenging task.
Typically, the authors of new algorithms limit themselves to small efficiency studies against only a few related algorithms (see, e.g., \cite{Abdollahzadeh2021,Alsattar2020,Azizi2021,Stripinis2018a,Liu2015}).
Therefore, the complete picture of the algorithm's performance across state-of-the-art algorithms is unclear.
Unfortunately, detailed comparative analyzes in the field are rare, and there is a constant need to update them by including recent algorithms.

\subsection{Related surveys and benchmarks}\label{sec:review1}

\Cref{tab:comparison} summarizes related works devoted to solving derivative-free global optimization problems.
Here, the first three columns give the source of the research, the year of publication, and the impact expressed as the total number of citations (according to Google Scholar).
The fourth column indicates the total number of solvers, while the fifth shows the total number of test problems used in the experimental studies.
The last two columns give the range of the dimension and the function evaluation budget.

Most published surveys (see, e.g., \cite{Mongeau2000,Auger2009,Pham2011,Stripinis2021c,Kvasov2018,Liu2019}) are limited to relatively small-dimensional test problems.
The largest number of derivative-free global optimization solvers was investigated in our recent publication \cite{Stripinis2021c}, where the new \dgo{} toolbox was introduced.
Nevertheless, this paper considered only \direct-type~\cite{Jones1993} algorithms.
However, the most extensive study of derivative-free solvers was carried out in~\cite{Rios2013} and currently has the most significant impact.
The study was carried out using only a limited budget of objective function evaluations, typical for expensive global optimization problems.
Since then, the number of new derivative-free optimization solvers significantly increased, and it is not clear whether the best-performing algorithms are still the best?


\begin{table}[ht]
	\centering
	\caption{History of derivative-free global optimization numerical studies}
	\resizebox{\textwidth}{!}{
		\begin{tabular}{lrccccc}
			\toprule
			\multirow{2}{*}{Source} & \multirow{2}{*}{Year} & \multirow{2}{*}{Impact} & \multicolumn{2}{c}{Number of} & Range & Evaluation \\
			\cmidrule{4-5}
			&  &  & solvers & problems & of $n$  & budget \\
			\midrule
			Mongeau et al. \cite{Mongeau2000}       & $2000$ & $80$  & $6$  & $11$  & $1-20$  & $700-30,000$ \\
			Hansen et al. \cite{Hansen2010}         & $2010$ & $411$  & $31$ & $120$  & $2-40$  & $100,000,000n$ \\
			Auger et al. \cite{Auger2009}           & $2011$ & $98$  & $5$  & $3$   & $10-40$ & $10,000,000$  \\
			Pham et al. \cite{Pham2011}             & $2011$ & $58$  & $5$  & $16$  & $20$    & $100,000$  \\
			Rios et al. \cite{Rios2013}             & $2013$ & $1,260$ & $22$ & $502$ & $1-300$ & $2,500$ \\
			Kvasov et al. \cite{Kvasov2018}         & $2018$ & $78$  & $13$ & $134$ & $1$     & $10,000$ \\
			Liu et al. \cite{Liu2019}  			    & $2019$ & $17$   & $27$ & $126$ & $1$     & $1,000$ \\
			Stripinis et al. \cite{Stripinis2021c}  & $2022$ & $1$   & $39$ & $128$ & $2-15$  & $2,000,000$ \\
			\midrule
			This paper  							& $2022$ & $0$   & $64$ & $1,197$ & $2-100$ & $2,500-500,000$ \\
			\bottomrule
		\end{tabular}
	}
	\label{tab:comparison}
\end{table}

\subsection{Main contributions}\label{sec:contributions}
The main contributions of this research are the following:
\begin{itemize}
	\item[1.] The largest set of solvers ($64$) considered in the experimental study of derivative-free global optimization.
	\item[2.] Large and diverse benchmark pool ($1,197$) consisting of:
	\begin{itemize}
		\item[a.] The \gkls{}-type multi-modal box-constrained ``black-box'' global optimization problems~\cite{Gaviano2003} (recently  extended to generally-constrained problems~\cite{Sergeyev2021gkls}).
		\item[b.] Test and engineering global optimization problems from \directgolib~\cite{DIRECTGOLib2022}.
	\end{itemize}
	\item[3.] All solvers tested using two different stopping conditions:
	\begin{itemize}
		\item[a.] Assuming expensive global optimization and limiting function evaluation budget to $2,500$.
		\item[b.] Assuming cheap problems and using $500,000$ function evaluation budget.
	\end{itemize}
	\item[4.] Significantly expanded \directgolib~\cite{DIRECTGOLib2022}.
\end{itemize}

The rest of the paper is organized as follows.
A brief review of considered algorithms is described in~\Cref{sec:review}.
\Cref{sec:comparison} describes the experimental design.
Extensive experimental results are reported in \cref{sec:benchmarking}, while \cref{sec:results} summarizes the performance of the $64$ considered algorithms.
\Cref{sec:conclusion} concludes the paper.

\section{Brief review of considered algorithms}\label{sec:review}

This section briefly reviews derivative-free global optimization techniques considered in this study.

\subsection{Deterministic global search algorithms}\label{sec:deterministic}
In this subsection, we briefly review classes of deterministic derivative-free algorithms from which algorithms were selected for this study (see~\Cref{tab:direct}).
In \Cref{tab:direct}, the first two columns give a short algorithm's description and source.
The third column indicates the year of publication, while the fourth shows the total number of citations so far in Google Scholar.
The fifth column indicates the class of the algorithm, while the sixth shows the acronym used in this paper.
The last two columns give the source for the implementation and the exact version used in this study.



\subsubsection{Branch-and-bound (BB) algorithms}
Branch-and-bound algorithms sequentially partition the optimization domain and determine lower and upper bounds for the optimum.
For example, Lipschitzian-based approaches (see, e.g., \cite{Shubert1972,Pinter1996book,Paulavicius2006,Paulavicius2007}) construct and optimize a function that underestimates the original one in a piecewise manner.
The major drawback of the most traditional Lipschitz optimization algorithms is the requirement of knowing the Lipschitz constant.
There exist solvers operating with an a priori given an estimate of the Lipschitz constant \cite{Sergeyev1998a}, its adaptive estimates \cite{GERGEL1997,GERGEL1999163,Sergeyev1998a}, and adaptive estimates of local Lipschitz constants \cite{Sergeyev1998a,Sergeyev2008:book}.
Among all algorithms considered in this study, only the \lgo{} algorithm belongs to this class (see~\Cref{tab:direct}).

\subsubsection{\direct-type algorithms}
The \direct{} algorithm developed by Jones et al.~\cite{Jones1993} extends classical Lipschitz optimization~\cite{Paulavicius2006,Paulavicius2007,Paulavicius2008,Paulavicius2009b,Pinter1996book,Piyavskii1967,Sergeyev2011,Shubert1972}, where the need for the Lipschitz constant is eliminated.
This feature made \direct-type algorithms especially attractive for solving various real-world optimization problems~(see, e.g., \cite{Baker2000,Bartholomew2002,Carter2001,Cox2001,Serafino2011,Gablonsky2001,Liuzzi2010,Paulavicius2019:eswa,Paulavicius2014:book,Stripinis2018b} and the references given therein).

\begin{figure}[htbp]
	\centering
	\resizebox{\textwidth}{!}{
		\begin{tikzpicture}[node distance=1.5cm]
			\tikzstyle{startstop} = [rectangle, rounded corners, text width=7em, minimum width=2cm, minimum height=1cm,text centered, draw=blue!60, fill=blue!60, line width=0.4mm]
			\tikzstyle{procesion}   = [rectangle, text width=15em, minimum width=3cm, minimum height=1cm, text centered, draw=DarkGreen!70, fill=DarkGreen!70, line width=0.4mm]
			\tikzstyle{decision}  = [diamond, aspect=2, inner sep=0pt, text width=5em, text centered, draw=red!80, fill=red!80, line width=0.4mm]
			
			\node (start)   [startstop] {\textcolor{white}{\textbf{Start}}};
			\node (Phase1)  [procesion, below of=start, yshift=-0.25cm]   {\textcolor{white}{\textbf{Initialization}}};
			\node (Phase2)  [procesion, below of=Phase1, yshift=-0.25cm]   {\textcolor{white}{\textbf{Selection}: identify the most promising regions in the domain}};
			\node (Phase3)  [procesion, below of=Phase2, yshift=-0.25cm]   {\textcolor{white}{\textbf{Sampling}: evaluate objective function in selected regions}};
			\node (Phase4)  [procesion, below of=Phase3, yshift=-0.25cm]   {\textcolor{white}{\textbf{Subdivision}: divide selected regions}};
			\node (term)    [startstop, left of=start, xshift=-3.5cm] {\textcolor{white}{\textbf{Termination}}};
			\node (stop)    [decision, below of=term, yshift=-3.25cm] {\textcolor{white}{\textbf{Stop?}}};
			\node (ii)      [coordinate, left of=Phase2, xshift=-3.5cm] {};
			
			\draw[->, line width=0.5mm, draw=black!50] (start)  -- (Phase1);
			\draw[->, line width=0.5mm, draw=black!50] (Phase1) -- (Phase2);
			\draw[->, line width=0.5mm, draw=black!50] (Phase2) -- (Phase3);
			\draw[->, line width=0.5mm, draw=black!50] (Phase3) -- (Phase4);
			\draw[->, line width=0.5mm, draw=black!50] (Phase4) -| (stop);
			\draw[->, line width=0.5mm, draw=black!50] (stop)   -- node[anchor=east] {\textcolor{black!50}{\textbf{Yes}}} (term);
			\draw[-, line width=0.5mm, draw=black!50]  (stop)   -- (ii);
			\draw[->, line width=0.5mm, draw=black!50] (ii)     -- node[anchor=north] {\textcolor{black!50}{\textbf{No}}} (Phase2);
		\end{tikzpicture}
		\hspace{1cm}
		\begin{tikzpicture}[node distance=1.5cm]
			
			\tikzstyle{process}   = [rectangle, text width=5em, minimum width=3cm, minimum height=1cm, text centered, draw=DarkGreen!70, fill=DarkGreen!70, line width=0.4mm]
			\tikzstyle{startstop} = [rectangle, rounded corners, text width=7em, minimum width=2cm, minimum height=1cm,text centered, draw=blue!60, fill=blue!60, line width=0.4mm]
			\tikzstyle{procesion}   = [rectangle, text width=15em, minimum width=3cm, minimum height=1cm, text centered, draw=DarkGreen!70, fill=DarkGreen!70, line width=0.4mm]
			\tikzstyle{decision}  = [diamond, aspect=2, inner sep=0pt, text width=5em, text centered, draw=red!80, fill=red!80, line width=0.4mm]
			
			\node (start)   [startstop] {\textcolor{white}{\textbf{Start}}};
			\node (Phase1)  [procesion, below of=start, yshift=-0.25cm]   {\textcolor{white}{Generate initial popolation}};
			\node (Phase2)  [procesion, below of=Phase1, yshift=-0.25cm]  {\textcolor{white}{Evaluate population}};
			\node (term)    [startstop, left of=start, xshift=-3.5cm] {\textcolor{white}{\textbf{Termination}}};
			\node (stop)    [decision, below of=term, yshift=-1.85cm] {\textcolor{white}{\textbf{Stop?}}};
			\node (Phase3)  [process, below of=stop, yshift=-0.5cm]   {\textcolor{white}{\textbf{Selection}}};
			\node (Phase4)  [process, below of=Phase3, yshift=-0.25cm]   {\textcolor{white}{\textbf{Crossover, Mutation}}};
			\node (Phase5)  [procesion, right of=Phase4, xshift=3.5cm]   {\textcolor{white}{Generate new population}};
			
			\draw[->, line width=0.5mm, draw=black!50] (start)  -- (Phase1);
			\draw[->, line width=0.5mm, draw=black!50] (Phase1) -- (Phase2);
			\draw[->, line width=0.5mm, draw=black!50] (stop)   -- node[anchor=east] {\textcolor{black!50}{\textbf{Yes}}} (term);
			\draw[->, line width=0.5mm, draw=black!50] (Phase2)   --  (stop);
			\draw[->, line width=0.5mm, draw=black!50] (stop) -- node[anchor=east] {\textcolor{black!50}{\textbf{No}}} (Phase3);
			\draw[->, line width=0.5mm, draw=black!50] (Phase3)   --  (Phase4);
			\draw[->, line width=0.5mm, draw=black!50] (Phase4)   --  (Phase5);
			\draw[->, line width=0.5mm, draw=black!50] (Phase5)   --  (Phase2);
		\end{tikzpicture}
	}
	\caption{The basic structures of \direct{}-type (left side) and  population-based (right side) algorithms.}
	\label{fig:flowchart_direct}
\end{figure}
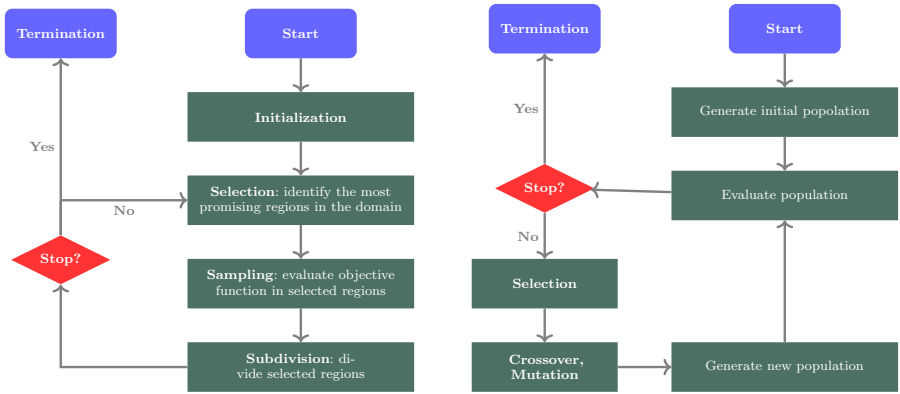

Many \direct-type algorithms have been proposed, but almost all share the same basic structure (see left side of \Cref{fig:flowchart_direct}).
The main three steps are selection, sampling, and partitioning.
First, a specific \direct-type algorithm identifies the so-called potentially optimal hyper-rectangles (POHs).
Then explores these POHs by performing new samples and subdividing them using some partitioning scheme.

Among various proposals, lots of attention was paid to improving the selection of POHs (see, e.g., ~\cite{Baker2000,Gablonsky2001:phd,Mockus2017,Paulavicius2019:eswa,Stripinis2018a}).
Other authors (see, e.g.,~\cite{Jones2001,Liu2015b,Paulavicius2016:jogo,Paulavicius2013:jogo,Sergeyev2006}) have shown that applying different partitioning techniques can have a positive impact on the performance.
Our recent analysis in \cite{Stripinis2021b} revealed that better \direct-type algorithms are obtained by creating new combinations of existing  selection and partitioning schemes.
For a detailed, comprehensive descriptions and review of \direct-type algorithms published in the last 25 years, we refer to\cite{Jones2021}.
In total, $26$ \direct-type algorithms were chosen in this study (see \Cref{tab:direct}).


\begin{table*}
	\caption{List of the deterministic algorithms (their implementations) in chronological order starting from the oldest.}
	\resizebox{\textwidth}{!}{
		\begin{tabular}{lccccccc}
			\toprule
			Algorithm & Reference & Year & Impact 					& Class			& Acronym     & Implementation & Version \\
			\midrule
			The original \direct & \cite{Jones1993} & $1993$ & $2,366$ & \direct-type & \direct & \dgo~\cite{DIRECTGO2022} & $1.1.0$ \\
			Multilevel coordinate search & \cite{Huyer1999} & $1999$ & $658$ & MCS & \mcs & \mcs~\cite{Arnold1999} & $2.0$ \\
			Revised version of the \direct & \cite{Jones2001} & $2001$ & $401$ & \direct-type & \directrev & \dgo~\cite{DIRECTGO2022} & $1.1.0$ \\
			Locally-biased form of the \direct & \cite{Gablonsky2001} & $2001$ & $365$ & \direct-type & \directl & \dgo~\cite{DIRECTGO2022} & $1.1.0$ \\
			Aggressive version of the \direct & \cite{Baker2000} & $2001$ & $73$ & \direct-type & \directaggress & \dgo~\cite{DIRECTGO2022} & $1.1.0$ \\
			Adaptive setting of the \direct{} parameters & \cite{Finkel2004aa} & $2004$ & $20$ & \direct-type & \directrest & \dgo~\cite{DIRECTGO2022} & $1.1.0$ \\
			Clustering technique-based hybridized \direct & \cite{Holmstrom2004} & $2004$ & $34$ & \direct-type & \glccluster & \tomlab~\cite{Holmstrom2010} & $8.7$ \\
			Implementation of original \direct{} algorithm & \cite{Holmstrom2004} & $2004$ & $34$ & \direct-type & \glbsolve & \tomlab~\cite{Holmstrom2010} & $8.7$ \\
			Multi-start scatter search algorithm & \cite{Holmstrom2004} & $2004$ & $34$ & MS & \oqnlp & \tomlab~\cite{Holmstrom2010} & $8.7$ \\
			Multi-start scatter search algorithm & \cite{Holmstrom2004} & $2004$ & $34$ & MS & \mqnlp & \tomlab~\cite{Holmstrom2010} & $8.7$ \\
			Hybridized branch-and-bound algorithm & \cite{Holmstrom2004} & $2004$ & $34$ & BB & \lgo & \tomlab~\cite{Holmstrom2010} & $8.7$ \\
			Additive scaling in \direct{} algorithm & \cite{Finkel2006} & $2006$ & $64$ & \direct-type & \directm & \dgo~\cite{DIRECTGO2022} & $1.1.0$ \\
			Adaptive diagonal curves & \cite{Sergeyev2006} & $2006$ & $118$ & \direct-type & \adc & \dgo~\cite{DIRECTGO2022} & $1.1.0$ \\
			Aggressively hybridized version of the \direct & \cite{Liuzzi2010} & $2010$ & $64$ & \direct-type & \dirmin & \dgo~\cite{DIRECTGO2022} & $1.1.0$ \\
			Linear scaling and the \direct{} algorithm & \cite{Liu2013} & $2013$ & $17$ & \direct-type & \directa & \dgo~\cite{DIRECTGO2022} & $1.1.0$ \\
			Simplicial partitioning technique based \direct & \cite{Paulavicius2013:jogo} & $2013$ & $58$ & \direct-type & \disimplc & \dgo~\cite{DIRECTGO2022} & $1.1.0$ \\
			Simplicial partitioning technique based \direct & \cite{Paulavicius2013:jogo} & $2013$ & $58$ & \direct-type & \disimplv & \dgo~\cite{DIRECTGO2022} & $1.1.0$ \\
			Multilevel \direct{} algorithm & \cite{Liu2015} & $2015$ & $30$ & \direct-type & \directmro & \dgo~\cite{DIRECTGO2022} & $1.1.0$ \\
			\direct{} using bisection of hyper-rectangles & \cite{Paulavicius2016:jogo} & $2016$ & $20$ & \direct-type & \birect & \dgo~\cite{DIRECTGO2022} & $1.1.0$ \\
			Pareto-Lipschitzian optimization & \cite{Mockus2017} & $2017$ & $16$ & \direct-type & \plor & \dgo~\cite{DIRECTGO2022} & $1.1.0$ \\
			Pareto selection based \direct & \cite{Stripinis2018a} & $2018$ & $16$ & \direct-type & \directg & \dgo~\cite{DIRECTGO2022} & $1.1.0$ \\
			Two-step (Global-Local), selection based \direct & \cite{Stripinis2018a} & $2018$ & $16$ & \direct-type & \directgl & \dgo~\cite{DIRECTGO2022} & $1.1.0$ \\
			Globally-biased hybridized version of the \birect & \cite{Paulavicius2019:eswa} & $2019$ & $30$ & \direct-type & \birmin & \dgo~\cite{DIRECTGO2022} & $1.1.0$ \\
			Globally-biased version of the \birect & \cite{Paulavicius2019:eswa} & $2019$ & $30$ & \direct-type & \birectgb & \dgo~\cite{DIRECTGO2022} & $1.1.0$ \\
			Improved aggressive version of the \directrev & \cite{Stripinis2021b} & $2022$ & $0$ & \direct-type & \dtca & \dgo~\cite{DIRECTGO2022} & $1.1.0$ \\
			Improved aggressive version of the \birect & \cite{Stripinis2021b} & $2022$ & $0$ & \direct-type & \dbdva & \dgo~\cite{DIRECTGO2022} & $1.1.0$ \\
			Improved original selection based \directrev & \cite{Stripinis2021b} & $2022$ & $0$ & \direct-type & \dtco & \dgo~\cite{DIRECTGO2022} & $1.1.0$ \\
			Improved original selection based \birect & \cite{Stripinis2021b} & $2022$ & $0$ & \direct-type & \dbdvo & \dgo~\cite{DIRECTGO2022} & $1.1.0$ \\
			Two-step (Global-Local), selection based \directrev & \cite{Stripinis2021b} & $2022$ & $0$ & \direct-type & \dtcg & \dgo~\cite{DIRECTGO2022} & $1.1.0$ \\
			Two-step (Global-Local), selection based \birect & \cite{Stripinis2021b} & $2022$ & $0$ & \direct-type & \dbdvg & \dgo~\cite{DIRECTGO2022} & $1.1.0$ \\
			\bottomrule
			\multicolumn{8}{l}{* -- Algorithm is hybridized with a local search procedures.}
	\end{tabular}}
	\label{tab:direct}
\end{table*}


\subsubsection{Multilevel coordinate search (MCS) algorithms}\label{sec:multilevel}

Similar to \direct{}, \texttt{MCS} partitions the search space into smaller hyper-rectangles.
Each contains a distinguished point, the so-called base point, in which the objective function value is evaluated.
The hyper-rectangular partitioning procedure is not uniform, and regions are preferred where better function values are expected to be found.
Same as in \direct, \texttt{MCS} solver typically combine global search (a subdivision of the largest hyper-rectangles) and local search (a subdivision of the hyper-rectangles with good function values).
The \mcs{} algorithm carries the multilevel technique for balancing global and local search in each iteration.
\mcs{} was the only one included in our study of this class (see~\Cref{tab:direct}).

\subsubsection{Multi-start (MS) search algorithms}\label{sec:multi}

The multi-start technique \cite{Hey1979} is a well-known stochastic approach that attempts to run local search procedures from a set of random starting points. 
When the initial points are chosen in a non-random way, the algorithm acquires the property of determinism.
From this class, we consider two solvers (\oqnlp{} and \mqnlp{}), designed to find global optima of smooth constrained and mixed-integer nonlinear problems.

To summarize, the main advantages of deterministic techniques are the guarantee of the same result (output) and a typically low number of algorithmic input parameters~\cite{Jones2021}.
The main disadvantage is the limited eﬃciently on higher-dimensional global optimization problems~\cite{Jones2021,Rios2013}.

\subsection{Stochastic global search algorithms}\label{sec:stochastics}


Due to the non-deterministic nature, the guarantee of the optimal solution in stochastic algorithms is provided only in a probabilistic sense.
The convergence of most stochastic techniques is based on the classical probability theory, ``zero-one law.''
In the following subsections, we will review several well-known classes of stochastic algorithms from which algorithms were chosen for this study (see \Cref{tab:EAs}).
The structure of \Cref{tab:direct} is the same as for \Cref{tab:EAs}.
The only difference is that in the last column, we specify the nature of the algorithm.

\subsubsection{Pure Random Search (PRS)}\label{sec:PRS}

The most straightforward stochastic global optimization is pure random search (PRS), implemented as Monte-Carlo algorithm (see \Cref{tab:EAs}).
A set of uniformly distributed random points is generated in every iteration, and the objective function values at these points are recorded.
The generation of the new points is constructed independently of previous sampled points, i.e., PRS does not use information from previous iterations.
For this reason, the performance of the PRS technique is considered unrelated to the structure of the objective function.

\subsubsection{Markovian Global Search (MGS)}\label{sec:MGS}

MGS algorithms are yet simple but cleverer than PRS.
The distribution of the points in MGS depends on the previously sampled point and the current best-recorded function value.
The main criticism of MGS class algorithms is that they are often too myopic and do not make efficient use of information about the objective function achieved earlier.
The most popular in this class is the well-known simulated annealing algorithm, which also has the highest impact score among all the considered algorithms (see \Cref{tab:EAs}).

\begin{table*}
	\caption{List of the stochastic algorithms (their implementations) in chronological order starting from the oldest.}
	\resizebox{\textwidth}{!}{
		\begin{tabular}{lcccccccc}
			\toprule
			Algorithm & Reference  & Year & Impact & Class & Acronym & Implementation  & Version & Based on  \\
			\midrule
			Monte-Carlo algorithm & \cite{Metropolis1949} & $1955$ & $8,117$ & PRS & \mc & \mic~\cite{Tche2022} & $1.0$ & Random \\
			Genetic algorithm & \cite{Holland1975} & $1975$ & $38,629$ & PBS & \ga & \matlab~\cite{MATLAB2022} & $9.12.0$ & Evolutionary \\ 
			Simulated annealing & \cite{Kirkpatrick1983a} & $1983$ & $53,379$ & MGS & \sa & \matlab~\cite{MATLAB2022} & $9.12.0$ & Nature \\
			Adaptive Monte-Carlo algorithm & \cite{Patel1989} & $1989$ & $112$ & PAS & \amc & \mic~\cite{Tche2022} & $1.0$ & Random \\
			Ant colony optimization & \cite{Dorigo1992O} & $1992$ & $16,032$ & PBS & \aco & \ypea~\cite{Kalami2020} & $1.1.0.4$ & Swarm \\
			Cultural algorithm & \cite{Reynolds2008} & $1994$ & $1,258$ & PBS & \ca & \ypea~\cite{Kalami2020} & $1.1.0.4$ & Evolutionary \\
			Particle swarm optimization & \cite{Eberhart1995} & $1995$ & $17,922$ & PBS & \pa & \ypea~\cite{Kalami2020} & $1.1.0.4$ & Swarm \\
			Differential evolution & \cite{Storn1997} & $1997$ & $29,527$ & PBS & \de & \ypea~\cite{Kalami2020} & $1.1.0.4$ & Population \\
			Surrogate optimization algorithm & \cite{Gutmann2001} & $2001$ & $863$ & MBO & \surogate & \matlab~\cite{MATLAB2022} & $9.12.0$ & Model-based \\
			Harmony search & \cite{Geem2001} & $2001$ & $6,457$ & PBS & \hs & \ypea~\cite{Kalami2020} & $1.1.0.4$ & Nature \\
			Shuffled Frog Leaping Algorithm & \cite{Muzaffar2003} & $2003$ & $1,850$ & PBS & \sfla & \sfla~\cite{Yarpiz2022} & $1.0.0.0$ & Bio \\
			Multi-start global random search & \cite{Holmstrom2004} & $2004$ & $34$ & RMS & \multimin & \tomlab~\cite{Holmstrom2010} & $8.7$ & Random \\
			Artificial bee colony & \cite{Karaboga2005ANIB} & $2005$ & $7,669$ & PBS & \abc & \ypea~\cite{Kalami2020} & $1.1.0.4$ & Swarm \\
			Bees algorithm & \cite{Pham2006} & $2006$ & $1,535$ & PBS & \be & \ypea~\cite{Kalami2020} & $1.1.0.4$ & Swarm \\
			Invasive Weed Optimization & \cite{Mehrabian2006} & $2006$ & $1,441$ & PBS & \iwo & \ypea~\cite{Kalami2020} & $1.1.0.4$ & Swarm \\
			Imperialist Competitive Algorithm & \cite{Gargari2007} & $2007$ & $2,912$ & PBS & \ica & \ypea~\cite{Kalami2020} & $1.1.0.4$ & Evolutionary \\
			Multi-start global random search & \cite{Ugray2007} & $2007$ & $697$ & RMS & \multistart & \matlab~\cite{MATLAB2022} & $9.12.0$ & Random \\
			Multi-start global random search & \cite{Ugray2007} & $2007$ & $697$ & RMS & \globalsearch & \matlab~\cite{MATLAB2022} & $9.12.0$ & Random \\
			Stochastic Radial Basis Function & \cite{Regis2007} & $2007$ & $422$ & MBO & \StochasticRBF & \StochasticRBF~\cite{Julie2006} & $1.0.0.0$ & Model-based \\
			Biogeography-based Optimization & \cite{Dan2008} & $2008$ & $3,659$ & PBS & \bbo & \ypea~\cite{Kalami2020} & $1.1.0.4$ & Nature \\
			Firefly algorithm & \cite{Yang2009as} & $2009$ & $4,306$ & PBS & \fa & \ypea~\cite{Kalami2020} & $1.1.0.4$ & Swarm \\
			Cuckoo search & \cite{Yang2009a} & $2009$ & $6,772$ & PBS & \cs & \cs~\cite{Yang2022a} & $2.0.0$ & Swarm \\
			Bat-Inspired Algorithm & \cite{Yang2010as} & $2010$ & $5,179$ & PBS & \bat & \bat~\cite{Yang2022b} & $2.0.0$ & Swarm \\
			Covariance matrix adaptation strategy & \cite{Iruthayarajan2010} & $2010$ & $72$ & PBS & \cmaes & \ypea~\cite{Kalami2020} & $2.0.0$ & Evolutionary \\
			Teaching-Learning-based Optimization & \cite{Rao2011} & $2011$ & $3,307$ & PBS & \tlbo & \ypea~\cite{Kalami2020} & $1.1.0.4$ & Population \\
			Improved State Transition Algorithm & \cite{Saravanakumar2015} & $2015$ & $4$ & PBS & \sta & \sta~\cite{sta2022} & $1.0$ & Population \\
			Bayesian Adaptive Direct Search & \cite{Acerbi2017} & $2017$ & $152$ & MBO & \bads & \bads~\cite{Acerbi2022} & $1.0.8$ & Model-based \\
			Chaotic neural network algorithm & \cite{Sadollah2018} & $2018$ & $125$ & PBS & \cnna & \cnna~\cite{Zhang2022} & $1.0.0$ & Bio \\
			Neural Network Algorithm & \cite{Sadollah2018dynamic} & $2018$ & $130$ & PBS & \nna & \nna~\cite{Sadollah2018dynamic} & $1.0.0$ & Bio \\
			Bald Eagle search optimization & \cite{Alsattar2020} & $2020$ & $100$ & PBS & \beso & \beso~\cite{Hassan2022} & $1.0.0$ & Nature \\
			Fuzzy self-tuning differential evolution & \cite{Tsafarakis2020} & $2020$ & $19$ & PBS & \fstde & \fstde~\cite{Tsafarakis2020} & $2.0.0$ & Population \\
			African Vulture Optimization algorithm & \cite{Abdollahzadeh2021} & $2021$ & $101$ & PBS & \avoa & \avoa~\cite{Abdollahzadeh2021} & $1.0.2$ & Nature \\
			Atomic Orbital Search & \cite{Azizi2021} & $2021$ & $33$ & PBS & \aos & \aos~\cite{Azizi2022} & $1.0.0$ & Quantum \\
			Sand Cat swarm optimization & \cite{Seyyedabbasi2022} & $2022$ & $0$ & PBS & \scso & \scso~\cite{amir2022} & $1.0.0$ & Nature \\
			\bottomrule
			\multicolumn{9}{l}{* -- Algorithm is hybridized with local search procedures.}
	\end{tabular}}
	\label{tab:EAs}
\end{table*}

\subsubsection{Pure Adaptive Search (PAS)}\label{sec:PAS}

The PAS approach starts by generating a point uniformly distributed within the optimization domain, where the objective function value is recorded.
The next point is generated from a uniform distribution over the region constructed by the intersection of the whole optimization region with the open level set of points with objective function values less than the best-recorded value.
The approach proceeds iteratively in this manner until some stopping criterion is satisfied.
This way, PAS performs samples in that part of the optimization domain that gives a strictly improving objective function value at each iteration.
PAS technique has been implemented in the adaptive Monte-Carlo algorithm (see \Cref{tab:EAs}).

\subsubsection{Random Multi-Start (RMS)}\label{sec:RMS}

The multi-start calls a local search procedure with different initial values and returns the solutions and the optimal point from each started point.
Random multi-start (RMS) is a well-known and well known global optimization algorithm in practice.
Local search procedures are performed from several random points in the optimization domain.
The point can be generated using PRS (as in the \multimin, \multistart, and \globalsearch{} algorithms) or PAS techniques.

\subsubsection{Population-Based Search (PBS)}\label{sec:meta} 

The population (set of points) evolves in PBS algorithms rather than single points.
There are many publications on PBS, the majority of which deal with meta-heuristics rather than theory and generic methodology.
The term ``meta-heuristic'' has been known for almost four decades~\cite{GLOVER1986533}.
Meta-heuristics are high-level random search techniques designed to intelligently locate optimal or at least near-optimal solutions for complex optimization problems~\cite{Blum2003}.
As many as $25$ meta-heuristic algorithms were used in this study (see \Cref{tab:EAs}).

In~\cite{Gharehchopogh2019a,Gharehchopogh2020}, the authors concluded that the design and development of meta-heuristic algorithms had been studied more than other optimization techniques due to four main factors: i) meta-heuristic algorithms are inspired by relatively simple concepts from nature that make them easy to implement; ii) these algorithms have many input parameters to control the performance efficiency without altering the algorithm’s structure, making them flexible solving various optimization problems; iii) most meta-heuristic algorithms are derivation-free; iv) meta-heuristic algorithms often escape from local optimum better than other algorithms.

The majority of meta-heuristic algorithms are based on biological evolution principles.
In particular, they are based on simulations of various biological metaphors which differ in the nature of the representation schemes (structure, components, etc.).
The two most common paradigms are evolutionary and swarm systems~\cite{BLUM20114135,BOUSSAID201382,Zavala2014}.

\textbf{Evolutionary algorithms.} In the last decades, there has been a significant interest in evolutionary algorithms, which are based on the principles of natural evolution~\cite{JangaReddy2020}.
Evolutionary Algorithms simulate the biological progression of evolution at the cellular level using selection, crossover, mutation, and reproduction operators to generate increasingly better candidate solutions (chromosomes).
The algorithms consist of a population of individuals, each representing a search point in the optimization domain.
They are exposed to a collective learning process that proceeds from one generation to another.
The initial population is randomly generated and then subjected to selection, crossover, and mutation procedures over several generations to force the newly produced generations to move into more favorable regions of the optimization domain.
The progress in the search is achieved by evaluating the fitness (objective function) of all individuals in the population, selecting individuals with a better fitness value, and combining them to create new individuals with a higher probability of improving the previous fitness.
In the long sequence of the generations, the algorithms converge, and the best individual represents the solution.
Many different evolutionary approaches have been developed, but the basic algorithmic structure is very much the same for all the algorithms, see right side of \Cref{fig:flowchart_direct}.

The most popular evolutionary approach is the genetic algorithm (GA)~\cite{Holland1975} (see \Cref{tab:EAs}).
The GA is a population-based probabilistic search algorithm based on natural selection and genetics mechanics.
A genetic algorithm begins its search with a set of individuals, called population, randomly generated to cover the entire search space uniformly.
Individuals are associated with identity genes that define a fitness measure.
Then, evaluation of the fitness, selection of parents, applying genetic operations, crossover operator for creating offspring, and mutation operation for perturbing the individuals to produce a new population is performed.
The selection operator utilizes the ``survival of the fittest'' concept from Darwinian evolution theory and uses probabilistic rules to select the fittest candidate solutions (best in terms of the objective function) in the current population.
The iterative process is repeated until one of the stopping conditions is satisfied.



\textbf{Swarm intelligence.} Swarm intelligence (or bio-inspired computation) is an integral part of the field of artificial intelligence.
Mainly motivated by biological systems, swarm intelligence adopts the collective behavior of an organized group of animals as they seek to survive.
The main algorithms that fall under swarm intelligence approaches include ant colony optimization (ACO)~\cite{Dorigo1992O}, particle swarm optimization (PSO)~\cite{Eberhart1995}, artificial bee colony (ABC)~\cite{Karaboga2005ANIB}, bees algorithm (BE)~\cite{Pham2006}, invasive weed optimization (IWO)~\cite{Mehrabian2006}, firefly algorithms (FA)~\cite{Yang2009as}, cuckoo search (CS)~\cite{Yang2009a}, bat-inspired algorithm (BAT)~\cite{Yang2010as}.
Similar to evolutionary algorithms, swarm intelligence models are population-based iterative solvers.
The search begins with a randomly initialized population of individuals.
These individuals are then iteratively manipulated and evolved by mimicking the behavior of insects or animals to find the optimum solution.


\subsubsection{Model-based optimization (MBO)}\label{sec:model}

Model-based optimization algorithms generate a population of new points by sampling from a model (or a distribution).
The model (or a distribution) guides structural properties of the underlying real objective function $f$.
Model-based optimization algorithms are based on the concept that the search is directed into regions with improved solutions by adopting the model (or the distribution).
One of the essential ideas in model-based optimization is to substitute the expensive evaluations of the real objective function $f$ with evaluations of a cheap, grained model $\hat{f}$.

Bayesian optimization \cite{Jones1998} is one of the most popular model-based state-of-the-art machine learning frameworks for optimizing expensive black-box functions.
The Bayesian optimization constructs a Gaussian process to approximate the objective function.
A built model is a relatively low-cost surrogate to help guide the search toward promising/unknown regions.
Three MBO technique-based algorithms are considered in this survey: the Bayesian adaptive direct search \cite{Acerbi2017}, surrogate optimization algorithm \cite{Gutmann2001}, and stochastic radial basis Function \cite{Regis2007} (see \Cref{tab:EAs}).

To sum up, the main advantages of stochastic techniques are robust performance on different optimization problems, including high-dimensional ones, and relative simplicity of implementation.
At the same time, as main disadvantages are dependency on numerous input parameters~\cite{zilinskas_zhigljavsky_2016} and the possibility of very slow convergence.

\subsection{Additionally tested solvers}\label{sec:aditional_solvers}

We have additionally investigated more solvers in this study.
Here we briefly mention them and why they were not included in the final comparison:
\begin{enumerate}
	\item \tomlab\texttt{/EGO} -- Implementation of efficient global optimization algorithm for expensive ``black-box'' functions \cite{Jones1998}.
	\item \texttt{\tomlab/rbfSolve} -- Implementation of the radial basis function-based approach \cite{Bjorkman2000} that can handle expensive box-constrained ``black-box'' global optimization problems.
	\item \texttt{\matlab/bayesopt} -- Implement the Bayesian optimization \cite{Bjorkman2000} that can handle expensive box-constrained ``black-box'' global optimization problems.
	\item \texttt{1-DTDV-IA} Implementation of the improved aggressive version of the \adc{} algorithm \cite{Stripinis2021b}.
	\item \texttt{1-DTDV-GL} Implementation of the two-step-based (Global-Local) Pareto selection-based \adc{} algorithm \cite{Stripinis2021b}.
\end{enumerate}

The \texttt{EGO}, \texttt{rbfSolve}, and \texttt{bayesopt} algorithms are designed specifically for expensive objective functions, but the amount of computation involved in algorithmic steps has made these algorithms extremely slow.
Similarly, the \adc{} algorithm's partitioning scheme is considered slow, especially in high dimensions.
The above-mentioned algorithms were not considered in further comparison studies for these reasons.

\section{Experimental design}\label{sec:comparison}

\subsection{Test problems}\label{sec:testproblems}

This study uses two problem suites: \directgolib~\cite{DIRECTGOLib2022} and \gkls{} test problem generator~\cite{Gaviano2003}.
Problem attributes of uni-modality or multi-modality, low or high dimensional, and convex or non-convex geometries are covered.

\subsubsection{\gkls{} test problems}\label{sec:gkls}

The \gkls{} software~\cite{Gaviano2003} generates non-differentiable, continuously differentiable, and twice continuously differentiable classes of test functions for multi-modal, multi-dimensional box-constrained global optimization.
For each generated problem, all local and global minima are known.
The test problems are constructed by defining a convex quadratic function (paraboloid) systematically distorted by polynomials to produce local minima.
Each class of problems consists of 100 test functions.
The complexity of the class is determined using the following parameters: problem dimension $(n)$, number of local minima $(h)$, the value of the global minimum $(f^*)$, radius $(r)$ of the attraction region of the global minimizer and the distance $(d)$ from the global minimizer to the vertex of the quadratic function.
The complete repeatability of experiments is an essential feature of the generator.
If the same five parameters are provided to \gkls{}, the identical class of functions will be produced each time the generator is executed.

We use standard eight different complexity classes (see \Cref{tab:paramet}), which are the most widely and commonly used in other numerical studies \cite{Sergeyev2006,Paulavicius2014:jogo,Paulavicius2019:eswa,Stripinis2021b}.
The dimension ($n$) of generated test classes and other parameters are the same as in other mentioned studies.
For each dimension $n$, two test classes were considered: the ``simple'' class and the ``hard'' one.
For third and fourth dimensional classes, the difficulty is increased by enlarging the distance $d$ from the global optimum point $(\mathbf{x}^*)$ to the paraboloid vertex.
For second and fifth dimensional classes, this is achieved by decreasing the radius $r$.

\begin{table}[h]
	\centering
	\normalsize
	\caption{\gkls{} test classes used in numerical experiments.}
	\label{tab:paramet}
	\begin{tabular*}{\textwidth}{@{\extracolsep{\fill}}cccccccc}
		\toprule
		Class & $\#$ & Difficulty & $n$ & $f^*$ & $d$ & $r$ & $h$ \\
		\midrule
		$1$ & $100$ & simple 	 & $2$ & $-1$ & $0.90$ & $0.20$ & $10$  \\
		$2$ & $100$ & hard 		 & $2$ & $-1$ & $0.90$ & $0.10$ & $10$  \\
		$3$ & $100$ & simple 	 & $3$ & $-1$ & $0.66$ & $0.20$ & $10$  \\
		$4$ & $100$ & hard 		 & $3$ & $-1$ & $0.90$ & $0.20$ & $10$  \\
		$5$ & $100$ & simple 	 & $4$ & $-1$ & $0.66$ & $0.20$ & $10$  \\
		$6$ & $100$ & hard 		 & $4$ & $-1$ & $0.90$ & $0.20$ & $10$  \\
		$7$ & $100$ & simple     & $5$ & $-1$ & $0.66$ & $0.30$ & $10$  \\
		$8$ & $100$ & hard 		 & $5$ & $-1$ & $0.66$ & $0.20$ & $10$  \\
		\bottomrule
	\end{tabular*}
\end{table}

\subsubsection{\directgolib{} test set}\label{sec:directlib}

The second set of test problems chosen for the experiments comprises more traditional cases collected from the literature.
Test problems from the most recent and significantly expanded ($69$ new test problems, including 18 published this year \cite{Layeb2022}) \directgolib{}~\cite{DIRECTGOLib2022} library are used to evaluate considered algorithms.

Overview of all employed box-constrained optimization test problems from \directgolib{} and their properties are given in Appendix~\ref{apendix:figures} \Cref{tab:tests,tab:test}.
Here, the main features are reported: problem number (\#), name of the problem, source, dimensionality ($ n $), optimization domain ($ D $), problem type, and the known minimum ($ f^* $).
Moreover, the original domains $D$ for some problems are perturbed ($\tilde{D}$) that the solutions are not located in their midpoints or other points favorable for any deterministic algorithms.
Finally, newly added test problems are marked with the symbol $\beta$. 

Some of these test problems have several variants, e.g., \textit{Ackley}, \textit{Hartman}, \textit{Shekel}, and some of them, like \textit{Alpine}, \textit{Csendes}, \textit{Griewank}, etc., can be tested for varying dimensionality.
The $67$ test problems listed in \Cref{tab:tests} have fixed dimensions, while \Cref{tab:test} presents $55$ test benchmarks that can be used by specifying any dimension size ($n$).
For these test problems instances with $ n = 2, 5, 10, 20, 50,$ and $100$ variables were considered, leading in total to $397$ test problems (see summary in~\Cref{tab:directlib}).

Various subsets (e.g., non-convex, multi-modal, etc.) of the whole set were used to deepen the investigation.
The low-dimensional ($n \leq 10$) test set is designed to evaluate the performance of solvers on relatively small dimensional problems.
The high-dimensional ($n \geq 11$) test set is used to test the efficiency of solvers on higher dimensions.

\begin{table}[ht]
	\centering
	\caption{Characteristics of \directgolib{} test problems.}
		\begin{tabular*}{\textwidth}{@{\extracolsep{\fill}}lcccc}
				\toprule
				Dimension / type & convex & non-convex & uni-modal & multi-modal \\
				\midrule
				$n \leq 10$  & $55$ & $176$ & $41$ & $190$ \\
				$n \geq 11$ & $42$ & $124$ & $36$ & $130$ \\
				\midrule
				Total 						   & $97$ & $300$ & $77$ & $320$ \\
				\bottomrule
			\end{tabular*}
		\label{tab:directlib}
	\end{table}
	
	\subsection{Experimental setup and basis of solver comparisons}\label{sec:stop}
	
	All computations were performed on Intel R Core$^\textit{TM}$ i5-10400 @ 2.90GHz Processor running \matlab{} R2022a.
	The solutions returned by the solvers against the globally optimal solution for each problem were compared.
	A solver was considered to have successfully solved the test problem during a run if it returned a solution with an objective function value within $1 \%$.
	For all analytical test cases with a priori known global optima $ f^* $, the used stopping criterion is based on the percent error $(pe)$:
	\begin{equation}
		\label{eq:pe}
		\ pe = 100 \% \times
		\begin{cases}
			\frac{f({\mathbf{x}}) - f^*}{\mid f^* \mid},  & f^* \neq 0, \\
			f({\mathbf{x}}),  & f^* = 0,
		\end{cases}
	\end{equation}
	where $ f^* $ is the known global optimum.
	The algorithms were stopped when the percent error became smaller than the prescribed value $\varepsilon_{\rm pe} = 1$ or when the number of function evaluations exceeded the prescribed limit ($M_{\rm max}$).
	We note that when the optimal value is large, the algorithm will terminate when the distance to the optimum is relatively large.
	However, there were only a few such tasks in the test set.
	
	We used two different values of $M_{\rm max}$.
	First, like in \cite{Rios2013}, we use $2,500$ to evaluate which algorithms perform best for expensive objective problems.
	Additionally, $M_{\rm max} = 500,000$ is used to evaluate the performance when the objective function evaluations are cheap.
	As all algorithms were implemented in the same environment (\matlab),  additionally we included the execution time in the comparison.
	A limit ($T_{\rm max}$) of $3,600$ CPU seconds was imposed on each run.
	
	\subsection{Algorithmic settings}\label{sec:settings}
	
	The number of chosen algorithms is extensive, they are very different, and various algorithmic control structures determine their effectiveness. Algorithms also have different input parameters that can have a significant impact their performance.
	For example, \direct-type algorithms have a much smaller number of input parameters than meta-heuristic solvers, some of whom have none at all (e.g., \directaggress, \directg{} and \directgl{} from \dgo{} toolbox~\cite{Stripinis2021c,DIRECTGO2022}).
	The input parameters should be set in such a way as to ensure the best performance.
	Unfortunately, different input parameter values can drastically impact the algorithm's performance.
	Thus, following the same idea used in the \cite{Rios2013,Stripinis2021c} study, comparisons were carried out using the default parameter values for each algorithm.
	
	Among the $64$ algorithms involved in the study, $11$ are hybridized with local search procedures.
	All of them are distinguished by adding the * symbol at the end of the title.
	Most of these algorithms (\directrev, \dirmin, \birmin, \multistart, and \globalsearch) are hybridized using the interior-point \cite{Byrd2000} algorithm, implemented in \matlab{} nonlinear programming solver -- \texttt{fmincon}.
	\tomlab's global optimization algorithms hybridize with \texttt{SNPOT} or \texttt{NPSOL} \cite{Holmstrom2010} local minimizers.
	No derivative information was provided for them.
	Thus, if needed, the local solvers could only use the function values to approximate gradients, e.g., by a finite difference technique.
	However, the use of finite differences has been broadly dismissed in the derivative-free literature as expensive in terms of function evaluations~\cite{Shi2021}.
	
	\subsection{Benchmarking derivative-free optimization solvers}\label{sec:dataprofiles}
	
	To analyze and compare the algorithms' performance, we applied the data profiles~\cite{More2009} to the convergence test~\eqref{eq:pe}.
	The data profile is a popular and widely used tool for bench-marking and evaluating the performance of several algorithms (solvers) when run on a large problem set.
	Benchmark results are generated by running a certain algorithm $v$ (from a set of algorithms $\mathcal{V}$ under consideration) for each problem $u$ from a benchmark set $\mathcal{U}$ and recording the performance measure of interest.
	The performance measure could be, for example, the number of function evaluations, the computation time, the number of iterations, or the memory used.
	We used a number of function evaluations and the execution (computation) time criteria.
	
	The data profiles provide the percentage of problems that can be solved with a given budget of the desired performance measure.
	The data profile is defined
	
	\begin{equation}
		\lambda_{v}(\alpha) = \frac{1}{card(\mathcal{U})}\textrm{size} \left\{ u \in \mathcal{U} : t_{u,v} \le \alpha \right\},
		\label{eq:data-profile}
	\end{equation}
	where $t_{u,v} > 0$ is the number of performance measure required to solve problem $u$ by the algorithm $v$, and  $card(\mathcal{U})$ is the cardinality of $\mathcal{U}$.
	In our case, the $\lambda_{v}(\alpha)$ shows the percentage of problems that can be solved within $\alpha$ function evaluations, or seconds.

	\section{Numerical results}
	\label{sec:benchmarking}
	
	\subsection{Computational results solving \gkls{} test problems}\label{sec:gklstest}
	
	Numerical results on eight \gkls{} test classes from \Cref{tab:paramet} are reported in \Cref{tab:results0,tab:results1} for ``simple'' and ``hard''  classes of \gkls{} test problems using two different $M_{\rm max}$ values.
	The best results are given in bold.
	In both tables, the algorithms are ranked based on the average success rate (S.R.) reported in the sixth and eleventh columns.
	
	\begin{table}
		\caption{The average number ($f_{\rm avg.}$) and standard deviation ($f_{\rm s.d.}$) of function evaluations, the average value ($t_{\rm avg.}$) and standard deviation ($t_{\rm s.d.}$) of execution time, and the success rate (S. R.) for ``simple'' \gkls{} test problems using $M_{\rm max}=2,500$ and $M_{\rm max}=500,000$.}
		\resizebox{1\textwidth}{!}{
			\begin{tabular}[tb]{@{\extracolsep{\fill}}l|rrrrr|rrrrr}
				\toprule
				Evaluation budget & \multicolumn{5}{c|}{$M_{\rm max}=2,500$} & \multicolumn{5}{c}{$M_{\rm max}=500,000$} \\
				\midrule
				Criteria & \multicolumn{1}{c}{$f_{\rm avg.}$} & \multicolumn{1}{c}{$f_{\rm s.d.}$} & \multicolumn{1}{c}{$t_{\rm avg.}$} & \multicolumn{1}{c}{$t_{\rm s.d.}$} & \multicolumn{1}{c|}{S. R.} & \multicolumn{1}{c}{$f_{\rm avg.}$} & \multicolumn{1}{c}{$f_{\rm s.d.}$} & \multicolumn{1}{c}{$t_{\rm avg.}$} & \multicolumn{1}{c}{$t_{\rm s.d.}$} & \multicolumn{1}{c}{S. R.} \\
				\midrule
				\birmin 	& $	 \mathbf{1,163} 	$ & $	 - 	$ & $	0.21	$ & $	 - 	$ & $	 \mathbf{0.73} 	$ & $	 \mathbf{1,991} 	$ & $	 - 	$ & $	0.40	$ & $	 - 	$ & $	 \mathbf{1.00} 	$	\\
				\birectgb 	& $	1,235	$ & $	 - 	$ & $	0.20	$ & $	 - 	$ & $	0.69	$ & $	2,243	$ & $	 - 	$ & $	0.44	$ & $	 - 	$ & $	 \mathbf{1.00} 	$	\\
				\mqnlp 	& $	1,663	$ & $	 - 	$ & $	0.30	$ & $	 - 	$ & $	0.69	$ & $	7,361	$ & $	 - 	$ & $	1.84	$ & $	 - 	$ & $	0.99	$	\\
				\directa 	& $	1,219	$ & $	 - 	$ & $	0.11	$ & $	 - 	$ & $	0.67	$ & $	4,154	$ & $	 - 	$ & $	0.34	$ & $	 - 	$ & $	 \mathbf{1.00} 	$	\\
				\directm 	& $	1,233	$ & $	 - 	$ & $	0.12	$ & $	 - 	$ & $	0.67	$ & $	5,080	$ & $	 - 	$ & $	0.42	$ & $	 - 	$ & $	 \mathbf{1.00} 	$	\\
				\directmro 	& $	1,288	$ & $	 - 	$ & $	0.24	$ & $	 - 	$ & $	0.67	$ & $	7,785	$ & $	 - 	$ & $	1.42	$ & $	 - 	$ & $	 \mathbf{1.00} 	$	\\
				\directg 	& $	1,305	$ & $	 - 	$ & $	0.07	$ & $	 - 	$ & $	0.67	$ & $	5,587	$ & $	 - 	$ & $	0.27	$ & $	 - 	$ & $	 \mathbf{1.00} 	$	\\
				\direct 	& $	1,231	$ & $	 - 	$ & $	0.18	$ & $	 - 	$ & $	0.66	$ & $	4,284	$ & $	 - 	$ & $	0.37	$ & $	 - 	$ & $	 \mathbf{1.00} 	$	\\
				\glbsolve 	& $	1,233	$ & $	 - 	$ & $	0.13	$ & $	 - 	$ & $	0.66	$ & $	4,279	$ & $	 - 	$ & $	0.59	$ & $	 - 	$ & $	 \mathbf{1.00} 	$	\\
				\directrev 	& $	1,288	$ & $	 - 	$ & $	0.19	$ & $	 - 	$ & $	0.66	$ & $	3,839	$ & $	 - 	$ & $	0.54	$ & $	 - 	$ & $	 \mathbf{1.00} 	$	\\
				\birect 	& $	1,312	$ & $	 - 	$ & $	0.20	$ & $	 - 	$ & $	0.64	$ & $	4,139	$ & $	 - 	$ & $	0.59	$ & $	 - 	$ & $	 \mathbf{1.00} 	$	\\
				\dbdvo 	& $	1,312	$ & $	 - 	$ & $	0.21	$ & $	 - 	$ & $	0.64	$ & $	4,139	$ & $	 - 	$ & $	0.62	$ & $	 - 	$ & $	 \mathbf{1.00} 	$	\\
				\dbdvg 	& $	1,459	$ & $	 - 	$ & $	0.56	$ & $	 - 	$ & $	0.64	$ & $	6,018	$ & $	 - 	$ & $	0.79	$ & $	 - 	$ & $	 \mathbf{1.00} 	$	\\
				\dtco 	& $	1,314	$ & $	 - 	$ & $	0.19	$ & $	 - 	$ & $	0.63	$ & $	4,145	$ & $	 - 	$ & $	0.72	$ & $	 - 	$ & $	 \mathbf{1.00} 	$	\\
				\directrest 	& $	1,197	$ & $	 - 	$ & $	4.65	$ & $	 - 	$ & $	0.70	$ & $	46,770	$ & $	 - 	$ & $	163.30	$ & $	 - 	$ & $	0.92	$	\\
				\directgl 	& $	1,465	$ & $	 - 	$ & $	0.09	$ & $	 - 	$ & $	0.62	$ & $	7,010	$ & $	 - 	$ & $	0.36	$ & $	 - 	$ & $	 \mathbf{1.00} 	$	\\
				\dtcg 	& $	1,391	$ & $	 - 	$ & $	0.08	$ & $	 - 	$ & $	0.61	$ & $	8,879	$ & $	 - 	$ & $	0.32	$ & $	 - 	$ & $	 \mathbf{1.00} 	$	\\
				\dbdva 	& $	1,537	$ & $	 - 	$ & $	1.65	$ & $	 - 	$ & $	0.58	$ & $	10,543	$ & $	 - 	$ & $	1.41	$ & $	 - 	$ & $	 \mathbf{1.00} 	$	\\
				\glccluster 	& $	1,620	$ & $	 - 	$ & $	0.76	$ & $	 - 	$ & $	0.57	$ & $	4,704	$ & $	 - 	$ & $	2.14	$ & $	 - 	$ & $	 \mathbf{1.00} 	$	\\
				\oqnlp 	& $	1,532	$ & $	 - 	$ & $	0.44	$ & $	 - 	$ & $	0.57	$ & $	32,194	$ & $	 - 	$ & $	10.88	$ & $	 - 	$ & $	0.99	$	\\
				\mcs 	& $	1,601	$ & $	 - 	$ & $	1.15	$ & $	 - 	$ & $	0.55	$ & $	7,879	$ & $	 - 	$ & $	4.73	$ & $	 - 	$ & $	 \mathbf{1.00} 	$	\\
				\adc 	& $	1,565	$ & $	 - 	$ & $	5.19	$ & $	 - 	$ & $	0.54	$ & $	4,429	$ & $	 - 	$ & $	12.16	$ & $	 - 	$ & $	 \mathbf{1.00} 	$	\\
				\dtca 	& $	1,524	$ & $	 - 	$ & $	0.06	$ & $	 - 	$ & $	0.54	$ & $	16,281	$ & $	 - 	$ & $	1.27	$ & $	 - 	$ & $	 \mathbf{1.00} 	$	\\
				\directl 	& $	1,440	$ & $	 - 	$ & $	1.67	$ & $	 - 	$ & $	0.55	$ & $	27,988	$ & $	 - 	$ & $	11.44	$ & $	 - 	$ & $	0.98	$	\\
				\disimplv 	& $	1,528	$ & $	 - 	$ & $	3.69	$ & $	 - 	$ & $	0.50	$ & $	18,793	$ & $	 - 	$ & $	103.62	$ & $	 - 	$ & $	0.98	$	\\
				\disimplc 	& $	1,664	$ & $	 - 	$ & $	5.31	$ & $	 - 	$ & $	0.45	$ & $	49,439	$ & $	 - 	$ & $	488.99	$ & $	 - 	$ & $	0.97	$	\\
				\directaggress 	& $	1,809	$ & $	 - 	$ & $	0.09	$ & $	 - 	$ & $	0.43	$ & $	34,903	$ & $	 - 	$ & $	1.37	$ & $	 - 	$ & $	0.98	$	\\
				\dirmin 	& $	1,901	$ & $	 - 	$ & $	0.32	$ & $	 - 	$ & $	0.41	$ & $	23,477	$ & $	 - 	$ & $	3.72	$ & $	 - 	$ & $	 \mathbf{1.00} 	$	\\
				\plor 	& $	1,577	$ & $	 - 	$ & $	4.16	$ & $	 - 	$ & $	0.47	$ & $	88,448	$ & $	 - 	$ & $	293.94	$ & $	 - 	$ & $	0.88	$	\\
				\sa 	& $	1,972	$ & $	432	$ & $	1.24	$ & $	4.41	$ & $	0.36	$ & $	21,497	$ & $	16,336	$ & $	21.14	$ & $	82.60	$ & $	0.99	$	\\
				\multistart 	& $	2,000	$ & $	410	$ & $	3.49	$ & $	4.18	$ & $	0.36	$ & $	21,515	$ & $	16,827	$ & $	16.32	$ & $	26.75	$ & $	0.99	$	\\
				\StochasticRBF	& $	2,010	$ & $	318	$ & $	2.18	$ & $	7.54	$ & $	0.31	$ & $	34,851	$ & $	21,410	$ & $	147.17	$ & $	259.58	$ & $	0.99	$	\\
				\scso 	& $	1,825	$ & $	306	$ & $	1.33	$ & $	5.50	$ & $	0.41	$ & $	80,896	$ & $	53,683	$ & $	10.02	$ & $	13.71	$ & $	0.86	$	\\
				\surogate 	& $	1,486	$ & $	251	$ & $	72.91	$ & $	55.00	$ & $	0.56	$ & $	159,431	$ & $	64,750	$ & $	1,156.32	$ & $	1,620.79	$ & $	0.68	$	\\
				\be 	& $	1,978	$ & $	278	$ & $	0.51	$ & $	0.18	$ & $	0.42	$ & $	123,802	$ & $	105,942	$ & $	16.33	$ & $	28.34	$ & $	0.76	$	\\
				\sta 	& $	2,409	$ & $	78	$ & $	1.34	$ & $	5.93	$ & $	0.31	$ & $	89,792	$ & $	71,245	$ & $	10.02	$ & $	13.71	$ & $	0.84	$	\\
				\avoa 	& $	1,984	$ & $	340	$ & $	0.03	$ & $	0.01	$ & $	0.34	$ & $	125,968	$ & $	70,637	$ & $	1.73	$ & $	2.72	$ & $	0.80	$	\\
				\lgo 	& $	2,267	$ & $	 - 	$ & $	1.25	$ & $	 - 	$ & $	0.14	$ & $	135,346	$ & $	 - 	$ & $	15.87	$ & $	 - 	$ & $	0.98	$	\\
				\multimin 	& $	2,470	$ & $	17	$ & $	0.33	$ & $	0.17	$ & $	0.13	$ & $	25,702	$ & $	4,193	$ & $	3.26	$ & $	7.51	$ & $	0.98	$	\\
				\aos 	& $	2,026	$ & $	373	$ & $	0.13	$ & $	0.05	$ & $	0.32	$ & $	162,610	$ & $	69,751	$ & $	10.02	$ & $	13.71	$ & $	0.76	$	\\
				\hs 	& $	2,112	$ & $	331	$ & $	0.31	$ & $	0.05	$ & $	0.32	$ & $	171,004	$ & $	79,453	$ & $	22.66	$ & $	29.10	$ & $	0.72	$	\\
				\amc 	& $	2,015	$ & $	147	$ & $	0.11	$ & $	0.04	$ & $	0.45	$ & $	205,638	$ & $	100,081	$ & $	10.30	$ & $	12.26	$ & $	0.59	$	\\
				\globalsearch 	& $	2,500	$ & $	0	$ & $	0.52	$ & $	0.00	$ & $	0.00	$ & $	29,892	$ & $	19,291	$ & $	6.48	$ & $	10.03	$ & $	0.99	$	\\
				\nna	& $	2,292	$ & $	150	$ & $	0.08	$ & $	0.01	$ & $	0.22	$ & $	140,758	$ & $	71,521	$ & $	1.29	$ & $	1.79	$ & $	0.75	$	\\
				\cnna 	& $	2,313	$ & $	151	$ & $	0.07	$ & $	0.01	$ & $	0.15	$ & $	116,336	$ & $	62,272	$ & $	2.56	$ & $	4.30	$ & $	0.81	$	\\
				\pa 	& $	1,899	$ & $	364	$ & $	0.26	$ & $	0.11	$ & $	0.46	$ & $	259,194	$ & $	135,680	$ & $	35.71	$ & $	34.56	$ & $	0.48	$	\\
				\ica 	& $	1,926	$ & $	395	$ & $	0.63	$ & $	0.33	$ & $	0.39	$ & $	257,268	$ & $	171,481	$ & $	77.32	$ & $	76.06	$ & $	0.49	$	\\
				\abc 	& $	2,382	$ & $	123	$ & $	0.34	$ & $	0.05	$ & $	0.12	$ & $	178,138	$ & $	58,343	$ & $	24.97	$ & $	30.87	$ & $	0.73	$	\\
				\cs 	& $	2,373	$ & $	127	$ & $	0.03	$ & $	0.01	$ & $	0.15	$ & $	154,038	$ & $	139,577	$ & $	1.84	$ & $	2.72	$ & $	0.69	$	\\
				\beso 	& $	2,230	$ & $	192	$ & $	0.07	$ & $	0.02	$ & $	0.24	$ & $	254,856	$ & $	112,329	$ & $	7.62	$ & $	7.41	$ & $	0.49	$	\\
				\sfla 	& $	1,861	$ & $	520	$ & $	1.22	$ & $	5.21	$ & $	0.36	$ & $	320,261	$ & $	145,928	$ & $	16.34	$ & $	12.26	$ & $	0.36	$	\\
				\tlbo 	& $	2,215	$ & $	260	$ & $	0.32	$ & $	0.09	$ & $	0.25	$ & $	287,069	$ & $	126,127	$ & $	39.51	$ & $	33.79	$ & $	0.44	$	\\
				\ga 	& $	2,055	$ & $	408	$ & $	4.06	$ & $	9.33	$ & $	0.31	$ & $	341,106	$ & $	146,438	$ & $	25.41	$ & $	17.82	$ & $	0.32	$	\\
				\bads	& $	2,309	$ & $	222	$ & $	22.63	$ & $	21.33	$ & $	0.12	$ & $	269,637	$ & $	186,611	$ & $	1,955.32	$ & $	1,372.54$ & $	0.47	$	\\
				\fstde 	& $	2,167	$ & $	319	$ & $	0.33	$ & $	3.12	$ & $	0.24	$ & $	327,605	$ & $	159,995	$ & $	3.66	$ & $	12.63	$ & $	0.35	$	\\
				\bbo 	& $	2,113	$ & $	464	$ & $	0.34	$ & $	0.12	$ & $	0.22	$ & $	392,402	$ & $	125,295	$ & $	63.88	$ & $	34.24	$ & $	0.22	$	\\
				\iwo 	& $	2,096	$ & $	380	$ & $	0.23	$ & $	0.09	$ & $	0.21	$ & $	393,382	$ & $	109,053	$ & $	43.65	$ & $	22.97	$ & $	0.21	$	\\
				\fa 	& $	2,075	$ & $	510	$ & $	0.37	$ & $	0.15	$ & $	0.20	$ & $	398,582	$ & $	123,555	$ & $	64.66	$ & $	33.29	$ & $	0.20	$	\\
				\aco 	& $	2,267	$ & $	285	$ & $	0.32	$ & $	0.08	$ & $	0.14	$ & $	405,410	$ & $	113,734	$ & $	65.94	$ & $	33.39	$ & $	0.19	$	\\
				\mc 	& $	2,465	$ & $	62	$ & $	0.10	$ & $	0.01	$ & $	0.03	$ & $	365,777	$ & $	23,752	$ & $	15.04	$ & $	8.78	$ & $	0.30	$	\\
				\de 	& $	2,171	$ & $	415	$ & $	0.50	$ & $	0.18	$ & $	0.15	$ & $	422,807	$ & $	98,451	$ & $	98.60	$ & $	43.70	$ & $	0.15	$	\\
				\bat 	& $	2,304	$ & $	267	$ & $	0.04	$ & $	0.01	$ & $	0.09	$ & $	452,293	$ & $	66,491	$ & $	7.77	$ & $	2.82	$ & $	0.09	$	\\
				\cmaes	& $	2,461	$ & $	61	$ & $	0.03	$ & $	0.00	$ & $	0.03	$ & $	495,015	$ & $	23,339	$ & $	6.28	$ & $	2.11	$ & $	0.03	$	\\
				\ca 	& $	2,466	$ & $	36	$ & $	0.32	$ & $	0.03	$ & $	0.02	$ & $	487,905	$ & $	12,731	$ & $	62.01	$ & $	10.49	$ & $	0.02	$	\\
				
				\bottomrule
		\end{tabular}}
		\label{tab:results0}
	\end{table}
	
	\subsubsection{Comparison on ``simple'' \gkls{} problems with $M_{\rm max} = 2,500$}
	
	The second, third, and sixth columns in \Cref{tab:results0} show the average number of function evaluations ($f_{\rm avg.}$), standard deviation ($f_{\rm s.d.}$), and success rate (S. R.) for all $400$ ``simple'' \gkls{}-type problems using $M_{\rm max} = 2,500$.
	The hybridized \birmin{} algorithm showed the best average number of function evaluations, solving $292$ problems within a budget of $2,500$ function evaluations.
	The best \birmin{} algorithm is closely followed by the \directrest{}, which solved $280$, while \birectgb{} and \mqnlp{} solved $276$ test problems.
	
	\begin{table}
		\caption{The average number ($f_{\rm avg.}$) and standard deviation ($f_{\rm s.d.}$) of function evaluations, the average value ($t_{\rm avg.}$) and standard deviation ($t_{\rm s.d.}$) of execution time, and the success rate (S. R.) for ``hard'' \gkls{} test problems using $M_{\rm max}=2,500$ and $M_{\rm max}=500,000$.}
		\resizebox{1\textwidth}{!}{
			\begin{tabular}[tb]{@{\extracolsep{\fill}}l|rrrrr|rrrrr}
				\toprule
				Evaluation budget & \multicolumn{5}{c|}{$M_{\rm max}=2,500$} & \multicolumn{5}{c}{$M_{\rm max}=500,000$} \\
				\midrule
				Criteria & \multicolumn{1}{c}{$f_{\rm avg.}$} & \multicolumn{1}{c}{$f_{\rm s.d.}$} & \multicolumn{1}{c}{$t_{\rm avg.}$} & \multicolumn{1}{c}{$t_{\rm s.d.}$} & \multicolumn{1}{c|}{S. R.} & \multicolumn{1}{c}{$f_{\rm avg.}$} & \multicolumn{1}{c}{$f_{\rm s.d.}$} & \multicolumn{1}{c}{$t_{\rm avg.}$} & \multicolumn{1}{c}{$t_{\rm s.d.}$} & \multicolumn{1}{c}{S. R.} \\
				\midrule
				\birectgb 	& $	 \mathbf{1,594} 	$ & $	 - 	$ & $	0.28	$ & $	 - 	$ & $	 \mathbf{0.54} 	$ & $	7,482	$ & $	 - 	$ & $	1.91	$ & $	 - 	$ & $	 \mathbf{1.00} 	$	\\
				\birmin 	& $	1,638	$ & $	 - 	$ & $	0.31	$ & $	 - 	$ & $	 \mathbf{0.54} 	$ & $	 \mathbf{7,384} 	$ & $	 - 	$ & $	1.94	$ & $	 - 	$ & $	 \mathbf{1.00} 	$	\\
				\birect 	& $	1,795	$ & $	 - 	$ & $	0.28	$ & $	 - 	$ & $	0.46	$ & $	18,735	$ & $	 - 	$ & $	4.22	$ & $	 - 	$ & $	 \mathbf{1.00} 	$	\\
				\dbdvo 	& $	1,795	$ & $	 - 	$ & $	0.28	$ & $	 - 	$ & $	0.46	$ & $	18,735	$ & $	 - 	$ & $	3.97	$ & $	 - 	$ & $	 \mathbf{1.00} 	$	\\
				\dtco 	& $	1,876	$ & $	 - 	$ & $	0.17	$ & $	 - 	$ & $	0.44	$ & $	21,633	$ & $	 - 	$ & $	3.17	$ & $	 - 	$ & $	 \mathbf{1.00} 	$	\\
				\directa 	& $	1,838	$ & $	 - 	$ & $	0.18	$ & $	 - 	$ & $	0.44	$ & $	24,121	$ & $	 - 	$ & $	2.51	$ & $	 - 	$ & $	0.99	$	\\
				\glbsolve 	& $	1,847	$ & $	 - 	$ & $	0.21	$ & $	 - 	$ & $	0.44	$ & $	25,604	$ & $	 - 	$ & $	10.05	$ & $	 - 	$ & $	0.99	$	\\
				\direct 	& $	1,847	$ & $	 - 	$ & $	0.28	$ & $	 - 	$ & $	0.44	$ & $	31,770	$ & $	 - 	$ & $	73.64	$ & $	 - 	$ & $	0.99	$	\\
				\directrev 	& $	1,907	$ & $	 - 	$ & $	0.30	$ & $	 - 	$ & $	0.42	$ & $	31,176	$ & $	 - 	$ & $	5.50	$ & $	 - 	$ & $	0.99	$	\\
				\dbdvg 	& $	1,880	$ & $	 - 	$ & $	1.09	$ & $	 - 	$ & $	0.42	$ & $	52,264	$ & $	 - 	$ & $	3.18	$ & $	 - 	$ & $	0.98	$	\\
				\directm 	& $	1,892	$ & $	 - 	$ & $	0.20	$ & $	 - 	$ & $	0.40	$ & $	32,943	$ & $	 - 	$ & $	3.62	$ & $	 - 	$ & $	0.99	$	\\
				\glccluster 	& $	2,071	$ & $	 - 	$ & $	1.03	$ & $	 - 	$ & $	0.38	$ & $	23,717	$ & $	 - 	$ & $	102.25	$ & $	 - 	$ & $	0.99	$	\\
				\adc 	& $	2,011	$ & $	 - 	$ & $	6.87	$ & $	 - 	$ & $	0.36	$ & $	17,156	$ & $	 - 	$ & $	302.90	$ & $	 - 	$ & $	 \mathbf{1.00} 	$	\\
				\directg 	& $	2,007	$ & $	 - 	$ & $	0.13	$ & $	 - 	$ & $	0.32	$ & $	50,997	$ & $	 - 	$ & $	3.14	$ & $	 - 	$ & $	0.98	$	\\
				\mcs 	& $	2,065	$ & $	 - 	$ & $	1.61	$ & $	 - 	$ & $	0.32	$ & $	48,359	$ & $	 - 	$ & $	400.72	$ & $	 - 	$ & $	0.98	$	\\
				\dbdva 	& $	2,073	$ & $	 - 	$ & $	2.89	$ & $	 - 	$ & $	0.32	$ & $	59,344	$ & $	 - 	$ & $	6.25	$ & $	 - 	$ & $	0.98	$	\\
				\directgl 	& $	1,995	$ & $	 - 	$ & $	0.14	$ & $	 - 	$ & $	0.35	$ & $	70,499	$ & $	 - 	$ & $	4.50	$ & $	 - 	$ & $	0.95	$	\\
				\directl 	& $	1,961	$ & $	 - 	$ & $	2.25	$ & $	 - 	$ & $	0.40	$ & $	99,258	$ & $	 - 	$ & $	70.67	$ & $	 - 	$ & $	0.88	$	\\
				\dtcg 	& $	2,017	$ & $	 - 	$ & $	0.14	$ & $	 - 	$ & $	0.30	$ & $	59,837	$ & $	 - 	$ & $	3.25	$ & $	 - 	$ & $	0.97	$	\\
				\directmro 	& $	2,066	$ & $	 - 	$ & $	0.43	$ & $	 - 	$ & $	0.28	$ & $	67,746	$ & $	 - 	$ & $	19.93	$ & $	 - 	$ & $	0.97	$	\\
				\mqnlp 	& $	2,176	$ & $	 - 	$ & $	0.49	$ & $	 - 	$ & $	0.31	$ & $	73,834	$ & $	 - 	$ & $	18.49	$ & $	 - 	$ & $	0.89	$	\\
				\disimplv 	& $	1,973	$ & $	 - 	$ & $	4.69	$ & $	 - 	$ & $	0.37	$ & $	102,617	$ & $	 - 	$ & $	818.48	$ & $	 - 	$ & $	0.82	$	\\
				\plor 	& $	1,847	$ & $	 - 	$ & $	4.78	$ & $	 - 	$ & $	0.41	$ & $	138,981	$ & $	 - 	$ & $	570.12	$ & $	 - 	$ & $	0.78	$	\\
				\oqnlp 	& $	2,039	$ & $	 - 	$ & $	0.64	$ & $	 - 	$ & $	0.31	$ & $	101,334	$ & $	 - 	$ & $	34.49	$ & $	 - 	$ & $	0.87	$	\\
				\directrest 	& $	1,763	$ & $	 - 	$ & $	7.33	$ & $	 - 	$ & $	0.45	$ & $	150,320	$ & $	 - 	$ & $	435.87	$ & $	 - 	$ & $	0.72	$	\\
				\sa 	& $	2,265	$ & $	307	$ & $	2.17	$ & $	19.58	$ & $	0.18	$ & $	58,740	$ & $	42,115	$ & $	110.45	$ & $	267.37	$ & $	0.98	$	\\
				\dtca 	& $	2,160	$ & $	 - 	$ & $	0.11	$ & $	 - 	$ & $	0.22	$ & $	87,158	$ & $	 - 	$ & $	5.28	$ & $	 - 	$ & $	0.94	$	\\
				\disimplc 	& $	2,068	$ & $	 - 	$ & $	6.72	$ & $	 - 	$ & $	0.30	$ & $	129,657	$ & $	 - 	$ & $	989.05	$ & $	 - 	$ & $	0.81	$	\\
				\dirmin 	& $	2,286	$ & $	 - 	$ & $	0.42	$ & $	 - 	$ & $	0.15	$ & $	95,598	$ & $	 - 	$ & $	14.39	$ & $	 - 	$ & $	0.94	$	\\
				\multistart 	& $	2,320	$ & $	249	$ & $	4.13	$ & $	6.94	$ & $	0.15	$ & $	89,525	$ & $	51,552	$ & $	74.58	$ & $	115.32	$ & $	0.94	$	\\
				\StochasticRBF	& $	2,321	$ & $	171	$ & $	2.04	$ & $	0.82	$ & $	0.13	$ & $	141,541	$ & $	73,122	$ & $	765.35	$ & $	1,073.29	$ & $	0.92	$	\\
				\lgo 	& $	2,406	$ & $	 - 	$ & $	1.34	$ & $	 - 	$ & $	0.06	$ & $	224,132	$ & $	 - 	$ & $	24.19	$ & $	 - 	$ & $	0.98	$	\\
				\scso 	& $	2,174	$ & $	234	$ & $	0.60	$ & $	3.04	$ & $	0.23	$ & $	152,263	$ & $	60,166	$ & $	13.75	$ & $	14.66	$ & $	0.75	$	\\
				\multimin 	& $	2,475	$ & $	21	$ & $	0.39	$ & $	0.13	$ & $	0.09	$ & $	90,724	$ & $	4,012	$ & $	13.51	$ & $	24.20	$ & $	0.88	$	\\
				\globalsearch 	& $	2,500	$ & $	0	$ & $	0.51	$ & $	0.00	$ & $	0.00	$ & $	77,313	$ & $	47,342	$ & $	17.26	$ & $	25.92	$ & $	0.96	$	\\
				\directaggress 	& $	2,268	$ & $	 - 	$ & $	0.12	$ & $	 - 	$ & $	0.15	$ & $	149,890	$ & $	 - 	$ & $	6.13	$ & $	 - 	$ & $	0.79	$	\\
				\surogate 	& $	1,910	$ & $	143	$ & $	95.54	$ & $	47.70	$ & $	0.39	$ & $	238,198	$ & $	14,053	$ & $	1,698.33	$ & $	1,732.56	$ & $	0.53	$	\\
				\sta 	& $	2,459	$ & $	52	$ & $	0.78	$ & $	4.29	$ & $	0.14	$ & $	173,779	$ & $	77,085	$ & $	13.75	$ & $	14.66	$ & $	0.72	$	\\
				\avoa 	& $	2,250	$ & $	240	$ & $	0.03	$ & $	0.01	$ & $	0.17	$ & $	199,707	$ & $	70,184	$ & $	2.74	$ & $	3.08	$ & $	0.68	$	\\
				\be 	& $	2,305	$ & $	228	$ & $	0.63	$ & $	0.26	$ & $	0.21	$ & $	207,543	$ & $	84,412	$ & $	27.13	$ & $	32.86	$ & $	0.59	$	\\
				\amc 	& $	2,156	$ & $	177	$ & $	0.11	$ & $	0.03	$ & $	0.34	$ & $	279,752	$ & $	97,760	$ & $	14.33	$ & $	12.66	$ & $	0.44	$	\\
				\cnna 	& $	2,436	$ & $	80	$ & $	0.07	$ & $	0.01	$ & $	0.06	$ & $	171,677	$ & $	65,071	$ & $	3.69	$ & $	4.69	$ & $	0.71	$	\\
				\nna	& $	2,425	$ & $	97	$ & $	0.13	$ & $	0.47	$ & $	0.10	$ & $	198,713	$ & $	69,676	$ & $	1.82	$ & $	2.36	$ & $	0.66	$	\\
				\aos 	& $	2,352	$ & $	195	$ & $	0.15	$ & $	0.03	$ & $	0.11	$ & $	232,893	$ & $	59,090	$ & $	13.75	$ & $	14.66	$ & $	0.63	$	\\
				\hs 	& $	2,397	$ & $	134	$ & $	0.31	$ & $	0.03	$ & $	0.11	$ & $	271,837	$ & $	79,842	$ & $	35.62	$ & $	30.14	$ & $	0.54	$	\\
				\abc 	& $	2,479	$ & $	35	$ & $	0.36	$ & $	0.02	$ & $	0.02	$ & $	316,453	$ & $	46,668	$ & $	43.34	$ & $	30.77	$ & $	0.44	$	\\
				\cs 	& $	2,463	$ & $	60	$ & $	0.03	$ & $	0.00	$ & $	0.05	$ & $	301,850	$ & $	146,932	$ & $	3.55	$ & $	2.88	$ & $	0.40	$	\\
				\pa 	& $	2,285	$ & $	233	$ & $	0.31	$ & $	0.07	$ & $	0.19	$ & $	384,786	$ & $	106,816	$ & $	53.36	$ & $	29.67	$ & $	0.23	$	\\
				\beso 	& $	2,429	$ & $	98	$ & $	0.08	$ & $	0.01	$ & $	0.08	$ & $	358,556	$ & $	87,810	$ & $	10.91	$ & $	6.86	$ & $	0.30	$	\\
				\bads	& $	2,441	$ & $	102.19	$ & $	23.96	$ & $	26.56	$ & $	0.05	$ & $	342,512	$ & $	165,686	$ & $	2,244.44$ & $	1,537.23	$ & $	0.32	$	\\
				\ica 	& $	2,277	$ & $	290	$ & $	0.78	$ & $	0.24	$ & $	0.15	$ & $	408,300	$ & $	115,624	$ & $	121.72	$ & $	60.15	$ & $	0.18	$	\\
				\mc 	& $	2,488	$ & $	21	$ & $	0.12	$ & $	0.43	$ & $	0.01	$ & $	379,183	$ & $	31,939	$ & $	15.71	$ & $	8.14	$ & $	0.30	$	\\
				\sfla 	& $	2,279	$ & $	282	$ & $	0.12	$ & $	0.02	$ & $	0.14	$ & $	428,388	$ & $	89,571	$ & $	21.73	$ & $	8.93	$ & $	0.14	$	\\
				\tlbo 	& $	2,465	$ & $	45	$ & $	0.36	$ & $	0.05	$ & $	0.04	$ & $	394,521	$ & $	79,375	$ & $	55.15	$ & $	28.56	$ & $	0.22	$	\\
				\ga 	& $	2,339	$ & $	233	$ & $	2.46	$ & $	8.31	$ & $	0.11	$ & $	441,384	$ & $	84,982	$ & $	33.34	$ & $	17.85	$ & $	0.12	$	\\
				\fstde 	& $	2,428	$ & $	104	$ & $	0.46	$ & $	5.76	$ & $	0.07	$ & $	439,054	$ & $	77,954	$ & $	14.10	$ & $	27.40	$ & $	0.13	$	\\
				\iwo 	& $	2,318	$ & $	252	$ & $	0.26	$ & $	0.06	$ & $	0.09	$ & $	448,337	$ & $	72,598	$ & $	49.32	$ & $	17.36	$ & $	0.10	$	\\
				\fa 	& $	2,347	$ & $	238	$ & $	0.41	$ & $	0.10	$ & $	0.07	$ & $	463,779	$ & $	56,774	$ & $	72.76	$ & $	22.25	$ & $	0.07	$	\\
				\bbo 	& $	2,397	$ & $	165	$ & $	0.38	$ & $	0.07	$ & $	0.06	$ & $	469,301	$ & $	48,668	$ & $	73.59	$ & $	19.68	$ & $	0.06	$	\\
				\aco 	& $	2,458	$ & $	63	$ & $	0.35	$ & $	0.03	$ & $	0.03	$ & $	472,440	$ & $	40,534	$ & $	77.53	$ & $	20.99	$ & $	0.06	$	\\
				\de 	& $	2,414	$ & $	137	$ & $	0.56	$ & $	0.09	$ & $	0.04	$ & $	477,775	$ & $	35,432	$ & $	107.72	$ & $	26.06	$ & $	0.04	$	\\
				\bat 	& $	2,436	$ & $	100	$ & $	0.04	$ & $	0.01	$ & $	0.03	$ & $	483,269	$ & $	26,205	$ & $	8.10	$ & $	1.85	$ & $	0.03	$	\\
				\ca 	& $	2,494	$ & $	3	$ & $	0.32	$ & $	0.01	$ & $	0.01	$ & $	498,013	$ & $	1,477	$ & $	60.25	$ & $	5.63	$ & $	0.01	$	\\
				\cmaes	& $	2,495	$ & $	5	$ & $	0.03	$ & $	0.00	$ & $	0.01	$ & $	498,754	$ & $	2,605	$ & $	6.44	$ & $	1.92	$ & $	0.01	$	\\
				
				\bottomrule
		\end{tabular}}
		\label{tab:results1}
	\end{table}
	
	Next, based on the obtained results (success rates), all $64$ algorithms are divided into four categories, Q1-Q4, sixteen in each.
	Q1 shows the top $16$ performing algorithms, Q2 the $17$th to $32$th places, Q3 the $33$st to $48$th, and Q4 the $49$th to $64$th.
	The data profiles on the left side of \Cref{fig:simple} illustrate the algorithmic performance  on ``simple'' \gkls{}-type problems (see \Cref{tab:paramet}) using $M_{\rm max} = 2,500$.
	The horizontal axis indicates the progress of the specific solver as the number of objective function evaluations slowly reached $2,500$.
	The left upper graph of \Cref{fig:simple} revealed that within a low function evaluation budget ($\leq 800$), there are at least six algorithms (\direct, \directrest, \glbsolve, \directm, \birmin, \dtco) which perform very similarly and can solve about $50 \%$ of all ``simple'' problems.
	Further, when the budget increases ($\geq 800$), the \birmin{} algorithm is the most efficient and solves $73 \%$ of test problems.
	Unexpectedly, all sixteen algorithms in the Q1 category are deterministic and even fifteen of them are \direct-type algorithms.
	Furthermore, $10$ out of the top $16$ in category Q2 belong to the \direct-type class.
	
	The best average results were achieved among stochastic type algorithms using the \surogate{} algorithm.
	However, it ranks only $20$th among all the algorithms and solved only $56 \%$ of these test problems.
	The second, third, and fourth best solvers in the pool of stochastic algorithms are two meta-heuristic \pa{} and \be{} approaches and the adaptive Monte-Carlo algorithm (\amc).
	Compared to the deterministic derivative-free algorithms, the best meta-heuristic \pa{} algorithm outperform only four of them, the \disimplc, \directaggress, \dirmin, and \lgo{} algorithms.
	
	Another important observation is that among the top $30$ algorithms, only five (\birmin, \mqnlp, \directrev, \oqnlp, and \mcs) are hybridized.
	It shows that local search procedures are not always practical when used excessively, and function evaluation budgets are limited.
	Algorithms such as \birmin{} and \directrev, which only use local searches when some improvement in the objective function is obtained, are more efficient than other strategy-based multi-start algorithms.
	
	\begin{figure*}[htp]
		\resizebox{\textwidth}{!}{
}
		\caption{Data profiles: fraction of ``simple'' \gkls{} test problems solved within an allowable number of function evaluations show on the $x$-axis, while the $y$-axis shows the proportion of solved problems. On the left-side $M_{\rm max}=2,500$, while on the right-side $M_{\rm max}=500,000$ is used. Q1 shows the top-$16$ performing algorithms, Q2 the $17$th to $32$th places, Q3 the $33$st to $48$th, and Q4 the $49$th to $64$th.}
		\label{fig:simple}
	\end{figure*}
	
	The fourth and fifth columns in \Cref{tab:results0} presents the average number of execution time in seconds ($t_{\rm avg.}$) and standard deviation ($t_{\rm s.d.}$) for all $400$ ``simple'' \gkls{}-type problems using $M_{\rm max} = 2,500$.
	Despite the poor performance of meta-heuristic algorithms based on the function evaluation criteria, they are often much faster than deterministic ones (see the fourth column in \Cref{tab:results0}).
	Therefore, they can perform more function evaluations in the same time budget.
	The tree the fastest are the \avoa, \cs, and \cmaes{} algorithms.
	Among the fastest approaches, at least nine other approaches with an average speed of fewer than $0.1$ fractions of a second.
	Out of $64$ algorithms, only $18$ had an average time greater than one second.
	The model-based \surogate{} and \bads{} algorithms were the two slowest algorithms, of which the \surogate{} averaged more than one minute ($72.91$ s.).
	
	\subsubsection{Comparison on ``simple'' \gkls{} problems with $M_{\rm max} = 500,000$}
	
	From seventh to eleventh columns in \Cref{tab:results0} shows the average number of function evaluations ($f_{\rm avg.}$), standard deviation ($f_{\rm s.d.}$), average time ($f_{\rm avg.}$), standard deviation ($t_{\rm s.d.}$), and success rate for all $400$ GKLS ``simple'' class test problems using $M_{\rm max}=500,000$.
	As the maximal budget for function evaluation increased significantly, the number of fails has accordingly decreased for most of the algorithms.
	Even $21$ algorithms out of $64$ demonstrated a perfect success rate (all are deterministic).
	However, even with $M_{\rm max}=500,000$, six hybridized algorithms (\mqnlp, \oqnlp, \multistart, \lgo, \multimin, and \globalsearch) suffered from a lack of success and surprisingly are significantly outperformed by the standard \direct-type algorithms.
	Again, quite surprisingly, increasing the budget for objective function evaluations for some algorithms (\ca, \de, \sfla, \iwo, \bbo, \fa, \cmaes, and \bat) did not improve.
	The success ratio remained the same as within $M_{\rm max}=2,500$.
	
	Again, the best average results were achieved using the hybridized \birmin{} algorithm, closely followed by \birectgb.
	In contrast, the third \directrev{} and fourth \dbdvo{} best algorithms deliver approximately $48 \%$ and $52 \%$ worse overall average number of function evaluations.
	Among the stochastic approaches, the best average results were achieved using \sa{} algorithm, but it ranks only in $23$rd place among all the best-performing algorithms.
	Another two algorithms dominating the stochastic pool are two hybridized \multistart{} and \multimin{} approaches, requiring approximately $91.75 \%$ and $92.23 \%$ more function evaluations than \birmin.
	Also, when the number of function evaluations allowed is higher, twice as many stochastic algorithms make it into the top-$32$.
	
	The data profiles on the right-side of \Cref{fig:simple} illustrate the performances of the algorithms within $M_{\rm max}=500,000$.
	The data profiles in the top right graph of \Cref{fig:simple} reveal that all the best algorithms in the Q1 category are \direct-type and perform very competitively.
	When the budget of function evaluation increases ($\geq 800$), the \birmin{} algorithm gives slightly better performance efficiency and achieves an ideal success rate ($1.00$) within the $20,000$ budget.
	None of the $16$ Q1 algorithms needed more than $100,000$ function evaluations to reach the ideal success ratio.
	In the set of stochastic algorithms, the best performing algorithm was \sa, closely followed by \multistart{} and \globalsearch.
	All three stochastic algorithms were very close to the ideal success ratio ($0.99$), but the solvers needed to use a maximal budget ($M_{\rm max}=500,000$).
	
	\begin{figure*}[htp]
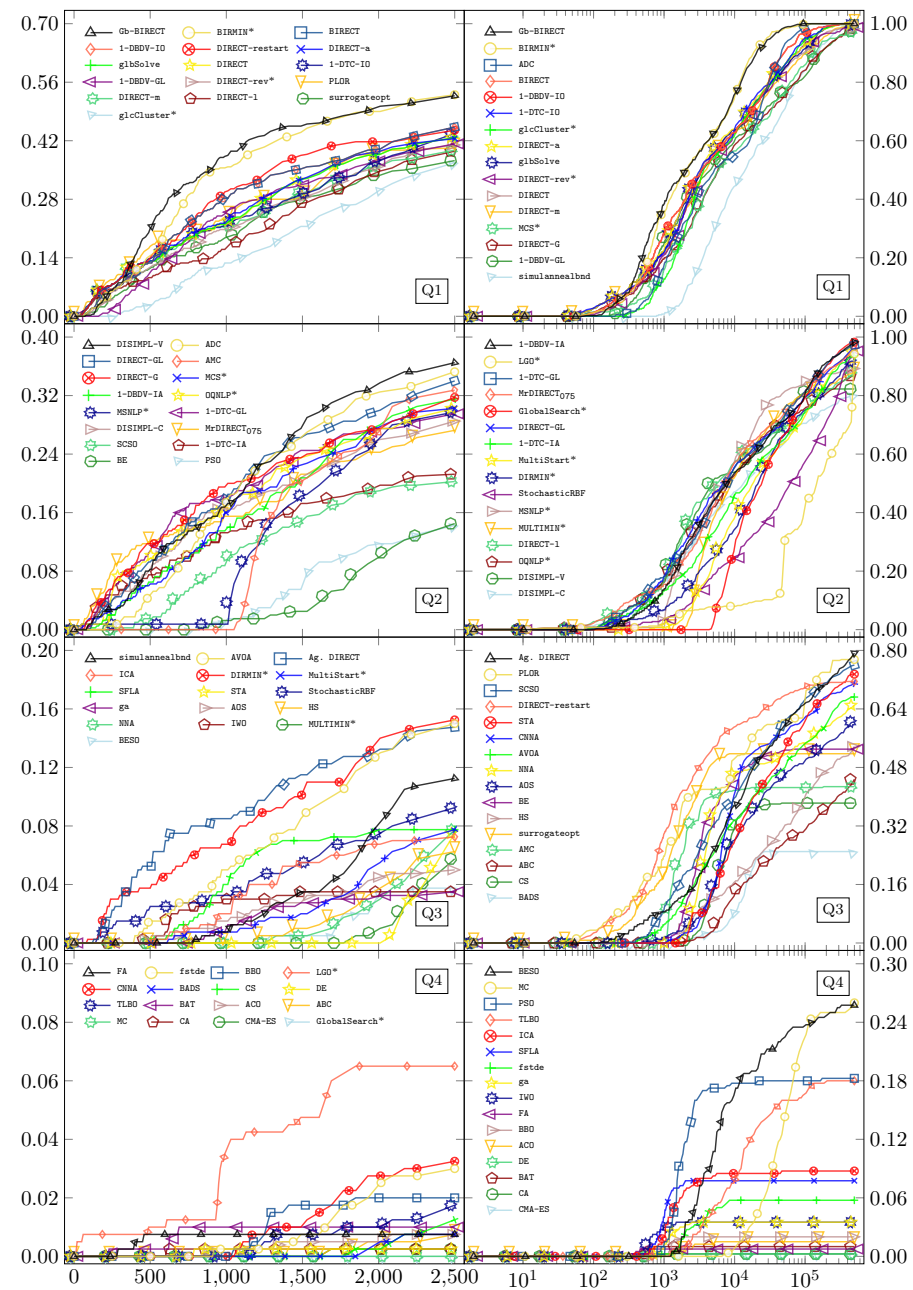

		\resizebox{\textwidth}{!}{
}
		\caption{Data profiles: fraction of ``hard'' \gkls{} test problems solved within an allowable number of function evaluations show on the $x$-axis, while the $y$-axis shows the proportion of solved problems. On the left-side $M_{\rm max}=2,500$, while on the right-side $M_{\rm max}=500,000$ is used. Q1 shows the top-$16$ performing algorithms, Q2 the $17$th to $32$th places, Q3 the $33$st to $48$th, and Q4 the $49$th to $64$th.}
		\label{fig:hard}
	\end{figure*}
	
	Among the $64$ algorithms, \dtcg{} was the algorithm that solved all the problems fastest.
	Only fourteen algorithms managed to fit into the one minute.
	All approaches belong to the \direct-type class.
	As the budget for evaluating objective functions is significantly higher than $M_{\rm max}=2,500$, some algorithms have proven to be very inefficient with respect to running time.
	Model-based algorithms are clearly not the best choice for solving problems with large budgets for evaluating objective functions.
	The average time to solve $400$ ``simple'' class problems with the \bads{} algorithm is as high as $1,955.32$ seconds.
	The standard deviation shows that the algorithm is often terminated due to a violation of the time-limited ($M_{\rm max} = 3,600$) stopping condition.
	
	\subsubsection{Comparison on ``hard'' \gkls{} problems with $M_{\rm max} = 2,500$}
	
	The summarized results are shown in the second to sixth columns of \Cref{tab:results1}.
	When the complexity increased, the performance efficiency significantly worsened, and the number of fails, was also increased.
	The \birmin{} and \birectgb{} algorithms share the largest number ($214$) of solved test problems.
	Quite surprisingly the non-hybridized \birectgb{} algorithm requires approximately $3\%$ fewer objective function evaluations.
	The third and fourth-best algorithms were \birect{} and \dbdvo.
	Both solved slightly less than half ($184$) of these test problems within $M_{\rm max}=2,500$.
	Interestingly, the top four algorithms are based on the same partitioning scheme~\cite{Paulavicius2016:jogo}.
	
	Among the stochastic algorithms, in the ``hard'' \gkls{} class, the \surogate{} algorithm is still the best but ranks only $15$th.
	The second and third best averages results and success rates were achieved by the \amc{} and \scso{} stochastic algorithms.
	The remaining stochastic algorithms could not solve more than a fifth of these test problems within $M_{\rm max}=2,500$.
	
	As in the previous study on $400$ ``simple'' test class \gkls{} problems, the fastest solvers (the lowest execution time) achieved with \avoa, \cs, and \cmaes{}.
	
	Finally, the data profiles on the left side of \Cref{fig:hard} illustrate the performances of algorithms on ``hard'' problems up to $M_{\rm max}=2,500$.
	The best two solvers, \birmin{} and \birectgb{}, share very similar performance and solved $54 \%$ of test instances, followed by \birect{} and \dbdvo{}, which solved $46 \%$.
	Among the stochastic algorithms, only the \surogate{} algorithm is the sole candidate in Q1.
	However, efficiency is one of the worst in this category.
	
	\subsubsection{Comparison on ``hard'' \gkls{} problems with $M_{\rm max} = 500,000$}
	
	Finally, the seventh to tenth columns in \Cref{tab:results1} presents the average and standard deviation values of function evaluations and execution times, while the last column --- success rate on all $400$ GKLS ``hard'' test problems using $M_{\rm max} = 500,000$.
	Only six algorithms achieved a perfect success rate compared to $21$ for the ``simple`` class.
	All of them are \direct-type algorithms.
	Once again, the two best-performing algorithms are \birmin{} and \birectgb.
	The \birmin{} required approximately $2.3$ and $2.5$ times fewer function evaluations than the third and fourth-best \adc{} and \dbdvo{} algorithms.
	The \sa{} approach ranks first among stochastic algorithms.
	However, compared to the best-performing \birmin{} algorithm, it required approximately $92 \%$ more objective function evaluations.
	
	This time, the average execution time revealed that no one algorithm could make it into the one-second time interval within $500,000$.
	The top three algorithms (including \birectgb{} and \birmin{}, which are the best performers in average function evaluations and success rates) had average times of just under two seconds.
	
	Finally, data profiles on the right-side of \Cref{fig:hard} show the performance on ``hard'' \gkls{} up to $M_{\rm max}=500,000$.
	The two \direct-type algorithms (\birmin{} and \birectgb) outperform the following best solvers in the Q1 group.
	As in the previous case, only one best stochastic algorithm (\sa) is in Q1.
	
	\subsection{Computational results solving \directgolib{} test problems}\label{sec:directgolib}
	
	\Cref{tab:results2} summarizes experimental results on $397$ test problems from \directgolib, where the best results are given in bold.
	As previously, the second to sixth columns presents the performance using with a low evaluation budget ($M_{\rm max}=2,500$), while from the seventh to eleventh columns show the performance when $M_{\rm max}=500,000$.
	
	\begin{table}
		\caption{The average number ($f_{\rm avg.}$) and standard deviation ($f_{\rm s.d.}$) of function evaluations, the average value ($t_{\rm avg.}$) and standard deviation ($t_{\rm s.d.}$) of execution time, and the success rate (S. R.) for \directgolib{} test problems using $M_{\rm max}=2,500$ and $M_{\rm max}=500,000$.}
		\resizebox{\textwidth}{!}{
			\begin{tabular}{@{\extracolsep{\fill}}l|rrrrr|rrrrr}
				\toprule
				Evaluation budget & \multicolumn{5}{c|}{$M_{\rm max}=2,500$} & \multicolumn{5}{c}{$M_{\rm max}=500,000$} \\
				\midrule
				Criteria & \multicolumn{1}{c}{$f_{\rm avg.}$} & \multicolumn{1}{c}{$f_{\rm s.d.}$} & \multicolumn{1}{c}{$t_{\rm avg.}$} & \multicolumn{1}{c}{$t_{\rm s.d.}$} & \multicolumn{1}{c|}{S. R.} & \multicolumn{1}{c}{$f_{\rm avg.}$} & \multicolumn{1}{c}{$f_{\rm s.d.}$} & \multicolumn{1}{c}{$t_{\rm avg.}$} & \multicolumn{1}{c}{$t_{\rm s.d.}$} & \multicolumn{1}{c}{S. R.} \\
				\midrule
				\avoa 	& $	 \mathbf{1,263} 	$ & $	180	$ & $	0.04	$ & $	0.17	$ & $	 \mathbf{0.62} 	$ & $	 \mathbf{113,651} 	$ & $	9,305	$ & $	1.67	$ & $	3.73	$ & $	 \mathbf{0.80} 	$	\\
				\oqnlp 	& $	1,315	$ & $	 - 	$ & $	0.68	$ & $	-	$ & $	0.54	$ & $	149,927	$ & $	 - 	$ & $	52.03	$ & $	-	$ & $	0.73	$	\\
				\mcs 	& $	1,389	$ & $	 - 	$ & $	7.01	$ & $	-	$ & $	0.53	$ & $	166,738	$ & $	 - 	$ & $	930.09	$ & $	-	$ & $	0.68	$	\\
				\beso 	& $	1,498	$ & $	56	$ & $	0.15	$ & $	0.11	$ & $	0.49	$ & $	178,163	$ & $	21,669	$ & $	14.43	$ & $	23.13	$ & $	0.67	$	\\
				\scso 	& $	1,458	$ & $	110	$ & $	0.19	$ & $	0.17	$ & $	0.52	$ & $	189,452	$ & $	12,209	$ & $	21.92	$ & $	42.57	$ & $	0.63	$	\\
				\mqnlp 	& $	1,480	$ & $	 - 	$ & $	0.78	$ & $	-	$ & $	0.50	$ & $	188,425	$ & $	 - 	$ & $	66.77	$ & $	-	$ & $	0.64	$	\\
				\lgo 	& $	1,528	$ & $	 - 	$ & $	2.25	$ & $	-	$ & $	0.44	$ & $	207,080	$ & $	 - 	$ & $	44.53	$ & $	-	$ & $	0.62	$	\\
				\dirmin 	& $	1,601	$ & $	 - 	$ & $	0.89	$ & $	-	$ & $	0.43	$ & $	202,417	$ & $	 - 	$ & $	70.78	$ & $	-	$ & $	0.61	$	\\
				\birmin 	& $	1,529	$ & $	 - 	$ & $	8.98	$ & $	-	$ & $	0.46	$ & $	218,943	$ & $	 - 	$ & $	1,376.18	$ & $	-	$ & $	0.57	$	\\
				\fstde 	& $	1,770	$ & $	61	$ & $	0.46	$ & $	0.28	$ & $	0.38	$ & $	185,126	$ & $	12,402	$ & $	45.49	$ & $	72.16	$ & $	0.65	$	\\
				\dtcg 	& $	1,747	$ & $	 - 	$ & $	1.92	$ & $	-	$ & $	0.37	$ & $	192,664	$ & $	 - 	$ & $	192.28	$ & $	-	$ & $	0.66	$	\\
				\directrev 	& $	1,597	$ & $	 - 	$ & $	18.08	$ & $	-	$ & $	0.45	$ & $	221,212	$ & $	 - 	$ & $	1,395.02	$ & $	-	$ & $	0.57	$	\\
				\aos 	& $	1,667	$ & $	106	$ & $	0.12	$ & $	0.08	$ & $	0.45	$ & $	227,009	$ & $	18,453	$ & $	12.41	$ & $	14.18	$ & $	0.57	$	\\
				\ca 	& $	1,536	$ & $	157	$ & $	0.64	$ & $	0.48	$ & $	0.48	$ & $	239,902	$ & $	23,049	$ & $	97.59	$ & $	165.04	$ & $	0.52	$	\\
				\dbdvg 	& $	1,838	$ & $	 - 	$ & $	1.26	$ & $	-	$ & $	0.33	$ & $	199,097	$ & $	 - 	$ & $	104.76	$ & $	-	$ & $	0.66	$	\\
				\multistart 	& $	1,705	$ & $	97	$ & $	0.26	$ & $	0.18	$ & $	0.40	$ & $	225,333	$ & $	8,417	$ & $	27.15	$ & $	54.90	$ & $	0.57	$	\\
				\dtca 	& $	1,799	$ & $	 - 	$ & $	0.87	$ & $	-	$ & $	0.35	$ & $	213,808	$ & $	 - 	$ & $	64.89	$ & $	-	$ & $	0.62	$	\\
				\dbdva 	& $	1,836	$ & $	 - 	$ & $	0.91	$ & $	-	$ & $	0.34	$ & $	216,331	$ & $	 - 	$ & $	60.69	$ & $	-	$ & $	0.62	$	\\
				\sta 	& $	2,421	$ & $	31	$ & $	0.35	$ & $	1.78	$ & $	0.27	$ & $	169,599	$ & $	21,235	$ & $	1.94	$ & $	3.92	$ & $	0.69	$	\\
				\dbdvo 	& $	1,713	$ & $	 - 	$ & $	7.71	$ & $	-	$ & $	0.39	$ & $	224,606	$ & $	 - 	$ & $	1,079.43	$ & $	-	$ & $	0.56	$	\\
				\dtco 	& $	1,704	$ & $	 - 	$ & $	8.69	$ & $	-	$ & $	0.38	$ & $	241,479	$ & $	 - 	$ & $	1,288.62	$ & $	-	$ & $	0.53	$	\\
				\directgl 	& $	1,860	$ & $	 - 	$ & $	0.47	$ & $	-	$ & $	0.31	$ & $	222,408	$ & $	 - 	$ & $	14.82	$ & $	-	$ & $	0.59	$	\\
				\plor 	& $	1,800	$ & $	 - 	$ & $	2.15	$ & $	-	$ & $	0.34	$ & $	233,973	$ & $	 - 	$ & $	251.50	$ & $	-	$ & $	0.55	$	\\
				\aco 	& $	1,861	$ & $	81	$ & $	1.16	$ & $	0.76	$ & $	0.35	$ & $	237,220	$ & $	20,908	$ & $	157.31	$ & $	246.80	$ & $	0.54	$	\\
				\bads	& $	1,576	$ & $	141	$ & $	454.06	$ & $	356.54	$ & $	0.44	$ & $	276,024	$ & $	10,680	$ & $	2,007.41	$ & $	1,754.28	$ & $	0.45	$	\\
				\tlbo 	& $	1,958	$ & $	68	$ & $	0.20	$ & $	0.10	$ & $	0.32	$ & $	238,584	$ & $	20,726	$ & $	19.92	$ & $	25.90	$ & $	0.56	$	\\
				\de 	& $	1,834	$ & $	96	$ & $	0.52	$ & $	0.29	$ & $	0.35	$ & $	243,004	$ & $	36,269	$ & $	63.24	$ & $	76.15	$ & $	0.53	$	\\
				\cnna 	& $	2,151	$ & $	84	$ & $	0.13	$ & $	0.06	$ & $	0.19	$ & $	187,284	$ & $	13,095	$ & $	8.05	$ & $	12.16	$ & $	0.68	$	\\
				\sfla 	& $	1,872	$ & $	80	$ & $	0.58	$ & $	0.37	$ & $	0.35	$ & $	245,699	$ & $	27,318	$ & $	77.09	$ & $	127.82	$ & $	0.52	$	\\
				\directg 	& $	1,871	$ & $	 - 	$ & $	0.44	$ & $	-	$ & $	0.31	$ & $	243,025	$ & $	 - 	$ & $	14.61	$ & $	-	$ & $	0.55	$	\\
				\glbsolve 	& $	1,855	$ & $	 - 	$ & $	1.78	$ & $	-	$ & $	0.30	$ & $	250,617	$ & $	 - 	$ & $	204.10	$ & $	-	$ & $	0.56	$	\\
				\pa 	& $	2,046	$ & $	75	$ & $	0.31	$ & $	0.13	$ & $	0.29	$ & $	227,322	$ & $	27,376	$ & $	29.57	$ & $	31.69	$ & $	0.56	$	\\
				\glccluster 	& $	2,108	$ & $	 - 	$ & $	13.61	$ & $	-	$ & $	0.25	$ & $	205,664	$ & $	 - 	$ & $	1,455.86	$ & $	-	$ & $	0.59	$	\\
				\ica 	& $	2,016	$ & $	73	$ & $	0.68	$ & $	0.32	$ & $	0.28	$ & $	235,327	$ & $	22,066	$ & $	76.64	$ & $	82.83	$ & $	0.54	$	\\
				\nna	& $	2,262	$ & $	62	$ & $	0.20	$ & $	0.07	$ & $	0.14	$ & $	182,670	$ & $	17,531	$ & $	12.70	$ & $	26.55	$ & $	0.67	$	\\
				\ga 	& $	2,020	$ & $	115	$ & $	0.16	$ & $	0.08	$ & $	0.27	$ & $	247,560	$ & $	35,022	$ & $	15.23	$ & $	20.07	$ & $	0.54	$	\\
				\cs 	& $	2,213	$ & $	66	$ & $	0.09	$ & $	0.05	$ & $	0.18	$ & $	219,341	$ & $	10,682	$ & $	5.49	$ & $	7.76	$ & $	0.62	$	\\
				\directl 	& $	1,835	$ & $	 - 	$ & $	9.62	$ & $	-	$ & $	0.31	$ & $	274,955	$ & $	 - 	$ & $	1,446.93	$ & $	-	$ & $	0.47	$	\\
				\bbo 	& $	1,953	$ & $	121	$ & $	1.28	$ & $	0.69	$ & $	0.31	$ & $	263,791	$ & $	29,302	$ & $	174.24	$ & $	233.95	$ & $	0.47	$	\\
				\birect 	& $	1,840	$ & $	 - 	$ & $	4.31	$ & $	-	$ & $	0.32	$ & $	280,776	$ & $	 - 	$ & $	644.96	$ & $	-	$ & $	0.45	$	\\
				\fa 	& $	1,957	$ & $	104	$ & $	0.36	$ & $	0.11	$ & $	0.28	$ & $	257,724	$ & $	26,112	$ & $	38.75	$ & $	37.66	$ & $	0.49	$	\\
				\birectgb 	& $	1,864	$ & $	 - 	$ & $	6.02	$ & $	-	$ & $	0.31	$ & $	284,990	$ & $	 - 	$ & $	914.77	$ & $	-	$ & $	0.45	$	\\
				\globalsearch 	& $	2,470	$ & $	12	$ & $	0.21	$ & $	0.04	$ & $	0.11	$ & $	198,720	$ & $	14,353	$ & $	13.44	$ & $	20.18	$ & $	0.65	$	\\
				\directmro 	& $	1,875	$ & $	 - 	$ & $	1.53	$ & $	-	$ & $	0.29	$ & $	281,020	$ & $	 - 	$ & $	180.80	$ & $	-	$ & $	0.46	$	\\
				\directaggress 	& $	1,991	$ & $	 - 	$ & $	0.47	$ & $	-	$ & $	0.26	$ & $	276,015	$ & $	 - 	$ & $	10.69	$ & $	-	$ & $	0.48	$	\\
				\hs 	& $	2,078	$ & $	99	$ & $	0.52	$ & $	0.23	$ & $	0.24	$ & $	268,195	$ & $	16,769	$ & $	63.72	$ & $	67.14	$ & $	0.50	$	\\
				\be 	& $	2,102	$ & $	77	$ & $	0.40	$ & $	0.10	$ & $	0.24	$ & $	262,122	$ & $	16,767	$ & $	39.29	$ & $	38.12	$ & $	0.49	$	\\
				\surogate 	& $	1,714	$ & $	74	$ & $	387.54	$ & $	313.01	$ & $	0.35	$ & $	316,040	$ & $	5,711	$ & $	2,367.07	$ & $	1,654.37	$ & $	0.37	$	\\
				\direct 	& $	1,900	$ & $	 - 	$ & $	1.61	$ & $	-	$ & $	0.28	$ & $	284,593	$ & $	 - 	$ & $	195.63	$ & $	-	$ & $	0.44	$	\\
				\directm 	& $	1,893	$ & $	 - 	$ & $	7.58	$ & $	-	$ & $	0.28	$ & $	286,145	$ & $	 - 	$ & $	112.62	$ & $	-	$ & $	0.44	$	\\
				\directrest 	& $	1,857	$ & $	 - 	$ & $	6.32	$ & $	-	$ & $	0.30	$ & $	295,712	$ & $	 - 	$ & $	1,001.60	$ & $	-	$ & $	0.42	$	\\
				\cmaes	& $	2,025	$ & $	69	$ & $	0.29	$ & $	0.24	$ & $	0.27	$ & $	285,322	$ & $	27,388	$ & $	36.74	$ & $	63.86	$ & $	0.45	$	\\
				\amc 	& $	2,048	$ & $	91	$ & $	0.16	$ & $	0.07	$ & $	0.26	$ & $	283,518	$ & $	16,866	$ & $	16.50	$ & $	23.23	$ & $	0.44	$	\\
				\directa 	& $	1,905	$ & $	 - 	$ & $	3.32	$ & $	-	$ & $	0.28	$ & $	303,547	$ & $	 - 	$ & $	160.28	$ & $	-	$ & $	0.41	$	\\
				\StochasticRBF	& $	1,838	$ & $	47	$ & $	19.20	$ & $	76.18	$ & $	0.30	$ & $	316,717	$ & $	2,573	$ & $	2,285.50	$ & $	1,715.81	$ & $	0.37	$	\\
				\iwo 	& $	2,111	$ & $	88	$ & $	0.31	$ & $	0.13	$ & $	0.22	$ & $	281,501	$ & $	26,616	$ & $	36.92	$ & $	32.58	$ & $	0.45	$	\\
				\abc 	& $	2,213	$ & $	67	$ & $	0.41	$ & $	0.14	$ & $	0.17	$ & $	276,920	$ & $	4,196	$ & $	45.86	$ & $	40.19	$ & $	0.49	$	\\
				\adc 	& $	1,876	$ & $	 - 	$ & $	13.20	$ & $	-	$ & $	0.29	$ & $	326,049	$ & $	 - 	$ & $	2,348.85	$ & $	-	$ & $	0.35	$	\\
				\disimplv 	& $	1,938	$ & $	 - 	$ & $	1,996.45	$ & $	-	$ & $	0.26	$ & $	323,428	$ & $	 - 	$ & $	2,335.31	$ & $	-	$ & $	0.36	$	\\
				\multimin 	& $	2,425	$ & $	15	$ & $	0.17	$ & $	0.04	$ & $	0.07	$ & $	238,069	$ & $	12,135	$ & $	11.99	$ & $	15.84	$ & $	0.53	$	\\
				\sa 	& $	2,046	$ & $	87	$ & $	1.14	$ & $	0.51	$ & $	0.23	$ & $	328,016	$ & $	9,129	$ & $	177.11	$ & $	181.02	$ & $	0.36	$	\\
				\disimplc 	& $	2,009	$ & $	 - 	$ & $	1,997.31	$ & $	-	$ & $	0.22	$ & $	349,459	$ & $	 - 	$ & $	2,256.41	$ & $	-	$ & $	0.32	$	\\
				\bat 	& $	2,008	$ & $	122	$ & $	0.06	$ & $	0.04	$ & $	0.25	$ & $	369,555	$ & $	28,028	$ & $	4.55	$ & $	4.28	$ & $	0.26	$	\\
				\mc 	& $	2,319	$ & $	61	$ & $	0.16	$ & $	0.05	$ & $	0.10	$ & $	402,594	$ & $	11,360	$ & $	19.98	$ & $	15.86	$ & $	0.22	$	\\
				
				\bottomrule
		\end{tabular}}
		\label{tab:results2}
	\end{table}
	
	Unlike in \Cref{sec:gklstest}, where deterministic algorithms dominate over stochastic, the stochastic solvers are very competitive and often perform better.
	There have been five algorithms (\avoa, \oqnlp, \mcs, \scso, and \mqnlp) that managed to achieve a success rate of $\ge 0.5$ within $M_{\rm max}=2,500$.
	Moreover, nine algorithms (\avoa, \oqnlp, \sta, \mcs, \cnna, \amc, \beso, \dtcg, and \dbdvg) were able to solve at least two-thirds of the test problems when $M_{\rm max}=500,000$.
	Among all problems, the meta-heuristic nature-inspired \avoa{} algorithm showed the best performance in both cases (for $M_{\rm max}=2,500$, and $M_{\rm max}=500,000$).
	The success rates of the \avoa{} algorithm (in the $M_{\rm max}=2,500$ and $M_{\rm max}=500,000$ cases) are $62 \%$ and $80 \%$.
	The \avoa{} also proves to be the fastest when comparing the average execution time.
	
	\begin{figure*}[htp]
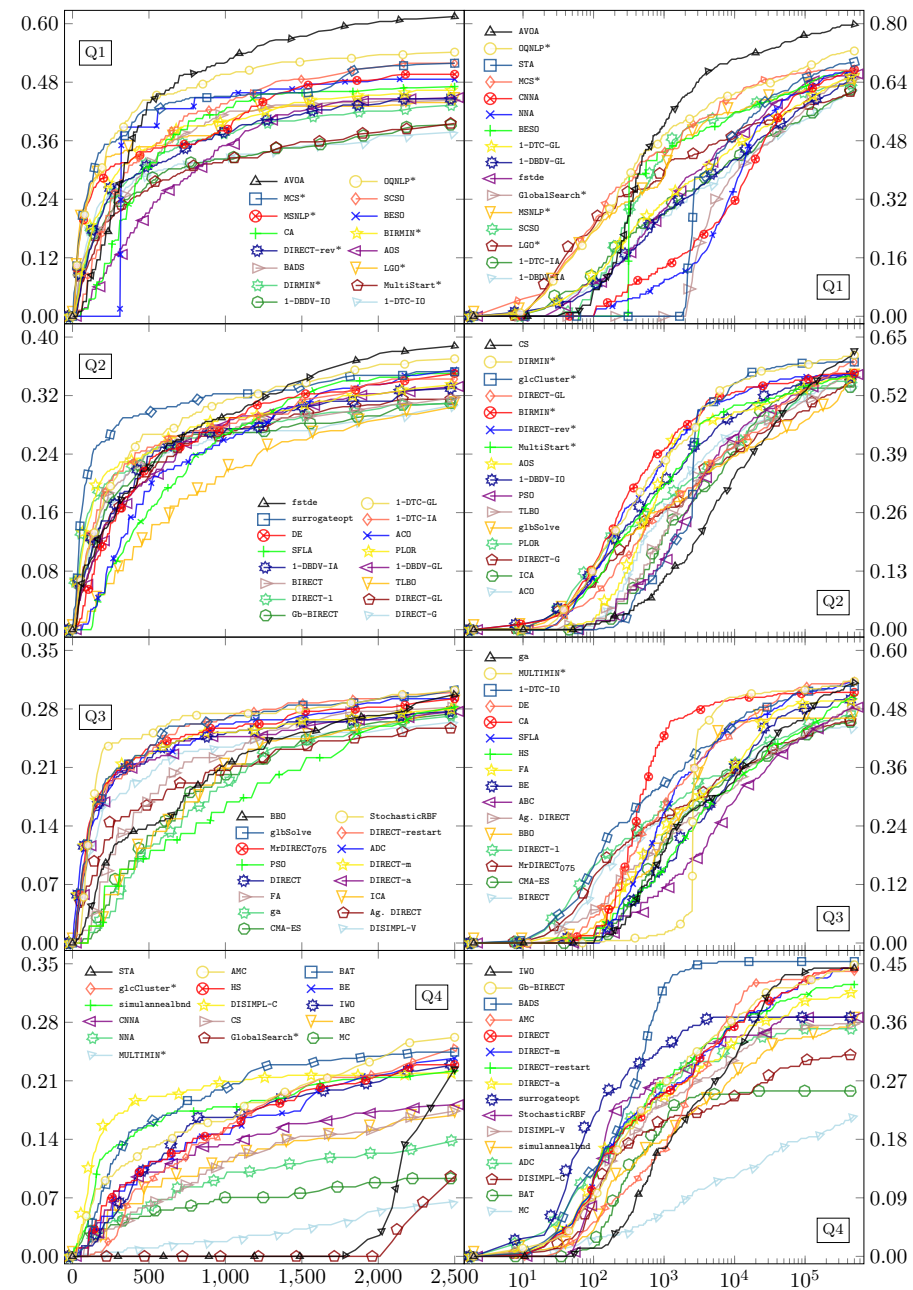

		\resizebox{\textwidth}{!}{
}
		\caption{Data profiles: fraction of $397$ \directgolib{} test problems solved as a function of allowable number of function evaluation. The horizontal coordinate equals a number of function evaluations or median value while the vertical coordinate equals a proportion of solved problems. The evaluation budget on the left is $2,500$, while on the right it is $500,000$. Also Q1 shows the top-$16$ performing algorithms, Q2 the $17$th to $32$th place, Q3 the $33$st to $48$th and Q4 the $49$th to $64$th.}
		\label{fig:overall}
	\end{figure*}
	
	The data profiles on the left-side of \Cref{fig:overall} illustrate the performance of algorithms up to $M_{\rm max}=2,500$.
	When the budget is meager ($\leq 500$), the best performing algorithms are \oqnlp{} and \mcs{} and solved approximately $ 45 \%$ of all problems from \directgolib{}.
	However, when the budget increases ($ > 500$), the \avoa{} algorithm starts outperforming others.
	Among the top-$32$ solvers, fifteen are \direct-type, while ten are meta-heuristics.
	
	Many hybridized algorithms did not perform well in the previous \Cref{sec:gklstest}.
	However, lots of the \directgolib{} test problems' dimensions are much higher than the \gkls, and the hybridized solvers performed much better this time.
	Among the sixteen best-performing algorithms (Q1), eight are hybridized.
	
	The data profiles on the right-side of \Cref{fig:overall} illustrate the performances using $M_{\rm max}=500,000$.
	The best two performing algorithms are the \avoa{} and \oqnlp. 
	Using a higher budget, the efficiency of non-hybridized solvers has improved significantly.
	Now eleven out of the top sixteen algorithms are non-hybridized.
	Only six of the sixteen algorithms made into the Q1 when $M_{\rm max}=2,500$ remained here using $M_{\rm max}=500,000$.
	
\begin{figure*}[htp]
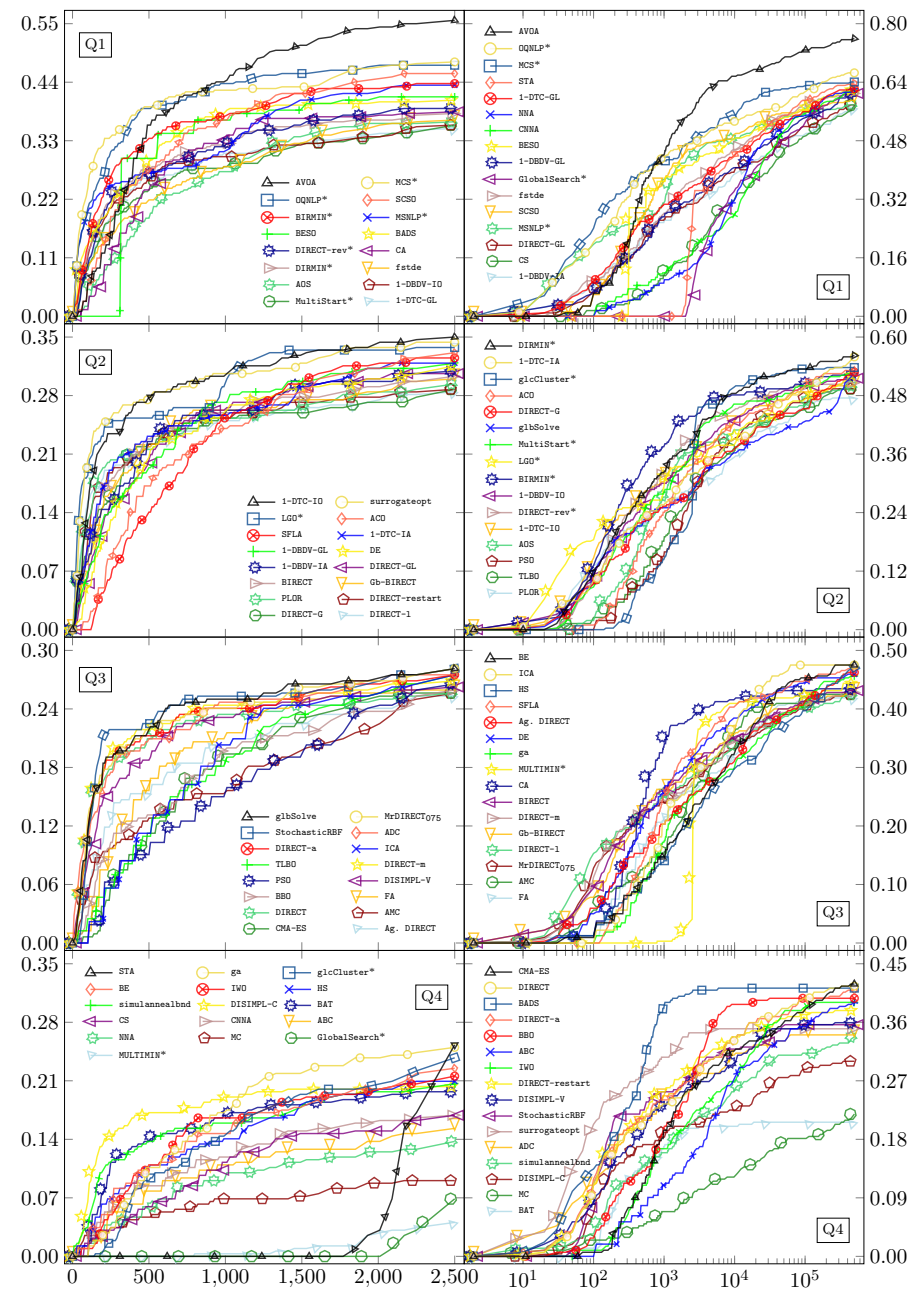

		\resizebox{\textwidth}{!}{
}
		\caption{Data profiles: fraction of $320$ multi-modal and/or non-convex \directgolib{} test problems solved within an allowable number of function evaluations show on the $x$-axis, while the $y$-axis shows the proportion of solved problems. On the left-side $M_{\rm max}=2,500$, while on the right-side $M_{\rm max}=500,000$ is used. Q1 shows the top-$16$ performing algorithms, Q2 the $17$th to $32$th places, Q3 the $33$st to $48$th, and Q4 the $49$th to $64$th.}
		\label{fig:multimodal}
	\end{figure*}

	Additionally, data profiles in \Cref{fig:multimodal,fig:unimodal,fig:low,fig:high} demonstrate a more detailed performance of the algorithms on different subsets of the \directgolib.
	The first observation is that four algorithms always make it into the Q1 --- \avoa, \oqnlp, \mcs{}, and \mqnlp.
	The \avoa{} algorithm demonstrated the most significant advantage in solving high-dimensional, multi-modal, and non-convex test problems.
	On simpler uni-modal and/or convex test problems, the \lgo{} was the best performing algorithm, followed closely by \oqnlp, \avoa{}, and \sta{}.
	For low-dimensional test problems, the best results are obtained using \direct-type algorithms, especially when $M_{\rm max}=500,000$.
	
	\begin{figure*}[htp]
		\resizebox{\textwidth}{!}{
			%
}
		\caption{Data profiles: fraction of $97$ uni-modal and/or convex \directgolib{} test problems solved within an allowable number of function evaluations show on the $x$-axis, while the $y$-axis shows the proportion of solved problems. On the left-side $M_{\rm max}=2,500$, while on the right-side $M_{\rm max}=500,000$ is used. Q1 shows the top-$16$ performing algorithms, Q2 the $17$th to $32$th places, Q3 the $33$st to $48$th, and Q4 the $49$th to $64$th.}
		\label{fig:unimodal}
	\end{figure*}
	
	\begin{figure*}[htp]
		\resizebox{\textwidth}{!}{
			%
}
		\caption{Data profiles: fraction of $231$ low-dimensional ($n \leq 10$) \directgolib{} test problems solved within an allowable number of function evaluations show on the $x$-axis, while the $y$-axis shows the proportion of solved problems. On the left-side $M_{\rm max}=2,500$, while on the right-side $M_{\rm max}=500,000$ is used. Q1 shows the top-$16$ performing algorithms, Q2 the $17$th to $32$th places, Q3 the $33$st to $48$th, and Q4 the $49$th to $64$th.}
		\label{fig:low}
	\end{figure*}
	
	\begin{figure*}[htp]
		\resizebox{\textwidth}{!}{
			%
}
		\caption{Data profiles: fraction of $166$ high-dimensional ($n > 10$) \directgolib{} test problems solved within an allowable number of function evaluations show on the $x$-axis, while the $y$-axis shows the proportion of solved problems. On the left-side $M_{\rm max}=2,500$, while on the right-side $M_{\rm max}=500,000$ is used. Q1 shows the top-$16$ performing algorithms, Q2 the $17$th to $32$th places, Q3 the $33$st to $48$th, and Q4 the $49$th to $64$th.}
		\label{fig:high}
	\end{figure*}
	
	\begin{table}
		\centering
		\caption{The average number of function evaluations ($f_{\rm avg.}$), standard deviation (S. D.), and success rate (S. R.) when the solution is: at zero-point ($122$ cases), and at other non-zero locations ($198$ cases).}
			\fontsize{7.5pt}{7.5pt}\selectfont
			\begin{tabular*}{\textwidth}{@{\extracolsep{\fill}}l|rrrr|rrrr}
				\toprule
				Evaluation budget & \multicolumn{4}{c|}{$M_{\rm max}=2,500$} & \multicolumn{4}{c}{$M_{\rm max}=500,000$} \\
				\midrule
				Minimum location & \multicolumn{2}{c}{Zero-point} & \multicolumn{2}{c|}{Non-zero} & \multicolumn{2}{c}{Zero-point} & \multicolumn{2}{c}{Non-zero} \\
				\midrule
				Criteria & $f_{\rm avg.}$ & S. R. & $f_{\rm avg.}$ & S. R. & $f_{\rm avg.}$ & S. R. & $f_{\rm avg.}$ & S. R. \\
				\midrule
				\avoa	& $	972	$ & $	0.81	$ & $	1,671	$ & $	0.40	$ & $	21,239	$ & $	0.97	$ & $	198,287	$ & $	0.63	$	\\
				\oqnlp	& $	1,525	$ & $	0.47	$ & $	1,423	$ & $	0.47	$ & $	157,898	$ & $	0.72	$ & $	196,557	$ & $	0.63	$	\\
				\mcs	& $	1,508	$ & $	0.47	$ & $	1,432	$ & $	0.50	$ & $	168,392	$ & $	0.69	$ & $	194,373	$ & $	0.62	$	\\
				\scso	& $	1,207	$ & $	0.70	$ & $	1,842	$ & $	0.31	$ & $	109,500	$ & $	0.80	$ & $	272,537	$ & $	0.46	$	\\
				\beso	& $	1,237	$ & $	0.62	$ & $	1,914	$ & $	0.30	$ & $	143,461	$ & $	0.76	$ & $	243,672	$ & $	0.53	$	\\
				\mqnlp	& $	1,677	$ & $	0.42	$ & $	1,615	$ & $	0.44	$ & $	226,343	$ & $	0.56	$ & $	215,091	$ & $	0.59	$	\\
				\fstde	& $	1,882	$ & $	0.36	$ & $	1,751	$ & $	0.36	$ & $	191,994	$ & $	0.67	$ & $	222,718	$ & $	0.56	$	\\
				\dtcg	& $	1,863	$ & $	0.33	$ & $	1,727	$ & $	0.37	$ & $	225,135	$ & $	0.58	$ & $	193,966	$ & $	0.65	$	\\
				\dirmin	& $	1,658	$ & $	0.42	$ & $	1,741	$ & $	0.36	$ & $	206,652	$ & $	0.61	$ & $	240,291	$ & $	0.53	$	\\
				\sta	& $	2,433	$ & $	0.32	$ & $	2,398	$ & $	0.26	$ & $	143,206	$ & $	0.76	$ & $	231,138	$ & $	0.56	$	\\
				\birmin	& $	1,647	$ & $	0.40	$ & $	1,534	$ & $	0.46	$ & $	252,237	$ & $	0.51	$ & $	237,463	$ & $	0.53	$	\\
				\ca	& $	1,359	$ & $	0.62	$ & $	1,917	$ & $	0.26	$ & $	173,093	$ & $	0.67	$ & $	339,296	$ & $	0.30	$	\\
				\dbdvg	& $	1,904	$ & $	0.30	$ & $	1,839	$ & $	0.33	$ & $	236,293	$ & $	0.57	$ & $	207,497	$ & $	0.64	$	\\
				\directrev	& $	1,714	$ & $	0.39	$ & $	1,702	$ & $	0.39	$ & $	231,957	$ & $	0.55	$ & $	261,467	$ & $	0.48	$	\\
				\aos	& $	1,771	$ & $	0.45	$ & $	1,847	$ & $	0.32	$ & $	233,569	$ & $	0.57	$ & $	275,762	$ & $	0.46	$	\\
				\dbdva	& $	1,906	$ & $	0.30	$ & $	1,846	$ & $	0.33	$ & $	250,008	$ & $	0.55	$ & $	229,805	$ & $	0.59	$	\\
				\lgo	& $	1,594	$ & $	0.41	$ & $	1,862	$ & $	0.30	$ & $	247,416	$ & $	0.54	$ & $	260,493	$ & $	0.52	$	\\
				\multistart	& $	1,941	$ & $	0.33	$ & $	1,687	$ & $	0.38	$ & $	260,109	$ & $	0.52	$ & $	236,205	$ & $	0.53	$	\\
				\dtca	& $	1,949	$ & $	0.30	$ & $	1,794	$ & $	0.34	$ & $	266,705	$ & $	0.52	$ & $	224,916	$ & $	0.59	$	\\
				\aco	& $	1,907	$ & $	0.35	$ & $	1,908	$ & $	0.32	$ & $	223,801	$ & $	0.59	$ & $	255,436	$ & $	0.49	$	\\
				\directgl	& $	1,955	$ & $	0.26	$ & $	1,827	$ & $	0.33	$ & $	246,293	$ & $	0.56	$ & $	221,610	$ & $	0.59	$	\\
				\dbdvo	& $	1,842	$ & $	0.31	$ & $	1,725	$ & $	0.39	$ & $	259,482	$ & $	0.50	$ & $	242,540	$ & $	0.53	$	\\
				\dtco	& $	1,783	$ & $	0.34	$ & $	1,731	$ & $	0.36	$ & $	242,585	$ & $	0.52	$ & $	262,260	$ & $	0.48	$	\\
				\cnna	& $	2,238	$ & $	0.20	$ & $	2,146	$ & $	0.17	$ & $	188,573	$ & $	0.71	$ & $	232,406	$ & $	0.56	$	\\
				\directg	& $	1,999	$ & $	0.25	$ & $	1,846	$ & $	0.32	$ & $	267,858	$ & $	0.52	$ & $	248,495	$ & $	0.54	$	\\
				\bads	& $	1,486	$ & $	0.42	$ & $	1,645	$ & $	0.39	$ & $	246,506	$ & $	0.43	$ & $	298,821	$ & $	0.39	$	\\
				\glbsolve	& $	1,973	$ & $	0.26	$ & $	1,852	$ & $	0.30	$ & $	283,702	$ & $	0.48	$ & $	260,338	$ & $	0.55	$	\\
				\nna	& $	2,286	$ & $	0.15	$ & $	2,243	$ & $	0.14	$ & $	124,519	$ & $	0.72	$ & $	227,581	$ & $	0.56	$	\\
				\glccluster	& $	2,205	$ & $	0.21	$ & $	2,078	$ & $	0.26	$ & $	215,249	$ & $	0.57	$ & $	245,672	$ & $	0.52	$	\\
				\tlbo	& $	2,102	$ & $	0.29	$ & $	1,960	$ & $	0.29	$ & $	283,226	$ & $	0.49	$ & $	261,589	$ & $	0.48	$	\\
				\sfla	& $	1,966	$ & $	0.34	$ & $	1,901	$ & $	0.30	$ & $	281,684	$ & $	0.47	$ & $	275,514	$ & $	0.44	$	\\
				\pa	& $	2,127	$ & $	0.28	$ & $	2,042	$ & $	0.27	$ & $	266,552	$ & $	0.51	$ & $	259,321	$ & $	0.48	$	\\
				\plor	& $	1,987	$ & $	0.25	$ & $	1,819	$ & $	0.32	$ & $	270,034	$ & $	0.48	$ & $	268,858	$ & $	0.48	$	\\
				\de	& $	1,804	$ & $	0.30	$ & $	1,858	$ & $	0.31	$ & $	195,629	$ & $	0.48	$ & $	279,593	$ & $	0.44	$	\\
				\cs	& $	2,271	$ & $	0.18	$ & $	2,183	$ & $	0.18	$ & $	255,804	$ & $	0.59	$ & $	232,084	$ & $	0.57	$	\\
				\ica	& $	2,020	$ & $	0.29	$ & $	2,012	$ & $	0.26	$ & $	200,028	$ & $	0.49	$ & $	262,589	$ & $	0.47	$	\\
				\ga	& $	2,079	$ & $	0.28	$ & $	2,038	$ & $	0.24	$ & $	262,609	$ & $	0.53	$ & $	286,766	$ & $	0.44	$	\\
				\birect	& $	1,937	$ & $	0.26	$ & $	1,831	$ & $	0.33	$ & $	297,452	$ & $	0.43	$ & $	285,983	$ & $	0.43	$	\\
				\birectgb	& $	1,978	$ & $	0.25	$ & $	1,825	$ & $	0.34	$ & $	315,525	$ & $	0.39	$ & $	278,499	$ & $	0.45	$	\\
				\globalsearch	& $	2,486	$ & $	0.10	$ & $	2,476	$ & $	0.07	$ & $	215,610	$ & $	0.66	$ & $	219,102	$ & $	0.59	$	\\
				\be	& $	2,112	$ & $	0.25	$ & $	2,095	$ & $	0.22	$ & $	239,363	$ & $	0.49	$ & $	279,700	$ & $	0.45	$	\\
				\hs	& $	2,048	$ & $	0.25	$ & $	2,101	$ & $	0.20	$ & $	231,553	$ & $	0.54	$ & $	296,494	$ & $	0.42	$	\\
				\directl	& $	1,886	$ & $	0.25	$ & $	1,795	$ & $	0.31	$ & $	277,431	$ & $	0.39	$ & $	273,043	$ & $	0.45	$	\\
				\bbo	& $	1,915	$ & $	0.29	$ & $	1,982	$ & $	0.26	$ & $	217,463	$ & $	0.45	$ & $	299,572	$ & $	0.39	$	\\
				\directaggress	& $	2,108	$ & $	0.20	$ & $	1,941	$ & $	0.29	$ & $	314,032	$ & $	0.39	$ & $	265,116	$ & $	0.51	$	\\
				\directmro	& $	1,999	$ & $	0.25	$ & $	1,854	$ & $	0.30	$ & $	314,277	$ & $	0.39	$ & $	282,801	$ & $	0.44	$	\\
				\surogate	& $	1,655	$ & $	0.35	$ & $	1,760	$ & $	0.33	$ & $	301,918	$ & $	0.36	$ & $	326,946	$ & $	0.34	$	\\
				\directm	& $	1,991	$ & $	0.25	$ & $	1,870	$ & $	0.28	$ & $	305,262	$ & $	0.40	$ & $	285,227	$ & $	0.44	$	\\
				\amc	& $	2,178	$ & $	0.25	$ & $	1,994	$ & $	0.26	$ & $	295,270	$ & $	0.44	$ & $	287,711	$ & $	0.42	$	\\
				\directa	& $	1,989	$ & $	0.25	$ & $	1,871	$ & $	0.29	$ & $	308,132	$ & $	0.40	$ & $	298,698	$ & $	0.41	$	\\
				\fa	& $	1,934	$ & $	0.27	$ & $	1,975	$ & $	0.25	$ & $	219,650	$ & $	0.41	$ & $	287,129	$ & $	0.42	$	\\
				\direct	& $	1,980	$ & $	0.25	$ & $	1,920	$ & $	0.26	$ & $	292,914	$ & $	0.43	$ & $	301,685	$ & $	0.41	$	\\
				\cmaes	& $	2,021	$ & $	0.25	$ & $	2,028	$ & $	0.25	$ & $	277,697	$ & $	0.39	$ & $	291,212	$ & $	0.44	$	\\
				\directrest	& $	1,988	$ & $	0.24	$ & $	1,827	$ & $	0.32	$ & $	335,458	$ & $	0.34	$ & $	300,210	$ & $	0.40	$	\\
				\StochasticRBF	& $	1,877	$ & $	0.24	$ & $	1,807	$ & $	0.30	$ & $	321,375	$ & $	0.33	$ & $	313,119	$ & $	0.37	$	\\
				\disimplv	& $	1,993	$ & $	0.23	$ & $	1,916	$ & $	0.28	$ & $	339,460	$ & $	0.33	$ & $	312,807	$ & $	0.38	$	\\
				\adc	& $	1,960	$ & $	0.25	$ & $	1,891	$ & $	0.30	$ & $	345,790	$ & $	0.31	$ & $	321,908	$ & $	0.36	$	\\
				\iwo	& $	2,117	$ & $	0.22	$ & $	2,106	$ & $	0.20	$ & $	251,732	$ & $	0.39	$ & $	304,491	$ & $	0.39	$	\\
				\sa	& $	2,018	$ & $	0.24	$ & $	2,068	$ & $	0.19	$ & $	319,207	$ & $	0.35	$ & $	334,819	$ & $	0.33	$	\\
				\abc	& $	2,218	$ & $	0.17	$ & $	2,209	$ & $	0.15	$ & $	251,156	$ & $	0.38	$ & $	296,818	$ & $	0.40	$	\\
				\multimin	& $	2,483	$ & $	0.08	$ & $	2,462	$ & $	0.02	$ & $	270,024	$ & $	0.50	$ & $	293,981	$ & $	0.40	$	\\
				\disimplc	& $	2,109	$ & $	0.18	$ & $	2,017	$ & $	0.22	$ & $	373,221	$ & $	0.27	$ & $	343,926	$ & $	0.33	$	\\
				\bat	& $	1,969	$ & $	0.20	$ & $	2,039	$ & $	0.21	$ & $	356,268	$ & $	0.20	$ & $	379,817	$ & $	0.22	$	\\
				\mc	& $	2,319	$ & $	0.13	$ & $	2,320	$ & $	0.07	$ & $	411,733	$ & $	0.22	$ & $	395,536	$ & $	0.22	$	\\
				\bottomrule
			\end{tabular*}
		\label{tab:results3}
	\end{table}
	
	\subsubsection{Investigating the impact of the solution lying at zero}
	Despite the outstanding performance of the \avoa{} and some other meta-heuristic algorithms, there is one crucial point that needs additional discussion.
	It is natural to assume that the algorithms are unbiased to the landscape of problems.
	However, during these experiments, it has been observed that some approaches give extra attention to exploring the region around the zero point.
	Of $397$ test problems, $173$ had solutions precisely at the zero point, and $122$ of these test problems were highly complex.
	Therefore, an additional study was carried out, and the results are summarized in \Cref{tab:results3}.
	Here we have selected two subgroups depending on where the solutions are: highly complex problems with solutions at the zero point ($122$ cases) and all problems with a solution at other locations ($198$ cases).
	
	Standard non-hybridized \direct-type algorithms make almost no difference whether the global solution is at the zero point or elsewhere.
	Furthermore, non-hybridized \direct-type solvers often perform better when the global solution is at a non-zero point.
	For example, the \directaggress{} algorithm manages to achieve $31 \%$ greater success rate on non-zero problems using $M_{\rm max}=2,500$, and about $23 \%$ when $M_{\rm max}=500,000$.
	
	The results of the hybridized algorithms are also very similar.
	Some algorithms work slightly better when the global solution is at the zero point, and in some cases, the opposite is true.
	Seven (out of eleven) hybridized algorithms better solved non-zero problems.
	When $M_{\rm max}=2,500$, the best performing algorithm on non-zero problems is \mcs, followed by \oqnlp{} and \birmin{} algorithms.
	When $M_{\rm max}=500,000$, the best results were achieved using the \oqnlp{} algorithm, closely followed by the \mcs{} and \mqnlp.
	However, in overall rating, the standard non-hybridized algorithms (\dtcg{} and \dbdvg) proved to be better at solving non-zero problems, especially when $M_{\rm max}=500,000$.
	
	The \avoa{} algorithm, which showed outstanding results in \Cref{tab:results2}, does not perform so well on non-zero problems.
	When $M_{\rm max}=2,500$, \avoa{} is outperformed by four algorithms, while when $M_{\rm max}=500,000$, it shares the third-fourth places.
	Completely different results are obtained with the \avoa{} algorithm on $122$ zero-point test problems.
	Using $M_{\rm max}=2,500$, the success ratio of the \avoa{} algorithm is $0.81$, i.e., more than twice that higher compared the success ratio ($0.40$) on non-zero problems.
	When $M_{\rm max}=500,000$, the \avoa{} success ratio (on zero-point problems) was close to the ideal ($0.97$), while on non-zero problems, it reached only $0.63$.
	The \avoa{} algorithm required about $89 \%$ fewer objective function evaluations on zero-point problems.
	Thus, even the \avoa{} algorithm's initial population is generated randomly over $D$, but the algorithm naturally tends to converge to a zero-point.
	
	Similar behavior has been observed for more meta-heuristic algorithms, such as \scso, \aos, \beso, \cnna, \ca, and \sta.
	Using $M_{\rm max}=2,500$ nine algorithms achieved more than $25 \%$ greater a success rate on zero-point problems.
	Eight of these nine algorithms are stochastic, and most of them ($5$) are population-based search (PBS) techniques.
	Similar trends persisted using $M_{\rm max}=500,000$.
	
	\section{Summary of numerical benchmarking}\label{sec:results}
	
	Two extensive different analyses on \gkls{} and \directgolib{} test datasets did not reveal a clear dominant algorithm.
	On the contrary, it has revealed rather conflicting results.
	While deterministic (especially \direct-type) algorithms performed better in solving \gkls{}-generated test problems, various stochastic algorithms dominated and performed better in solving \directgolib{} test problems.
	Therefore, to summarize all experimental results, we calculated the average success rates of all algorithms in \Cref{tab:results4}.
	Here, the second and fifth columns show the average success rates on all \gkls{} problems using two different budgets, while the third and sixth columns show on \directgolib{} test problems.
	The fourth and seventh columns present the combined average success rates on all $1,197$ test problems for each budget separately.
	Finally, the eighth column gives the overall average used to rank all considered algorithms.
	
	\begin{table}
		\centering
		\caption{The average number of success rates using different test sets and evaluation budgets.}
			\fontsize{7.5pt}{7.5pt}\selectfont
			\begin{tabular*}{\textwidth}{@{\extracolsep{\fill}}l|ccc|ccc|c}
					\toprule
					Evaluation budget & \multicolumn{3}{c|}{$M_{\rm max}=2,500$} & \multicolumn{3}{c|}{$M_{\rm max}=500,000$}  \\
					\cmidrule{0-6}
					Test set & \rotatebox{70}{\gkls} & \rotatebox{70}{\texttt{DIRECTGOLib}} & \rotatebox{70}{Average} & \rotatebox{70}{\gkls} & \rotatebox{70}{\texttt{DIRECTGOLib}} & \rotatebox{70}{Average} & \rotatebox{70}{Overall average} \\
					\midrule
					\birmin	&	$0.64$	&	$0.46$	&	$0.55$	&	$1.00$	&	$0.57$	&	$0.79$	&	$0.67$	\\
					\oqnlp	&	$0.44$	&	$0.54$	&	$0.49$	&	$0.93$	&	$0.73$	&	$0.83$	&	$0.66$	\\
					\mcs	&	$0.44$	&	$0.53$	&	$0.48$	&	$0.99$	&	$0.68$	&	$0.84$	&	$0.66$	\\
					\mqnlp	&	$0.50$	&	$0.50$	&	$0.50$	&	$0.94$	&	$0.64$	&	$0.79$	&	$0.65$	\\
					\directrev	&	$0.54$	&	$0.45$	&	$0.50$	&	$1.00$	&	$0.57$	&	$0.78$	&	$0.64$	\\
					\dbdvg	&	$0.53$	&	$0.33$	&	$0.43$	&	$0.99$	&	$0.66$	&	$0.83$	&	$0.63$	\\
					\dbdvo	&	$0.55$	&	$0.39$	&	$0.47$	&	$1.00$	&	$0.56$	&	$0.78$	&	$0.63$	\\
					\dtcg	&	$0.46$	&	$0.37$	&	$0.41$	&	$0.99$	&	$0.66$	&	$0.82$	&	$0.62$	\\
					\dtco	&	$0.54$	&	$0.38$	&	$0.46$	&	$1.00$	&	$0.53$	&	$0.77$	&	$0.61$	\\
					\avoa	&	$0.26$	&	$0.62$	&	$0.44$	&	$0.74$	&	$0.80$	&	$0.77$	&	$0.60$	\\
					\glbsolve	&	$0.55$	&	$0.30$	&	$0.43$	&	$1.00$	&	$0.56$	&	$0.78$	&	$0.60$	\\
					\dbdva	&	$0.45$	&	$0.34$	&	$0.40$	&	$0.99$	&	$0.62$	&	$0.81$	&	$0.60$	\\
					\birectgb	&	$0.62$	&	$0.31$	&	$0.46$	&	$1.00$	&	$0.45$	&	$0.73$	&	$0.59$	\\
					\directgl	&	$0.49$	&	$0.31$	&	$0.40$	&	$0.98$	&	$0.59$	&	$0.78$	&	$0.59$	\\
					\directg	&	$0.50$	&	$0.31$	&	$0.40$	&	$0.99$	&	$0.55$	&	$0.77$	&	$0.59$	\\
					\birect	&	$0.55$	&	$0.32$	&	$0.44$	&	$1.00$	&	$0.45$	&	$0.73$	&	$0.58$	\\
					\dtca	&	$0.38$	&	$0.35$	&	$0.37$	&	$0.97$	&	$0.62$	&	$0.80$	&	$0.58$	\\
					\glccluster	&	$0.48$	&	$0.25$	&	$0.36$	&	$1.00$	&	$0.59$	&	$0.79$	&	$0.58$	\\
					\dirmin	&	$0.28$	&	$0.43$	&	$0.36$	&	$0.97$	&	$0.61$	&	$0.79$	&	$0.57$	\\
					\scso	&	$0.32$	&	$0.52$	&	$0.42$	&	$0.81$	&	$0.63$	&	$0.72$	&	$0.57$	\\
					\direct	&	$0.55$	&	$0.28$	&	$0.42$	&	$1.00$	&	$0.44$	&	$0.72$	&	$0.57$	\\
					\directm	&	$0.54$	&	$0.28$	&	$0.41$	&	$1.00$	&	$0.44$	&	$0.72$	&	$0.56$	\\
					\directa	&	$0.56$	&	$0.28$	&	$0.42$	&	$1.00$	&	$0.41$	&	$0.70$	&	$0.56$	\\
					\directmro	&	$0.48$	&	$0.29$	&	$0.38$	&	$0.99$	&	$0.46$	&	$0.72$	&	$0.55$	\\
					\multistart	&	$0.26$	&	$0.40$	&	$0.33$	&	$0.97$	&	$0.57$	&	$0.77$	&	$0.55$	\\
					\directl	&	$0.48$	&	$0.31$	&	$0.39$	&	$0.93$	&	$0.47$	&	$0.70$	&	$0.55$	\\
					\plor	&	$0.44$	&	$0.34$	&	$0.39$	&	$0.83$	&	$0.55$	&	$0.69$	&	$0.54$	\\
					\lgo	&	$0.10$	&	$0.44$	&	$0.27$	&	$0.98$	&	$0.62$	&	$0.80$	&	$0.54$	\\
					\directrest	&	$0.58$	&	$0.30$	&	$0.44$	&	$0.82$	&	$0.42$	&	$0.62$	&	$0.53$	\\
					\adc	&	$0.45$	&	$0.29$	&	$0.37$	&	$1.00$	&	$0.35$	&	$0.68$	&	$0.52$	\\
					\sta	&	$0.23$	&	$0.27$	&	$0.25$	&	$0.78$	&	$0.69$	&	$0.74$	&	$0.49$	\\
					\disimplv	&	$0.44$	&	$0.26$	&	$0.35$	&	$0.90$	&	$0.36$	&	$0.63$	&	$0.49$	\\
					\aos	&	$0.22$	&	$0.45$	&	$0.33$	&	$0.70$	&	$0.57$	&	$0.63$	&	$0.48$	\\
					\directaggress	&	$0.29$	&	$0.26$	&	$0.28$	&	$0.89$	&	$0.48$	&	$0.68$	&	$0.48$	\\
					\sa	&	$0.27$	&	$0.23$	&	$0.25$	&	$0.99$	&	$0.36$	&	$0.67$	&	$0.46$	\\
					\StochasticRBF	&	$0.22$	&	$0.30$	&	$0.26$	&	$0.96$	&	$0.37$	&	$0.67$	&	$0.46$	\\
					\disimplc	&	$0.38$	&	$0.22$	&	$0.30$	&	$0.89$	&	$0.32$	&	$0.61$	&	$0.45$	\\
					\surogate	&	$0.48$	&	$0.35$	&	$0.41$	&	$0.61$	&	$0.37$	&	$0.49$	&	$0.45$	\\
					\globalsearch	&	$0.00$	&	$0.11$	&	$0.06$	&	$0.98$	&	$0.65$	&	$0.81$	&	$0.43$	\\
					\cnna	&	$0.11$	&	$0.19$	&	$0.15$	&	$0.76$	&	$0.68$	&	$0.72$	&	$0.43$	\\
					\be	&	$0.32$	&	$0.24$	&	$0.28$	&	$0.68$	&	$0.49$	&	$0.58$	&	$0.43$	\\
					\beso	&	$0.16$	&	$0.49$	&	$0.33$	&	$0.40$	&	$0.67$	&	$0.53$	&	$0.43$	\\
					\nna	&	$0.16$	&	$0.14$	&	$0.15$	&	$0.71$	&	$0.67$	&	$0.69$	&	$0.42$	\\
					\multimin	&	$0.11$	&	$0.07$	&	$0.09$	&	$0.93$	&	$0.53$	&	$0.73$	&	$0.41$	\\
					\amc	&	$0.40$	&	$0.26$	&	$0.33$	&	$0.52$	&	$0.44$	&	$0.48$	&	$0.40$	\\
					\hs	&	$0.22$	&	$0.24$	&	$0.23$	&	$0.63$	&	$0.50$	&	$0.57$	&	$0.40$	\\
					\pa	&	$0.33$	&	$0.29$	&	$0.31$	&	$0.36$	&	$0.56$	&	$0.46$	&	$0.38$	\\
					\cs	&	$0.10$	&	$0.18$	&	$0.14$	&	$0.55$	&	$0.62$	&	$0.58$	&	$0.36$	\\
					\ica	&	$0.27$	&	$0.28$	&	$0.28$	&	$0.34$	&	$0.54$	&	$0.44$	&	$0.36$	\\
					\fstde	&	$0.16$	&	$0.38$	&	$0.27$	&	$0.24$	&	$0.65$	&	$0.45$	&	$0.36$	\\
					\sfla	&	$0.25$	&	$0.35$	&	$0.30$	&	$0.25$	&	$0.52$	&	$0.39$	&	$0.34$	\\
					\tlbo	&	$0.15$	&	$0.32$	&	$0.23$	&	$0.33$	&	$0.56$	&	$0.45$	&	$0.34$	\\
					\bads	&	$0.08$	&	$0.44$	&	$0.26$	&	$0.40$	&	$0.45$	&	$0.42$	&	$0.34$	\\
					\abc	&	$0.07$	&	$0.17$	&	$0.12$	&	$0.59$	&	$0.49$	&	$0.54$	&	$0.33$	\\
					\ga	&	$0.21$	&	$0.27$	&	$0.24$	&	$0.22$	&	$0.54$	&	$0.38$	&	$0.31$	\\
					\aco	&	$0.09$	&	$0.35$	&	$0.22$	&	$0.13$	&	$0.54$	&	$0.33$	&	$0.28$	\\
					\de	&	$0.10$	&	$0.35$	&	$0.22$	&	$0.10$	&	$0.53$	&	$0.31$	&	$0.27$	\\
					\bbo	&	$0.14$	&	$0.31$	&	$0.23$	&	$0.14$	&	$0.47$	&	$0.31$	&	$0.27$	\\
					\fa	&	$0.14$	&	$0.28$	&	$0.21$	&	$0.14$	&	$0.49$	&	$0.31$	&	$0.26$	\\
					\ca	&	$0.02$	&	$0.48$	&	$0.25$	&	$0.02$	&	$0.52$	&	$0.27$	&	$0.26$	\\
					\iwo	&	$0.15$	&	$0.22$	&	$0.19$	&	$0.16$	&	$0.45$	&	$0.30$	&	$0.24$	\\
					\cmaes	&	$0.02$	&	$0.27$	&	$0.15$	&	$0.02$	&	$0.45$	&	$0.24$	&	$0.19$	\\
					\mc	&	$0.02$	&	$0.10$	&	$0.06$	&	$0.30$	&	$0.22$	&	$0.26$	&	$0.16$	\\
					\bat	&	$0.06$	&	$0.25$	&	$0.16$	&	$0.06$	&	$0.26$	&	$0.16$	&	$0.16$	\\
					
					\bottomrule
				\end{tabular*}
			\label{tab:results4}
		\end{table}
		
		According to the results presented in column fourth, the best algorithm for expensive global optimization $(M_{\rm max}=2,500)$ is \birmin, which solved $55 \%$ of test problems.
		\birmin{} was followed closely by \mqnlp, \directrev, \oqnlp, and \mcs{}, which solved $50 \%$, $50 \%$, $49 \%$, and $48 \%$ of test problems respectively.
		All five algorithms are hybridized in nature.
		The best performing non-hybridized algorithm is \dbdvo, closely followed by \dtco{} and \birectgb.
		We note that even $38$ algorithms solved less than $40 \%$.
		Among stochastic, only three solvers (\avoa, \scso, and \surogate{}) achieved a success ratio greater than $0.40$.
		The best from them (\avoa{} algorithm) shares only $9$--$11$ places in the overall ranking (together with \birect{} and \glbsolve{}).
		
		The seventh column in \Cref{tab:results4} presents the average success rate on all $1,197$ test problems using $M_{\rm max}=500,000$.
		In this case, the best performing algorithm is hybridized \mcs, which solved $84 \%$ of all test problems.
		The \oqnlp{} and non-hybridized version of \dbdvg{} achieved almost the same success rate, and both solved $83 \%$ of test problems.
		At least seventeen other algorithms solved more than three-quarters ($75 \%$) of problems.
		Among the top-20, only four belong to the class of stochastic algorithms.
		The \globalsearch{}, which solved $81 \%$, had the best success rate.
		
		The last column in \Cref{tab:results4} presents the combined overall average success rate.
		The best overall performance has been observed using the \birmin{} algorithm, closely followed by \mcs, \oqnlp, \mqnlp, and \directrev.
		Quite surprisingly, all five algorithms are deterministic and hybridized with local searches.
		The best non-hybridized \dbdvg{} and \dbdvo{} algorithms solved only $4 \%$ fewer problems than the best performing \birmin{} algorithm and, therefore, are very competitive with hybridized counterparts.
		The best stochastic solver was \avoa, which ranks only tenth in the overall rankings.
		However, in the name of fairness, we want to stress the frequent leadership of stochastic solvers in solving difficult high-dimensional problems.
		Such problems are very often encountered in practice.
		Also, while the often weaker performance of meta-heuristic solvers was observed based on the function evaluation criteria, they were often much faster (concerning execution time) than deterministic ones.
		Therefore, they can perform more function evaluations in the same time budget.
		
		Finally, a study a decade earlier \cite{Rios2013} using the same $(M_{\rm max}=2,500)$ concluded that the best-performing algorithms were \multimin, \glccluster, \mcs, and \lgo.
		Unexpectedly, only the \mcs{} algorithm remained among the best-performing in this study.
		In contrast, the \glccluster{} is on $28$th, while \lgo{} is only on $41$st.
		Therefore, these and similar experimental results should not be abstracted, but the best solver should be chosen according to the problem specifics.

\section{Conclusions and possible extensions}\label{sec:conclusion}

Sixty-four state-of-the-art stochastic and deterministic derivative-free optimization algorithms were tested on two large publicly available test sets (\gkls{} and \directgolib).
The computational results showed no single algorithm whose performance dominates all experiments.

For most stochastic algorithms, even low-dimensional \gkls{} test problems caused difficulties in finding solutions, even when the budget of function evaluation was large ($M_{\rm max}=500,000$).
Unlike stochastic algorithms, deterministic ones (especially \direct-type) showed superior efficiency on \gkls{} test problems.
The \birmin{} and \birectgb{} algorithms were the most effective using the low $(M_{\rm max}=2,500)$ and extensive ($M_{\rm max}=500,000$) budget of function evaluations.
Moreover, hybridized versions often cause excessive local search usage and worsen the performance compared to standard versions.
Additionally, hybridized version did not guarantee the best result even when $M_{\rm max}=500,000$.

Experimentation on the \directgolib{} test set revealed that stochastic algorithms are very competitive and often superior, especially in solving high-dimensional test problems.
When $M_{\rm max}=2,500$, the meta-heuristic-based \avoa, \scso, \beso, and \ca{} are the best-performing algorithms.
However, the performance of hybridized deterministic algorithms  (\oqnlp, \mcs, \mqnlp, and \birmin) is similar.
When $M_{\rm max}=500,000$, \avoa, \oqnlp, and \mcs{} remain the best-performing ones.
Moreover, this experimental part confirmed the well-known fact that \direct-type solvers are very effective for low-dimensional problems, but efficiency significantly decreases with dimensionality.


The summary of the entire set of test problems ($1,197$ in total) showed that, on average, the best-performing algorithms are deterministic hybrid versions: \birmin, \oqnlp, \mcs, \mqnlp, and \directrev.
In sixth to ninth places, the standard non-hybridized \direct-type algorithms, \dbdvg, \dbdvo, \dtcg, and \dtco{}, were introduced recently \cite{Stripinis2021b}.
The best stochastic algorithm, \avoa, ranks only tenth in this study.
Although low-dimensional, \gkls{} test problems were very unfavorable for stochastic solvers and had a significant impact on generalized final rankings.

Finally, this paper has dealt with box-constrained global optimization problems.
Therefore, it could be extended to a derivative-free constrained case.
At the same time, the family of existing stochastic algorithms is much larger than that of deterministic ones.
Therefore, including new ones could lead to discovering more competitive algorithms than those included in this study.
The results also show that the performance of particular algorithms can be highly dependent on the test functions.
Thus, the more extensive and diverse test problems are, the more reliable the results are.
Therefore, increasing the number and diversity of test problems used in benchmarks must be continuous.


\section*{Data statement}
\textbf{DIRECTGOLib} - \textbf{DIRECT} \textbf{G}lobal \textbf{O}ptimization test problems \textbf{Lib}rary is designed as a continuously-growing open-source GitHub repository to which anyone can easily contribute.
The exact data underlying this article from \directgolib{} can be accessed either on GitHub or at Zenodo:
\begin{itemize}
	\item \url{https://github.com/blockchain-group/DIRECTGOLib},
	\item \url{https://zenodo.org/record/6617799},
\end{itemize}
and used under the MIT license.
We welcome contributions and corrections to it.

\section*{Code availability}
All $64$ solvers used in this study are implemented in \matlab{} and references are given in \Cref{tab:direct,tab:EAs}.
Most solvers are open-source and freely available to everyone, except several solvers within the \tomlab{} toolbox.


\bibliographystyle{spmpsci}      
\bibliography{library}   

\begin{thebibliography}{100}
\providecommand{\url}[1]{{#1}}
\providecommand{\urlprefix}{URL }
\expandafter\ifx\csname urlstyle\endcsname\relax
  \providecommand{\doi}[1]{DOI~\discretionary{}{}{}#1}\else
  \providecommand{\doi}{DOI~\discretionary{}{}{}\begingroup
  \urlstyle{rm}\Url}\fi

\bibitem{Layeb2022}
Abdesslem, L.: New hard benchmark functions for global optimization (2022).
\newblock \urlprefix\url{https://www.mathworks.com/matlabcentral}.
\newblock MATLAB Central File Exchange. Retrieved February 18, 2022.

\bibitem{Abdollahzadeh2021}
Abdollahzadeh, B., Gharehchopogh, F.S., Mirjalili, S.: African vultures
  optimization algorithm: A new nature-inspired metaheuristic algorithm for
  global optimization problems.
\newblock Computers \& Industrial Engineering \textbf{158}, 107408 (2021).
\newblock \doi{https://doi.org/10.1016/j.cie.2021.107408}.
\newblock
  \urlprefix\url{https://www.sciencedirect.com/science/article/pii/S0360835221003120}

\bibitem{Acerbi2022}
Acerbi, L., Ma, W.J.: Bayesian adaptive direct search (bads) - v1.0.8 (2012).
\newblock
  \urlprefix\url{https://se.mathworks.com/matlabcentral/fileexchange/89162-bayesian-adaptive-direct-search-bads-optimizer?s_tid=srchtitle}.
\newblock MATLAB Central File Exchange. Retrieved June 22, 2022.

\bibitem{Acerbi2017}
Acerbi, L., Ma, W.J.: Practical bayesian optimization for model fitting with
  bayesian adaptive direct search (2017)

\bibitem{Hassan2022}
ALSattar, H.: Bald eagle search optimization algorithm (bes) (2022).
\newblock
  \urlprefix\url{https://www.mathworks.com/matlabcentral/fileexchange/86862-bald-eagle-search-optimization-algorithm-bes}.
\newblock MATLAB Central File Exchange. Retrieved February 18, 2022.

\bibitem{Alsattar2020}
Alsattar, H.A., Zaidan, A.A., Zaidan, B.B.: {Novel meta-heuristic bald eagle
  search optimisation algorithm}.
\newblock Artificial Intelligence Review \textbf{53}(3), 2237--2264 (2020).
\newblock \doi{10.1007/s10462-019-09732-5}.
\newblock \urlprefix\url{https://doi.org/10.1007/s10462-019-09732-5}

\bibitem{Gargari2007}
Atashpaz-Gargari, E., Lucas, C.: Imperialist competitive algorithm: An
  algorithm for optimization inspired by imperialistic competition.
\newblock In: 2007 IEEE Congress on Evolutionary Computation, pp. 4661--4667
  (2007).
\newblock \doi{10.1109/CEC.2007.4425083}

\bibitem{Auger2009}
Auger, A., Hansen, N., Perez~Zerpa, J.M., Ros, R., Schoenauer, M.: Experimental
  comparisons of derivative free optimization algorithms.
\newblock In: J.~Vahrenhold (ed.) Experimental Algorithms, pp. 3--15. Springer
  Berlin Heidelberg, Berlin, Heidelberg (2009)

\bibitem{Azizi2021}
Azizi, M.: Atomic orbital search: A novel metaheuristic algorithm.
\newblock Applied Mathematical Modelling \textbf{93}, 657--683 (2021).
\newblock \doi{https://doi.org/10.1016/j.apm.2020.12.021}.
\newblock
  \urlprefix\url{https://www.sciencedirect.com/science/article/pii/S0307904X20307198}

\bibitem{Azizi2022}
Azizi, M.: Matlab code for atomic orbital search: A novel metaheuristic (2022).
\newblock
  \urlprefix\url{https://www.mathworks.com/matlabcentral/fileexchange/91160-matlab-code-for-atomic-orbital-search-a-novel-metaheuristic}.
\newblock MATLAB Central File Exchange. Retrieved February 18, 2022.

\bibitem{Baker2000}
Baker, C.A., Watson, L.T., Grossman, B., Mason, W.H., Haftka, R.T.: Parallel
  global aircraft configuration design space exploration.
\newblock In: A.~Tentner (ed.) High Performance Computing Symposium 2000, pp.
  54--66. Soc. for Computer Simulation Internat (2000)

\bibitem{Bartholomew2002}
Bartholomew-Biggs, M.C., Parkhurst, S.C., Wilson, S.P.: Using {DIRECT} to solve
  an aircraft routing problem.
\newblock Computational Optimization and Applications \textbf{21}(3), 311--323
  (2002).
\newblock \doi{10.1023/A:1013729320435}

\bibitem{Bjorkman2000}
Bj{\"{o}}rkman, M., Holmstr{\"{o}}m, K.: {Global Optimization of Costly
  Nonconvex Functions Using Radial Basis Functions}.
\newblock Optimization and Engineering \textbf{1}(4), 373--397 (2000).
\newblock \doi{10.1023/A:1011584207202}.
\newblock \urlprefix\url{https://doi.org/10.1023/A:1011584207202}

\bibitem{BLUM20114135}
Blum, C., Puchinger, J., Raidl, G.R., Roli, A.: Hybrid metaheuristics in
  combinatorial optimization: A survey.
\newblock Applied Soft Computing \textbf{11}(6), 4135--4151 (2011).
\newblock \doi{https://doi.org/10.1016/j.asoc.2011.02.032}.
\newblock
  \urlprefix\url{https://www.sciencedirect.com/science/article/pii/S1568494611000962}

\bibitem{Blum2003}
Blum, C., Roli, A.: Metaheuristics in combinatorial optimization: Overview and
  conceptual comparison.
\newblock ACM Comput. Surv. \textbf{35}(3), 268–308 (2003).
\newblock \doi{10.1145/937503.937505}.
\newblock \urlprefix\url{https://doi.org/10.1145/937503.937505}

\bibitem{BOUSSAID201382}
Boussaïd, I., Lepagnot, J., Siarry, P.: A survey on optimization
  metaheuristics.
\newblock Information Sciences \textbf{237}, 82--117 (2013).
\newblock \doi{https://doi.org/10.1016/j.ins.2013.02.041}.
\newblock
  \urlprefix\url{https://www.sciencedirect.com/science/article/pii/S0020025513001588}.
\newblock Prediction, Control and Diagnosis using Advanced Neural Computations

\bibitem{Byrd2000}
Byrd, R.H., Gilbert, J.C., Nocedal, J.: {A trust region method based on
  interior point techniques for nonlinear programming}.
\newblock Mathematical Programming \textbf{89}(1), 149--185 (2000).
\newblock \doi{10.1007/PL00011391}.
\newblock \urlprefix\url{https://doi.org/10.1007/PL00011391}

\bibitem{Carter2001}
Carter, R.G., Gablonsky, J.M., Patrick, A., Kelley, C.T., Eslinger, O.J.:
  Algorithms for noisy problems in gas transmission pipeline optimization.
\newblock Optimization and Engineering \textbf{2}(2), 139--157 (2001).
\newblock \doi{10.1023/A:1013123110266}

\bibitem{Clerc1999}
Clerc, M.: {The swarm and the queen: Towards a deterministic and adaptive
  particle swarm optimization}.
\newblock In: Proceedings of the 1999 Congress on Evolutionary Computation, CEC
  1999 (1999).
\newblock \doi{10.1109/CEC.1999.785513}

\bibitem{Cox2001}
Cox, S.E., Haftka, R.T., Baker, C.A., Grossman, B., Mason, W.H., Watson, L.T.:
  A comparison of global optimization methods for the design of a high-speed
  civil transport.
\newblock Journal of Global Optimization \textbf{21}(4), 415--432 (2001).
\newblock \doi{10.1023/A:1012782825166}

\bibitem{Serafino2011}
{Di Serafino}, D., Liuzzi, G., Piccialli, V., Riccio, F., Toraldo, G.: A
  modified {DI}viding {RECT}angles algorithm for a problem in astrophysics.
\newblock Journal of Optimization Theory and Applications \textbf{151}(1),
  175--190 (2011).
\newblock \doi{10.1007/s10957-011-9856-9}

\bibitem{Dixon1978}
Dixon, L., Szeg\"{o}, C.: The global optimisation problem: An introduction.
\newblock In: L.~Dixon, G.~Szeg\"{o} (eds.) Towards Global Optimization,
  vol.~2, pp. 1--15. North-Holland Publishing Company (1978)

\bibitem{Dorigo1992O}
Dorigo, M., Maniezzo, V., Colorni, A.: Ant system: optimization by a colony of
  cooperating agents.
\newblock IEEE transactions on systems, man, and cybernetics. Part B,
  Cybernetics : a publication of the IEEE Systems, Man, and Cybernetics Society
  \textbf{26 1}, 29--41 (1996)

\bibitem{Eberhart1995}
Eberhart, R., Kennedy, J.: A new optimizer using particle swarm theory.
\newblock In: MHS'95. Proceedings of the Sixth International Symposium on Micro
  Machine and Human Science, pp. 39--43 (1995).
\newblock \doi{10.1109/MHS.1995.494215}

\bibitem{Muzaffar2003}
Eusuff, M.M., Lansey, K.E.: Optimization of water distribution network design
  using the shuffled frog leaping algorithm.
\newblock Journal of Water Resources Planning and Management \textbf{129}(3),
  210--225 (2003).
\newblock \doi{10.1061/(ASCE)0733-9496(2003)129:3(210)}

\bibitem{Ezugwu2021}
Ezugwu, A.E., Shukla, A.K., Nath, R., Akinyelu, A.A., Agushaka, J.O., Chiroma,
  H., Muhuri, P.K.: {Metaheuristics: a comprehensive overview and
  classification along with bibliometric analysis}.
\newblock Artificial Intelligence Review \textbf{54}(6), 4237--4316 (2021).
\newblock \doi{10.1007/s10462-020-09952-0}.
\newblock \urlprefix\url{https://doi.org/10.1007/s10462-020-09952-0}

\bibitem{Finkel2004aa}
Finkel, D.E., Kelley, C.T.: An adaptive restart implementation of direct.
\newblock Technical report CRSC-TR04-30, Center for Research in Scientific
  Computation, North Carolina State University, Raleigh  (2004)

\bibitem{Finkel2006}
Finkel, D.E., Kelley, C.T.: Additive scaling and the {DIRECT} algorithm.
\newblock Journal of Global Optimization \textbf{36}(4), 597--608 (2006).
\newblock \doi{10.1007/s10898-006-9029-9}

\bibitem{Floudas1999book}
Floudas, C.A.: Deterministic global optimization: theory, methods and
  applications, \emph{Nonconvex Optimization and Its Applications}, vol.~37.
\newblock Springer US (1999).
\newblock \doi{10.1007/978-1-4757-4949-6}

\bibitem{Gablonsky2001:phd}
Gablonsky, J.M.: Modifications of the {DIRECT} algorithm.
\newblock Ph.D. thesis, North Carolina State University (2001)

\bibitem{Gablonsky2001}
Gablonsky, J.M., Kelley, C.T.: A locally-biased form of the {DIRECT} algorithm.
\newblock Journal of Global Optimization \textbf{21}(1), 27--37 (2001).
\newblock \doi{10.1023/A:1017930332101}

\bibitem{Gavana2021}
Gavana, A.: Global optimization benchmarks and ampgo.
\newblock \url{http://infinity77.net/global_optimization/index.html}.
\newblock Online; accessed: 2021-07-22

\bibitem{Gaviano2003}
Gaviano, M., Kvasov, D.E., Lera, D., Sergeyev, Y.D.: Algorithm 829: {Software}
  for generation of classes of test functions with known local and global
  minima for global optimization.
\newblock ACM Transactions on Mathematical Software (TOMS) \textbf{29}(4),
  469--480 (2003).
\newblock \doi{10.1145/962437.962444}

\bibitem{Geem2001}
Geem, Z.W., Kim, J.H., Loganathan, G.: A new heuristic optimization algorithm:
  Harmony search.
\newblock SIMULATION \textbf{76}(2), 60--68 (2001).
\newblock \doi{10.1177/003754970107600201}.
\newblock \urlprefix\url{https://doi.org/10.1177/003754970107600201}

\bibitem{GERGEL1999163}
Gergel, V., Sergeyev, Y.: Sequential and parallel algorithms for global
  minimizing functions with lipschitzian derivatives.
\newblock Computers and Mathematics with Applications \textbf{37}(4), 163--179
  (1999).
\newblock \doi{https://doi.org/10.1016/S0898-1221(99)00067-X}.
\newblock
  \urlprefix\url{https://www.sciencedirect.com/science/article/pii/S089812219900067X}

\bibitem{GERGEL1997}
GERGEL, V.P.: {A Global Optimization Algorithm for Multivariate Functions with
  Lipschitzian First Derivatives}.
\newblock Journal of Global Optimization \textbf{10}(3), 257--281 (1997).
\newblock \doi{10.1023/A:1008290629896}.
\newblock \urlprefix\url{https://doi.org/10.1023/A:1008290629896}

\bibitem{Gharehchopogh2019a}
Gharehchopogh, F.S., Gholizadeh, H.: A comprehensive survey: Whale optimization
  algorithm and its applications.
\newblock Swarm and Evolutionary Computation \textbf{48}, 1--24 (2019).
\newblock \doi{https://doi.org/10.1016/j.swevo.2019.03.004}.
\newblock
  \urlprefix\url{https://www.sciencedirect.com/science/article/pii/S2210650218309350}

\bibitem{Gharehchopogh2020}
Gharehchopogh, F.S., Shayanfar, H., Gholizadeh, H.: {A comprehensive survey on
  symbiotic organisms search algorithms}.
\newblock Artificial Intelligence Review \textbf{53}(3), 2265--2312 (2020).
\newblock \doi{10.1007/s10462-019-09733-4}.
\newblock \urlprefix\url{https://doi.org/10.1007/s10462-019-09733-4}

\bibitem{GLOVER1986533}
Glover, F.: Future paths for integer programming and links to artificial
  intelligence.
\newblock Computers and Operations Research \textbf{13}(5), 533--549 (1986).
\newblock \doi{https://doi.org/10.1016/0305-0548(86)90048-1}.
\newblock
  \urlprefix\url{https://www.sciencedirect.com/science/article/pii/0305054886900481}.
\newblock Applications of Integer Programming

\bibitem{Gutmann2001}
Gutmann, H.M.: {A Radial Basis Function Method for Global Optimization}.
\newblock Journal of Global Optimization \textbf{19}(3), 201--227 (2001).
\newblock \doi{10.1023/A:1011255519438}.
\newblock \urlprefix\url{https://doi.org/10.1023/A:1011255519438}

\bibitem{Hansen2010}
Hansen, N., Auger, A., Ros, R., Finck, S., Po\v{s}\'{\i}k, P.: Comparing
  results of 31 algorithms from the black-box optimization benchmarking
  bbob-2009.
\newblock In: Proceedings of the 12th Annual Conference Companion on Genetic
  and Evolutionary Computation, GECCO '10, p. 1689–1696. Association for
  Computing Machinery, New York, NY, USA (2010).
\newblock \doi{10.1145/1830761.1830790}.
\newblock \urlprefix\url{https://doi.org/10.1145/1830761.1830790}

\bibitem{Hedar2005}
Hedar, A.: Test functions for unconstrained global optimization.
\newblock
  \url{http://www-optima.amp.i.kyoto-u.ac.jp/member/student/hedar/Hedar_files/TestGO.htm}
  (2005).
\newblock Online; accessed: 2017-03-22

\bibitem{Kalami2020}
Heris, M.K.: Yarpiz evolutionary algorithms toolbox for matlab (ypea) (2020).
\newblock
  \urlprefix\url{https://se.mathworks.com/matlabcentral/fileexchange/71373-ypea?s_tid=srchtitle}.
\newblock MATLAB Central File Exchange. Retrieved February 18, 2022.

\bibitem{Hey1979}
Hey, A.M.: {Towards Global Optimisation 2}.
\newblock Journal of the Operational Research Society \textbf{30}(9), 844
  (1979).
\newblock \doi{10.1057/jors.1979.201}.
\newblock \urlprefix\url{https://doi.org/10.1057/jors.1979.201}

\bibitem{Holland1975}
Holland, J.: Adaptation in Natural and Artificial Systems.
\newblock The University of Michigan Press, Ann Arbor (1975)

\bibitem{Holmstrom2004}
Holmstr{\"o}m, K., Edvall, M.M.: The TOMLAB Optimization Environment, pp.
  369--376.
\newblock Springer US, Boston, MA (2004).
\newblock \doi{10.1007/978-1-4613-0215-5_19}.
\newblock \urlprefix\url{https://doi.org/10.1007/978-1-4613-0215-5_19}

\bibitem{Holmstrom2010}
Holmstrom, K., Goran, A.O., Edvall, M.M.: User’s guide for tomlab 7 (2010).
\newblock \urlprefix\url{https://tomopt.com/}

\bibitem{Hooke1961}
Hooke, R., Jeeves, T.A.: `` direct search'' solution of numerical and
  statistical problems.
\newblock J. ACM \textbf{8}(2), 212–229 (1961).
\newblock \doi{10.1145/321062.321069}.
\newblock \urlprefix\url{https://doi.org/10.1145/321062.321069}

\bibitem{Horst1996:book}
Horst, R., Tuy, H.: Global Optimization: Deterministic Approaches.
\newblock Springer, Berlin (1996)

\bibitem{Huyer1999}
Huyer, W., Neumaier, A.: {Global Optimization by Multilevel Coordinate Search}.
\newblock Journal of Global Optimization \textbf{14}(4), 331--355 (1999).
\newblock \doi{10.1023/A:1008382309369}.
\newblock \urlprefix\url{https://doi.org/10.1023/A:1008382309369}

\bibitem{Iruthayarajan2010}
Iruthayarajan, M.W., Baskar, S.: Covariance matrix adaptation evolution
  strategy based design of centralized pid controller.
\newblock Expert systems with Applications \textbf{37}(8), 5775--5781 (2010)

\bibitem{Jamil2013}
Jamil, M., Yang, X.S.: A literature survey of benchmark functions for global
  optimisation problems.
\newblock International Journal of Mathematical Modelling and Numerical
  Optimisation \textbf{4}(2), 150--194 (2013).
\newblock \doi{10.1504/IJMMNO.2013.055204}.
\newblock
  \urlprefix\url{https://www.inderscienceonline.com/doi/abs/10.1504/IJMMNO.2013.055204}.
\newblock PMID: 55204

\bibitem{JangaReddy2020}
{Janga Reddy}, M., {Nagesh Kumar}, D.: {Evolutionary algorithms, swarm
  intelligence methods, and their applications in water resources engineering:
  a state-of-the-art review}.
\newblock H2Open Journal \textbf{3}(1), 135--188 (2020).
\newblock \doi{10.2166/h2oj.2020.128}.
\newblock \urlprefix\url{https://doi.org/10.2166/h2oj.2020.128}

\bibitem{Jones2001}
Jones, D.R.: The {{\sc Direct}} global optimization algorithm.
\newblock In: C.A. Floudas, P.M. Pardalos (eds.) The Encyclopedia of
  Optimization, pp. 431--440. Kluwer Academic Publishers, Dordrect (2001)

\bibitem{Jones2021}
Jones, D.R., Martins, J.R.R.A.: {The DIRECT algorithm: 25 years Later}.
\newblock Journal of Global Optimization \textbf{79}(3), 521--566 (2021).
\newblock \doi{10.1007/s10898-020-00952-6}.
\newblock \urlprefix\url{https://doi.org/10.1007/s10898-020-00952-6}

\bibitem{Jones1993}
Jones, D.R., Perttunen, C.D., Stuckman, B.E.: Lipschitzian optimization without
  the {L}ipschitz constant.
\newblock Journal of Optimization Theory and Application \textbf{79}(1),
  157--181 (1993).
\newblock \doi{10.1007/BF00941892}

\bibitem{Jones1998}
Jones, D.R., Schonlau, M., Welch, W.J.: {Efficient Global Optimization of
  Expensive Black-Box Functions}.
\newblock Journal of Global Optimization \textbf{13}(4), 455--492 (1998).
\newblock \doi{10.1023/A:1008306431147}.
\newblock \urlprefix\url{https://doi.org/10.1023/A:1008306431147}

\bibitem{Julie2006}
Julie: Stochastic radial basis function algorithm for global optimization
  (2022).
\newblock
  \urlprefix\url{https://www.mathworks.com/matlabcentral/fileexchange/42090-stochastic-radial-basis-function-algorithm-for-global-optimization}.
\newblock MATLAB Central File Exchange. Retrieved July 20, 2022.

\bibitem{Mishra2006}
K., M.S.: Some new test functions for global optimization and performance of
  repulsive particle swarm method (2006).
\newblock \doi{http://dx.doi.org/10.2139/ssrn.926132}.
\newblock \urlprefix\url{https://ssrn.com/abstract=926132}.
\newblock MATLAB Central File Exchange. Retrieved February 18, 2022.

\bibitem{Karaboga2005ANIB}
Karaboga, D.: An idea based on honey bee swarm for numerical optimization
  (2005)

\bibitem{Kirkpatrick1983a}
Kirkpatrick, S., Gelatt, C.D., Vecchi, M.P.: Optimization by simulated
  annealing.
\newblock Science \textbf{220}(4598), 671--680 (1983).
\newblock \doi{10.1126/science.220.4598.671}.
\newblock
  \urlprefix\url{https://www.science.org/doi/abs/10.1126/science.220.4598.671}

\bibitem{Kvasov2018}
Kvasov, D.E., Mukhametzhanov, M.S.: Metaheuristic vs. deterministic global
  optimization algorithms: The univariate case.
\newblock Applied Mathematics and Computation \textbf{318}, 245--259 (2018).
\newblock \doi{https://doi.org/10.1016/j.amc.2017.05.014}.
\newblock
  \urlprefix\url{https://www.sciencedirect.com/science/article/pii/S0096300317303028}.
\newblock Recent Trends in Numerical Computations: Theory and Algorithms

\bibitem{Liu2015b}
Liu, H., Xu, S., Wang, X., Wu, X., Song, Y.: A global optimization algorithm
  for simulation-based problems via the extended direct scheme.
\newblock Engineering Optimization \textbf{47}(11), 1441--1458 (2015).
\newblock \doi{10.1080/0305215X.2014.971777}

\bibitem{Liu2019}
Liu, J., Ploskas, N., Sahinidis, N.V.: {Tuning BARON using derivative-free
  optimization algorithms}.
\newblock Journal of Global Optimization \textbf{74}(4), 611--637 (2019).
\newblock \doi{10.1007/s10898-018-0640-3}.
\newblock \urlprefix\url{https://doi.org/10.1007/s10898-018-0640-3}

\bibitem{Liu2013}
Liu, Q.: Linear scaling and the direct algorithm.
\newblock Journal of Global Optimization \textbf{56}, 1233--1245 (2013).
\newblock \doi{10.1007/s10898-012-9952-x}

\bibitem{Liu2015}
Liu, Q., Zeng, J., Yang, G.: {MrDIRECT: a multilevel robust DIRECT algorithm
  for global optimization problems}.
\newblock Journal of Global Optimization \textbf{62}(2), 205--227 (2015).
\newblock \doi{10.1007/s10898-014-0241-8}

\bibitem{Tche2022}
LIU, T.: Mic (2022).
\newblock \urlprefix\url{https://github.com/TcheL/MIC}.
\newblock Retrieved March 6, 2022

\bibitem{Liuzzi2010}
Liuzzi, G., Lucidi, S., Piccialli, V.: A {DIRECT}-based approach exploiting
  local minimizations for the solution of large-scale global optimization
  problems.
\newblock Computational Optimization and Applications \textbf{45}, 353--375
  (2010).
\newblock \doi{10.1007/s10589-008-9217-2}

\bibitem{MATLAB2022}
MATLAB: Matlab and global optimization toolbox release 2022a (2022)

\bibitem{Mehrabian2006}
Mehrabian, A., Lucas, C.: A novel numerical optimization algorithm inspired
  from weed colonization.
\newblock Ecological Informatics \textbf{1}(4), 355--366 (2006).
\newblock \doi{https://doi.org/10.1016/j.ecoinf.2006.07.003}.
\newblock
  \urlprefix\url{https://www.sciencedirect.com/science/article/pii/S1574954106000665}

\bibitem{Metropolis1949}
Metropolis, N., Ulam, S.: The monte carlo method.
\newblock Journal of the American Statistical Association \textbf{44}(247),
  335--341 (1949).
\newblock \doi{10.1080/01621459.1949.10483310}.
\newblock
  \urlprefix\url{https://www.tandfonline.com/doi/abs/10.1080/01621459.1949.10483310}.
\newblock PMID: 18139350

\bibitem{sta2022}
Michael: An improved sta (2022).
\newblock
  \urlprefix\url{https://www.mathworks.com/matlabcentral/fileexchange/69721-an-improved-sta}.
\newblock MATLAB Central File Exchange. Retrieved May 9, 2022.

\bibitem{Mockus2017}
Mockus, J., Paulavi{\v{c}}ius, R., Rusakevi{\v{c}}ius, D., {\v{S}}e{\v{s}}ok,
  D., {\v{Z}}ilinskas, J.: {Application of Reduced-set Pareto-Lipschitzian
  Optimization to truss optimization}.
\newblock Journal of Global Optimization \textbf{67}(1-2), 425--450 (2017).
\newblock \doi{10.1007/s10898-015-0364-6}

\bibitem{Mongeau2000}
Mongeau, M., Karsenty, H., Rouzé, V., Hiriart-Urruty, J.: Comparison of
  public-domain software for black box global optimization.
\newblock Optimization Methods and Software \textbf{13}(3), 203--226 (2000).
\newblock \doi{10.1080/10556780008805783}.
\newblock \urlprefix\url{https://doi.org/10.1080/10556780008805783}

\bibitem{More2009}
{Mor\'{e}}, J.J., Wild, S.M.: Benchmarking derivative-free optimization
  algorithms.
\newblock SIAM Journal on Optimization \textbf{20}(1), 172--191 (2009).
\newblock \doi{10.1137/080724083}

\bibitem{Arnold1999}
Neumaier, A.: Global optimization by multilevel coordinate search (1999).
\newblock \urlprefix\url{https://www.mat.univie.ac.at/~neum/software/mcs/}.
\newblock Retrieved March 21, 2022

\bibitem{Patel1989}
Patel, N.R., Smith, R.L., Zabinsky, Z.B.: {Pure adaptive search in monte carlo
  optimization}.
\newblock Mathematical Programming \textbf{43}(1), 317--328 (1989).
\newblock \doi{10.1007/BF01582296}.
\newblock \urlprefix\url{https://doi.org/10.1007/BF01582296}

\bibitem{Paulavicius2016:jogo}
Paulavi{\v{c}}ius, R., Chiter, L., {\v{Z}}ilinskas, J.: {Global optimization
  based on bisection of rectangles, function values at diagonals, and a set of
  Lipschitz constants}.
\newblock Journal of Global Optimization \textbf{71}(1), 5--20 (2018).
\newblock \doi{10.1007/s10898-016-0485-6}

\bibitem{Paulavicius2014:jogo}
Paulavi{\v{c}}ius, R., Sergeyev, Y.D., Kvasov, D.E., {\v{Z}}ilinskas, J.:
  {Globally-biased DISIMPL algorithm for expensive global optimization}.
\newblock Journal of Global Optimization \textbf{59}(2-3), 545--567 (2014).
\newblock \doi{10.1007/s10898-014-0180-4}

\bibitem{Paulavicius2019:eswa}
Paulavi{\v{c}}ius, R., Sergeyev, Y.D., Kvasov, D.E., {\v{Z}}ilinskas, J.:
  {Globally-biased BIRECT algorithm with local accelerators for expensive
  global optimization}.
\newblock Expert Systems with Applications \textbf{144}, 11305 (2020).
\newblock \doi{10.1016/j.eswa.2019.113052}

\bibitem{Paulavicius2006}
Paulavi{\v{c}}ius, R., {\v{Z}}ilinskas, J.: {Analysis of different norms and
  corresponding Lipschitz constants for global optimization}.
\newblock Technological and Economic Development of Economy \textbf{36}(4),
  383--387 (2006).
\newblock \doi{10.1080/13928619.2006.9637758}

\bibitem{Paulavicius2007}
Paulavi{\v{c}}ius, R., {\v{Z}}ilinskas, J.: {Analysis of different norms and
  corresponding Lipschitz constants for global optimization in multidimensional
  case}.
\newblock Information Technology and Control \textbf{36}(4), 383--387 (2007)

\bibitem{Paulavicius2008}
Paulavi{\v{c}}ius, R., {\v{Z}}ilinskas, J.: {Improved Lipschitz bounds with the
  first norm for function values over multidimensional simplex}.
\newblock Mathematical Modelling and Analysis \textbf{13}(4), 553--563 (2008).
\newblock \doi{10.3846/1392-6292.2008.13.553-563}

\bibitem{Paulavicius2009b}
Paulavi{\v{c}}ius, R., {\v{Z}}ilinskas, J.: {Global optimization using the
  branch-and-bound algorithm with a combination of Lipschitz bounds over
  simplices}.
\newblock Technological and Economic Development of Economy \textbf{15}(2),
  310--325 (2009).
\newblock \doi{10.3846/1392-8619.2009.15.310-325}

\bibitem{Paulavicius2013:jogo}
Paulavi{\v{c}}ius, R., {\v{Z}}ilinskas, J.: {Simplicial Lipschitz optimization
  without the Lipschitz constant}.
\newblock Journal of Global Optimization \textbf{59}(1), 23--40 (2013).
\newblock \doi{10.1007/s10898-013-0089-3}

\bibitem{Paulavicius2014:book}
Paulavi{\v{c}}ius, R., {\v{Z}}ilinskas, J.: {Simplicial Global Optimization}.
\newblock SpringerBriefs in Optimization. Springer New York, New York, NY
  (2014).
\newblock \doi{10.1007/978-1-4614-9093-7}

\bibitem{Pham2006}
Pham, D., Ghanbarzadeh, A., Koç, E., Otri, S., Rahim, S., Zaidi, M.: The bees
  algorithm -- a novel tool for complex optimisation problems.
\newblock In: D.~Pham, E.~Eldukhri, A.~Soroka (eds.) Intelligent Production
  Machines and Systems, pp. 454--459. Elsevier Science Ltd, Oxford (2006).
\newblock \doi{https://doi.org/10.1016/B978-008045157-2/50081-X}.
\newblock
  \urlprefix\url{https://www.sciencedirect.com/science/article/pii/B978008045157250081X}

\bibitem{Pham2011}
Pham, N., Malinowski, A., Bartczak, T.: Comparative study of derivative free
  optimization algorithms.
\newblock IEEE Transactions on Industrial Informatics \textbf{7}(4), 592--600
  (2011).
\newblock \doi{10.1109/TII.2011.2166799}

\bibitem{Pinter1996book}
Pint{\'e}r, J.D.: Global optimization in action: continuous and Lipschitz
  optimization: algorithms, implementations and applications, \emph{Nonconvex
  Optimization and Its Applications}, vol.~6.
\newblock Springer US (1996).
\newblock \doi{10.1007/978-1-4757-2502-5}

\bibitem{Piyavskii1967}
Piyavskii, S.A.: An algorithm for finding the absolute minimum of a function.
\newblock Theory of Optimal Solutions \textbf{2}, 13--24 (1967).
\newblock \doi{10.1016/0041-5553(72)90115-2}.
\newblock In Russian

\bibitem{William2007}
Press, W.H., Teukolsky, S.A., Vetterling, W.T., Flannery, B.P.: Numerical
  Recipes 3rd Edition: The Art of Scientific Computing, 3 edn.
\newblock Cambridge University Press, USA (2007)

\bibitem{Rao2011}
Rao, R., Savsani, V., Vakharia, D.: Teaching–learning-based optimization: A
  novel method for constrained mechanical design optimization problems.
\newblock Computer-Aided Design \textbf{43}(3), 303--315 (2011).
\newblock \doi{https://doi.org/10.1016/j.cad.2010.12.015}.
\newblock
  \urlprefix\url{https://www.sciencedirect.com/science/article/pii/S0010448510002484}

\bibitem{Regis2007}
Regis, R.G., Shoemaker, C.A.: A stochastic radial basis function method for the
  global optimization of expensive functions.
\newblock INFORMS Journal on Computing \textbf{19}(4), 497--509 (2007).
\newblock \doi{10.1287/ijoc.1060.0182}.
\newblock \urlprefix\url{https://doi.org/10.1287/ijoc.1060.0182}

\bibitem{Reynolds2008}
Reynolds, R.G.: An introduction to cultural algorithms (2008)

\bibitem{Rios2013}
Rios, L.M., Sahinidis, N.V.: {Derivative-free optimization: a review of
  algorithms and comparison of software implementations}.
\newblock Journal of Global Optimization \textbf{56}(3), 1247--1293 (2007).
\newblock \doi{10.1007/s10898-012-9951-y}

\bibitem{Sadollah2018}
Sadollah, A., Sayyaadi, H., Yadav, A.: A dynamic metaheuristic optimization
  model inspired by biological nervous systems: Neural network algorithm.
\newblock Applied Soft Computing \textbf{71}, 747--782 (2018).
\newblock \doi{https://doi.org/10.1016/j.asoc.2018.07.039}.
\newblock
  \urlprefix\url{https://www.sciencedirect.com/science/article/pii/S1568494618304277}

\bibitem{Sadollah2018dynamic}
Sadollah, A., Sayyaadi, H., Yadav, A.: A dynamic metaheuristic optimization
  model inspired by biological nervous systems: Neural network algorithm.
\newblock Applied Soft Computing \textbf{71}, 747--782 (2018)

\bibitem{Saravanakumar2015}
Saravanakumar, G., Valarmathi, K., Rajasekaran, M.P., Srinivasan, S., Kadhar,
  K.M.A.: State transition algorithm (sta) based tuning of integer and
  fractional-order pid controller for benchmark system.
\newblock In: 2015 IEEE International Conference on Computational Intelligence
  and Computing Research (ICCIC), pp. 1--5 (2015).
\newblock \doi{10.1109/ICCIC.2015.7435814}

\bibitem{Sergeyev1998a}
Sergeyev, Y.D.: Global one-dimensional optimization using smooth auxiliary
  functions.
\newblock Mathematical Programming \textbf{81}(1), 127--146 (1998).
\newblock \doi{10.1007/BF01584848}

\bibitem{Sergeyev2006}
Sergeyev, Y.D., Kvasov, D.E.: Global search based on diagonal partitions and a
  set of {L}ipschitz constants.
\newblock SIAM Journal on Optimization \textbf{16}(3), 910--937 (2006).
\newblock \doi{10.1137/040621132}

\bibitem{Sergeyev2008:book}
Sergeyev, Y.D., Kvasov, D.E.: Diagonal Global Optimization Methods.
\newblock FizMatLit, Moscow (2008).
\newblock In Russian

\bibitem{Sergeyev2011}
Sergeyev, Y.D., Kvasov, D.E.: Lipschitz global optimization.
\newblock In: J.J. Cochran, L.A. Cox, P.~Keskinocak, J.P. Kharoufeh, J.C. Smith
  (eds.) Wiley Encyclopedia of Operations Research and Management Science (in 8
  volumes), vol.~4, pp. 2812--2828. John Wiley \& Sons, New York (2011)

\bibitem{Sergeyev2021gkls}
Sergeyev, Y.D., Kvasov, D.E., Mukhametzhanov, M.S.: A generator of
  multiextremal test classes with known solutions for black-box constrained
  global optimization.
\newblock IEEE Transactions on Evolutionary Computation pp. 1--1 (2021).
\newblock \doi{10.1109/TEVC.2021.3139263}

\bibitem{amir2022}
amir seyyedabbasi: Sand cat swarm optimization (2022).
\newblock
  \urlprefix\url{https://www.mathworks.com/matlabcentral/fileexchange/110185-sand-cat-swarm-optimization}.
\newblock MATLAB Central File Exchange. Retrieved May 9, 2022.

\bibitem{Seyyedabbasi2022}
Seyyedabbasi, A., Kiani, F.: {Sand Cat swarm optimization: a nature-inspired
  algorithm to solve global optimization problems}.
\newblock Engineering with Computers  (2022).
\newblock \doi{10.1007/s00366-022-01604-x}.
\newblock \urlprefix\url{https://doi.org/10.1007/s00366-022-01604-x}

\bibitem{Shi2021}
Shi, H.J.M., Xuan, M.Q., Oztoprak, F., Nocedal, J.: On the numerical
  performance of derivative-free optimization methods based on
  finite-difference approximations (2021).
\newblock \doi{10.48550/ARXIV.2102.09762}.
\newblock \urlprefix\url{https://arxiv.org/abs/2102.09762}

\bibitem{Shubert1972}
Shubert, B.O.: A sequential method seeking the global maximum of a function.
\newblock SIAM Journal on Numerical Analysis \textbf{9}, 379--388 (1972).
\newblock \doi{10.1137/0709036}

\bibitem{Dan2008}
Simon, D.: Biogeography-based optimization.
\newblock IEEE Transactions on Evolutionary Computation \textbf{12}(6),
  702--713 (2008).
\newblock \doi{10.1109/TEVC.2008.919004}

\bibitem{Storn1997}
Storn, R., Price, K.: {Differential Evolution – A Simple and Efficient
  Heuristic for global Optimization over Continuous Spaces}.
\newblock Journal of Global Optimization \textbf{11}(4), 341--359 (1997).
\newblock \doi{10.1023/A:1008202821328}.
\newblock \urlprefix\url{https://doi.org/10.1023/A:1008202821328}

\bibitem{DIRECTGO2022}
Stripinis, L., Paulavi{\v c}ius, R.: {DIRECTGO: A new DIRECT-type MATLAB
  toolbox for derivative-free global optimization} (2022).
\newblock \urlprefix\url{https://github.com/blockchain-group/DIRECTGO}

\bibitem{Stripinis2021b}
Stripinis, L., Paulavi{\v{c}}ius, R.: An empirical study of various candidate
  selection and partitioning techniques in the direct framework.
\newblock Journal of Global Optimization  (2022).
\newblock \doi{10.1007/s10898-022-01185-5}.
\newblock \urlprefix\url{https://doi.org/10.1007/s10898-022-01185-5}

\bibitem{Stripinis2018a}
Stripinis, L., Paulavi{\v{c}}ius, R., {\v{Z}}ilinskas, J.: Improved scheme for
  selection of potentially optimal hyper-rectangles in {DIRECT}.
\newblock Optimization Letters \textbf{12}(7), 1699--1712 (2018).
\newblock \doi{10.1007/s11590-017-1228-4}

\bibitem{Stripinis2018b}
Stripinis, L., Paulavi{\v{c}}ius, R., {\v{Z}}ilinskas, J.: Penalty functions
  and two-step selection procedure based \direct-type algorithm for constrained
  global optimization.
\newblock Structural and Multidisciplinary Optimization \textbf{59}(6),
  2155--2175 (2019).
\newblock \doi{10.1007/s00158-018-2181-2}

\bibitem{Stripinis2021c}
Stripinis, L., Paulavi\v{c}ius, R.: Directgo: A new direct-type matlab toolbox
  for derivative-free global optimization.
\newblock ACM Transactions on Mathematical Software  (2022).
\newblock \doi{10.1145/3559755}.
\newblock \urlprefix\url{https://doi.org/10.1145/3559755}

\bibitem{DIRECTGOLib2022}
Stripinis, L., Paulavi\v{c}ius, R.: {DIRECTGOLib - DIRECT Global Optimization
  test problems Library} (2022).
\newblock \doi{10.5281/zenodo.6617799}.
\newblock \urlprefix\url{https://doi.org/10.5281/zenodo.6617799}

\bibitem{Derek2013}
Surjanovic, S., Bingham, D.: Virtual library of simulation experiments: Test
  functions and datasets.
\newblock \url{http://www.sfu.ca/~ssurjano/index.html} (2013).
\newblock Online; accessed: 2017-03-22

\bibitem{Tsafarakis2020}
Tsafarakis, S., Zervoudakis, K., Andronikidis, A., Altsitsiadis, E.: Fuzzy
  self-tuning differential evolution for optimal product line design.
\newblock European Journal of Operational Research \textbf{287}(3), 1161--1169
  (2020).
\newblock \doi{https://doi.org/10.1016/j.ejor.2020.05.018}.
\newblock
  \urlprefix\url{https://www.sciencedirect.com/science/article/pii/S037722172030446X}

\bibitem{Ugray2007}
Ugray, Z., Lasdon, L., Plummer, J., Glover, F., Kelly, J., Martí, R.: {Scatter
  Search and Local NLP Solvers: A Multistart Framework for Global
  Optimization}.
\newblock INFORMS Journal on Computing \textbf{19}(3), 328--340 (2007).
\newblock \doi{10.1287/ijoc.1060.0175}.
\newblock \urlprefix\url{https://doi.org/10.1287/ijoc.1060.0175}

\bibitem{zilinskas_zhigljavsky_2016}
\v{Z}ilinskas, A., Zhigljavsky, A.: Stochastic global optimization: A review on
  the occasion of 25 years of informatica.
\newblock Informatica \textbf{27}(2), 229--256 (2016).
\newblock \doi{10.15388/Informatica.2016.83}

\bibitem{Xi2020}
Xi, M., Sun, W., Chen, J.: Survey of derivative-free optimization.
\newblock Numerical Algebra, Control \& Optimization \textbf{10}(4), 537--555
  (2020)

\bibitem{Yang2022b}
Yang, X.: The standard bat algorithm (ba) (2022).
\newblock
  \urlprefix\url{https://www.mathworks.com/matlabcentral/fileexchange/74768-the-standard-bat-algorithm-ba}.
\newblock MATLAB Central File Exchange. Retrieved February 18, 2022.

\bibitem{Yang2022a}
Yang, X.: The standard cuckoo search (cs) (2022).
\newblock
  \urlprefix\url{https://www.mathworks.com/matlabcentral/fileexchange/74767-the-standard-cuckoo-search-cs}.
\newblock MATLAB Central File Exchange. Retrieved February 18, 2022.

\bibitem{Yang2009as}
Yang, X.S.: Firefly algorithms for multimodal optimization.
\newblock In: O.~Watanabe, T.~Zeugmann (eds.) Stochastic Algorithms:
  Foundations and Applications, pp. 169--178. Springer Berlin Heidelberg,
  Berlin, Heidelberg (2009)

\bibitem{Yang2010as}
Yang, X.S.: A New Metaheuristic Bat-Inspired Algorithm, pp. 65--74.
\newblock Springer Berlin Heidelberg, Berlin, Heidelberg (2010).
\newblock \doi{10.1007/978-3-642-12538-6_6}.
\newblock \urlprefix\url{https://doi.org/10.1007/978-3-642-12538-6_6}

\bibitem{Yang2009a}
Yang, X.S., Deb, S.: Cuckoo search via lévy flights.
\newblock In: 2009 World Congress on Nature Biologically Inspired Computing
  (NaBIC), pp. 210--214 (2009).
\newblock \doi{10.1109/NABIC.2009.5393690}

\bibitem{Yarpiz2022}
Yarpiz: Shuffled frog leaping algorithm (sfla) (2022).
\newblock
  \urlprefix\url{https://www.mathworks.com/matlabcentral/fileexchange/52861-shuffled-frog-leaping-algorithm-sfla}

\bibitem{Zavala2014}
Zavala, G.R., Nebro, A.J., Luna, F., {Coello Coello}, C.A.: {A survey of
  multi-objective metaheuristics applied to structural optimization}.
\newblock Structural and Multidisciplinary Optimization \textbf{49}(4),
  537--558 (2014).
\newblock \doi{10.1007/s00158-013-0996-4}.
\newblock \urlprefix\url{https://doi.org/10.1007/s00158-013-0996-4}

\bibitem{Zhang2022}
yiying zhang: Cclnna for global optimization (2022).
\newblock
  \urlprefix\url{https://www.mathworks.com/matlabcentral/fileexchange/100601-cclnna-for-global-optimization}.
\newblock MATLAB Central File Exchange. Retrieved February 18, 2022.

\bibitem{Zhigljavsky2008:book}
Zhigljavsky, A., {\v Z}ilinskas, A.: Stochastic Global Optimization.
\newblock Springer, New York (2008)

\end{thebibliography}

\section{Appendix - \directgolib{} library}\label{apendix:figures}

\begin{table*}[htbp]
	\caption{Key characteristics of box-constrained global optimization problems (with varying $n$) from \directgolib~\cite{DIRECTGOLib2022}.}
	\resizebox{\textwidth}{!}{
		\fontsize{10pt}{8pt}\selectfont
		\begin{tabular}{clrrR{2.25cm}rrR{1.75cm}}
			\toprule
			\# & Name & Source & \multicolumn{1}{r}{$ D $}  &  \multicolumn{1}{r}{$ \tilde{D} $} & Type & No. of minima & $ f^* $ \\
			\midrule
			$1$ & \textit{Ackley}$^{\alpha}$ & \cite{Hedar2005,Derek2013} & $[-15, 30]^n$ & $[-18, 47]^n$ & non-convex & multi-modal & $0.0000$ \\
			$2$ & \textit{AlpineN1}$^{\alpha}$$^{\beta}$ & \cite{Gavana2021} & $[-10, 10]^n$ & $[-10, 7.5]^n$ & non-convex & multi-modal & $0.0000$ \\
			$3$ & \textit{Alpine}$^{\alpha}$ & \cite{Gavana2021} & $[0, 10]^n$ & $[\sqrt[i]{2}, 8 + \sqrt[i]{2}]^i$ & non-convex & multi-modal & $-2.8081^n$ \\
			$4$ & \textit{Brown}$^{\beta}$ & \cite{Jamil2013} & $[-1, 4]^n$ & $-$ & convex & uni-modal & $0.0000$ \\
			$5$ & \textit{ChungR}$^{\beta}$ & \cite{Jamil2013} & $[-100, 350]^n$ & $-$ & convex & uni-modal & $0.0000$ \\
			$6$ & \textit{Csendes}$^{\alpha}$ & \cite{Gavana2021} & $[-10, 10]^n$ & $[-10, 25]^n$ & convex & multi-modal & $0.0000$ \\
			$7$ & \textit{Cubic}$^{\beta}$ & \cite{Derek2013} & $[-4, 3]^n$ & $-$ & convex & uni-modal & $0.0000$ \\
			$8$ & \textit{Deb01}$^{\alpha}$ & \cite{Gavana2021} & $[-1, 1]^n$ & $[-0.55, 1.45]^n$ & non-convex & multi-modal & $-1.0000$ \\
			$9$ & \textit{Deb02}$^{\alpha}$ & \cite{Gavana2021} & $[0, 1]^n$ & $[0.225, 1.225]^n$ & non-convex & multi-modal & $-1.0000$ \\
			$10$ & \textit{Dixon\_and\_Price} & \cite{Hedar2005,Derek2013} & $[-10, 10]^n$ & $-$ & convex & multi-modal & $0.0000$ \\
			$11$ & \textit{Dejong}$^{\beta}$ & \cite{Derek2013} & $[-3, 7]^n$ & $-$ & convex & uni-modal & $0.0000$ \\
			$12$ & \textit{Exponential}$^{\beta}$ & \cite{Jamil2013} & $[-1, 4]^n$ & $-$ & non-convex & multi-modal & $-1.0000$ \\
			$13$ & \textit{Exponential2}$^{\beta}$ & \cite{Derek2013} & $[0, 7]^n$ & $-$ & non-convex & multi-modal & $0.0000$ \\
			$14$ & \textit{Exponential3}$^{\beta}$ & \cite{Derek2013} & $[-30, 20]^n$ & $-$ & non-convex & multi-modal & $0.0000$ \\
			$15$ & \textit{Griewank}$^{\alpha}$ & \cite{Hedar2005,Derek2013} & $[-600, 600 ]^i$ & $[-\sqrt{600i}, 600\sqrt{i}^{-1} ]^i$ & non-convex & multi-modal & $0.0000$ \\
			$16$ & \textit{Layeb01}$^{\alpha}$$^{\beta}$ & \cite{Layeb2022} & $[-100, 100]^n$ & $[-100, 90]^n$ & convex & uni-modal & $0.0000$ \\
			$17$ & \textit{Layeb02}$^{\beta}$ & \cite{Layeb2022} & $[-10, 10]^n$ & $-$ & convex & uni-modal & $0.0000$ \\
			$18$ & \textit{Layeb03}$^{\alpha}$$^{\beta}$ & \cite{Layeb2022} & $[-10, 10]^n$ & $[-10, 12]^n$ & non-convex & multi-modal & $-n+1$ \\
			$19$ & \textit{Layeb04}$^{\beta}$ & \cite{Layeb2022} & $[-10, 10]^n$ & $-$ & non-convex & multi-modal & $(ln(0.001)-1)$$(n-1)$ \\
			$20$ & \textit{Layeb05}$^{\beta}$ & \cite{Layeb2022} & $[-10, 10]^n$ & $-$ & non-convex & multi-modal & $(ln(0.001))$ $(n - 1)$ \\
			$21$ & \textit{Layeb06}$^{\beta}$ & \cite{Layeb2022} & $[-10, 10]^n$ & $-$ & non-convex & multi-modal & $0.0000$ \\
			$22$ & \textit{Layeb07}$^{\alpha}$$^{\beta}$ & \cite{Layeb2022} & $[-10, 10]^n$ & $[-10, 12]^n$ & non-convex & multi-modal & $0.0000$ \\
			$23$ & \textit{Layeb08}$^{\beta}$ & \cite{Layeb2022} & $[-10, 10]^n$ & $-$ & non-convex & multi-modal & $log(0.001)(n-1)$ \\
			$24$ & \textit{Layeb09}$^{\beta}$ & \cite{Layeb2022} & $[-10, 10]^n$ & $-$ & non-convex & multi-modal & $0.0000$ \\
			$25$ & \textit{Layeb10}$^{\beta}$ & \cite{Layeb2022} & $[-100, 100]^n$ & $-$ & non-convex & multi-modal & $0.0000$ \\
			$26$ & \textit{Layeb11}$^{\beta}$ & \cite{Layeb2022} & $[-10, 10]^n$ & $-$ & non-convex & multi-modal & $n-1$ \\
			$27$ & \textit{Layeb12}$^{\beta}$ & \cite{Layeb2022} & $[-5, 5]^n$ & $-$ & non-convex & multi-modal & $-(e+1)(n-1)$ \\
			$28$ & \textit{Layeb13}$^{\beta}$ & \cite{Layeb2022} & $[-5, 5]^n$ & $-$ & non-convex & multi-modal & $0.0000$ \\
			$29$ & \textit{Layeb14}$^{\beta}$ & \cite{Layeb2022} & $[-100, 100]^n$ & $-$ & non-convex & multi-modal & $0.0000$ \\
			$30$ & \textit{Layeb15}$^{\beta}$ & \cite{Layeb2022} & $[-100, 100]^n$ & $-$ & non-convex & multi-modal & $0.0000$ \\
			$31$ & \textit{Layeb16}$^{\beta}$ & \cite{Layeb2022} & $[-10, 10]^n$ & $-$ & non-convex & multi-modal & $0.0000$ \\
			$32$ & \textit{Layeb17}$^{\beta}$ & \cite{Layeb2022} & $[-10, 10]^n$ & $-$ & non-convex & multi-modal & $0.0000$ \\
			$33$ & \textit{Layeb18}$^{\beta}$ & \cite{Layeb2022} & $[-10, 10]^n$ & $-$ & non-convex & multi-modal & $ln(0.001)(n-1)$ \\
			$34$ & \textit{Levy} & \cite{Hedar2005,Derek2013} & $[-5, 5]^n$ & $[-10, 10]^n$ & non-convex & multi-modal & $0.0000$ \\
			$35$ & \textit{Michalewicz} & \cite{Hedar2005,Derek2013} & $[0, \pi]^n$ & $-$ & non-convex & multi-modal & $\chi$ \\
			$36$ & \textit{Pinter}$^{\alpha}$ & \cite{Gavana2021} & $[-10, 10]^n$ & $[-5.5, 14.5]^n$ & non-convex & multi-modal & $0.0000$ \\
			$37$ & \textit{Qing} & \cite{Gavana2021} & $[-500, 500]^n$ & $-$ & non-convex & multi-modal & $0.0000$ \\
			$38$ & \textit{Quadratic}$^{\beta}$ & \cite{Derek2013} & $[-2, 3]^n$ & $-$ & convex & uni-modal & $0.0000$ \\
			$39$ & \textit{Rastrigin}$^{\alpha}$ & \cite{Hedar2005,Derek2013} & $[-5.12, 5.12]^n$ & $[-5\sqrt[i]{2}, 7+\sqrt[i]{2}]^i$ & non-convex & multi-modal & $0.0000$ \\
			$40$ & \textit{Rosenbrock}$^{\alpha}$ & \cite{Hedar2005,Dixon1978} & $[-5, 10]^n$ & $[-5\sqrt{i}^{-1}, 10\sqrt{i} ]^i$ & non-convex & uni-modal & $0.0000$ \\
			$41$ & \textit{Rotated\_H\_Ellip}$^{\alpha}$ & \cite{Derek2013} & $[-65.536, 65.536]^n$ & $[-35, 95]^n$ & convex & uni-modal & $0.0000$ \\
			$42$ & \textit{Schwefel}$^{\alpha}$ & \cite{Hedar2005,Derek2013} & $[-500, 500]^n$ & $[-500 + 100\sqrt{i}^{-1},$ $500 - 40\sqrt{i}^{-1} ]^i$ & non-convex & multi-modal & $0.0000$ \\
			$43$ & \textit{SineEnvelope}$^{\beta}$ & \cite{Gavana2021} & $[-100, 100]^n$ & $-$ & non-convex & multi-modal & $-2.6535(n-1)$ \\
			$44$ & \textit{Sinenvsin}$^{\alpha}$$^{\beta}$ & \cite{Mishra2006} & $[-100, 100]^n$ & $[-100, 150]^n$ & non-convex & multi-modal & $0.0000$ \\
			$45$ & \textit{Sphere}$^{\alpha}$ & \cite{Hedar2005,Derek2013} & $[-5.12, 5.12]^n$ & $[-2.75, 7.25]^n$ & convex & uni-modal & $0.0000$ \\
			$46$ & \textit{Styblinski\_Tang} & \cite{Clerc1999} & $[-5, 5]^n$ & $[-5, 5+3^{1/i}]^n$ & non-convex & multi-modal & $-39.1661n$ \\
			$47$ & \textit{Sum\_Squares}$^{\alpha}$ & \cite{Clerc1999} & $[-10, 10]^n$ & $[-5.5, 14.5]^n$ & convex & uni-modal & $0.0000$ \\
			$48$ & \textit{Sum\_Of\_Powers}$^{\alpha}$ & \cite{Derek2013} & $[-1, 1]^n$ & $[-0.55, 1.45]^n$ & convex & uni-modal & $0.0000$ \\
			$49$ & \textit{Trid} & \cite{Hedar2005,Derek2013} & $[-100, 100]^n$ & $-$ & convex & multi-modal & $-\frac{1}{6}n^3 - \frac{1}{2}n^2 + \frac{2}{3}n$ \\
			$50$ & \textit{Trigonometric}$^{\alpha}$$^{\beta}$ & \cite{Hedar2005,Derek2013} & $[-100, 100]^n$ & $[-100, 150]^n$ & non-convex & multi-modal & $0.0000$ \\
			$51$ & \textit{Vincent} & \cite{Clerc1999} & $[0.25, 10]^n$ & $-$ & non-convex & multi-modal & $-n$ \\
			$52$ & \textit{WWavy}$^{\alpha}$$^{\beta}$ & \cite{Jamil2013} & $[-\pi, \pi]^n$ & $[-\pi, 3\pi]^n$ & non-convex & multi-modal & $0.0000$ \\
			$53$ & \textit{XinSheYajngN1}$^{\beta}$ & \cite{Jamil2013} & $[-11, 29]^n$ & $[-11, 29]^n$ & non-convex & multi-modal & $-1.0000$ \\
			$54$ & \textit{XinSheYajngN2}$^{\alpha}$$^{\beta}$ & \cite{Jamil2013} & $[-\pi, \pi]^n$ & $[-\pi, 3\pi]^n$ & non-convex & multi-modal & $0.0000$ \\
			$55$ & \textit{Zakharov}$^{\alpha}$ & \cite{Hedar2005,Derek2013} & $[-5, 10]^n$ & $[-1.625, 13.375]^n$ & convex & multi-modal & $0.0000$ \\
			\midrule
			\multicolumn{8}{l}{$i$ -- indexes used for variable bounds ($1, ..., n$)} \\
			\multicolumn{8}{l}{$\chi$ -- solution depend on problem dimension} \\
			\multicolumn{8}{l}{$\alpha$ -- domain $D$ was perturbed} \\
			\multicolumn{8}{l}{$\beta$ -- New test problem in \directgolib} \\
			\bottomrule
	\end{tabular}}
	\label{tab:test}
\end{table*}

\begin{table*}[h!]
	\centering
	\caption{Key characteristics of box-constrained global optimization problems (with fixed $n$) from \directgolib~\cite{DIRECTGOLib2022}.}
	\resizebox{\textwidth}{!}{
		\fontsize{10pt}{8pt}\selectfont
		\begin{tabular}{clrrR{2.25cm}R{2.25cm}rrr}
			\toprule
			\# & Name & Source & $n$ & \multicolumn{1}{r}{$ D $}  &  \multicolumn{1}{r}{$ \tilde{D} $} & Type & No. of minima & $ f^* $ \\
			\midrule
			$1$ & \textit{AckleyN2}$^{\alpha}$$^{\beta}$ & \cite{Jamil2013} & $2$ & $[-15, 30]^n$ & $[-18, 47]^n$ & uni-modal & convex & $-200$ \\
			$2$ & \textit{AckleyN3}$^{\alpha}$$^{\beta}$ & \cite{Jamil2013} & $2$ & $[-15, 30]^n$ & $[-18, 47]^n$ & uni-modal & convex & $-186.4112$ \\
			$3$ & \textit{AckleyN4}$^{\alpha}$$^{\beta}$ & \cite{Jamil2013} & $2$ & $[-15, 30]^n$ & $[-18, 47]^n$ & non-convex & multi-modal & $-4.5901$ \\
			$4$ & \textit{Adjiman}$^{\beta}$ & \cite{Jamil2013} & $2$ & $[-1, 2]^n$ & $-$ & non-convex & multi-modal & $-2.0218$ \\
			$5$ & \textit{BartelsConn}$^{\alpha}$$^{\beta}$ & \cite{Jamil2013} & $2$ & $[-500, 500]^n$ & $[-300, 700]^n$ & non-convex & multi-modal & $1.0000$ \\
			$6$ & \textit{Beale} & \cite{Hedar2005,Derek2013} & $2$ & $[-4.5, 4.5]^n$ & $-$ & non-convex & multi-modal & $0.0000$ \\
			$7$ & \textit{BiggsEXP2}$^{\beta}$ & \cite{Jamil2013} & $2$ & $[0, 20]^n$ & $-$ & non-convex & multi-modal & $0.000$ \\
			$8$ & \textit{BiggsEXP3}$^{\beta}$ & \cite{Jamil2013} & $3$ & $[0, 20]^n$ & $-$ & non-convex & multi-modal & $0.000$ \\
			$9$ & \textit{BiggsEXP4}$^{\beta}$ & \cite{Jamil2013} & $4$ & $[0, 20]^n$ & $-$ & non-convex & multi-modal & $0.000$ \\
			$10$ & \textit{BiggsEXP5}$^{\beta}$ & \cite{Jamil2013} & $5$ & $[0, 20]^n$ & $-$ & non-convex & multi-modal & $0.000$ \\
			$11$ & \textit{BiggsEXP6}$^{\beta}$ & \cite{Jamil2013} & $6$ & $[0, 20]^n$ & $-$ & non-convex & multi-modal & $0.000$ \\
			$12$ & \textit{Bird}$^{\beta}$ & \cite{Jamil2013} & $2$ & $[-2\pi, 2\pi]^n$ & $-$ & non-convex & multi-modal & $-106.7645$ \\
			$13$ & \textit{Bohachevsky1}$^{\alpha}$ & \cite{Hedar2005,Derek2013} & $2$ & $[-100, 100]^n$ & $[-55, 145]^n$ & convex & uni-modal & $0.0000$ \\
			$14$ & \textit{Bohachevsky2}$^{\alpha}$ & \cite{Hedar2005,Derek2013} & $2$ & $[-100, 100]^n$ & $[-55, 145]^n$ & non-convex & multi-modal & $0.0000$ \\
			$15$ & \textit{Bohachevsky3}$^{\alpha}$ & \cite{Hedar2005,Derek2013} & $2$ & $[-100, 100]^n$ & $[-55, 145]^n$ & non-convex & multi-modal & $0.0000$ \\
			$16$ & \textit{Booth} & \cite{Hedar2005,Derek2013} & $2$ & $[-10, 10]^n$ & $-$ & convex & uni-modal & $0.0000$ \\
			$17$ & \textit{Brad}$^{\beta}$ & \cite{Jamil2013} & $3$ & $[-0.25, 0.25] \times [0.01, 2.5]^2$ & $-$ & non-convex & multi-modal & $6.9352$ \\
			$18$ & \textit{Branin} & \cite{Hedar2005,Dixon1978} & $2$ & $[-5, 10] \times [10,15]$ & $-$ & non-convex & multi-modal & $0.3978$ \\
			$19$ & \textit{Bukin4}$^{\beta}$ & \cite{Jamil2013} & $2$ & $[-15, -5] \times [-3,3]$ & $-$ & convex & multi-modal & $0.0000$ \\
			$20$ & \textit{Bukin6} & \cite{Derek2013} & $2$ & $[-15, 5] \times [-3,3]$ & $-$ & convex & multi-modal & $0.0000$ \\
			$21$ & \textit{CarromTable}$^{\beta}$ & \cite{Gavana2021} & $2$ & $[-10, 10]^n$ & $-$ & non-convex & multi-modal & $-24.1568$ \\
			$22$ & \textit{ChenBird}$^{\beta}$ & \cite{Jamil2013} & $2$ & $[-500, 500]^n$ & $-$ & non-convex & multi-modal & $-2000.0000$ \\
			$23$ & \textit{ChenV}$^{\beta}$ & \cite{Jamil2013} & $2$ & $[-500, 500]^n$ & $-$ & non-convex & multi-modal & $-2000.0009$ \\
			$24$ & \textit{Chichinadze}$^{\beta}$ & \cite{Jamil2013} & $2$ & $[-30, 30]^n$ & $-$ & non-convex & multi-modal & $-42.9443$ \\
			$25$ & \textit{Cola}$^{\beta}$ & \cite{Jamil2013} & $17$ & $[-4, 4]^n$ & $-$ & non-convex & multi-modal & $12.0150$ \\
			$26$ & \textit{Colville} & \cite{Hedar2005,Derek2013} & $4$ & $[-10, 10]^n$ & $-$ & non-convex & multi-modal & $0.0000$ \\
			$27$ & \textit{Cross\_function}$^{\beta}$ & \cite{Gavana2021} & $2$ & $[-10, 10]^n$ & $-$ & non-convex & multi-modal & $0.00004$ \\
			$28$ & \textit{Cross\_in\_Tray} & \cite{Derek2013} & $2$ & $[0, 10]^n$ & $-$ & non-convex & multi-modal & $-2.0626$ \\
			$29$ & \textit{CrownedCross}$^{\beta}$ & \cite{Gavana2021} & $2$ & $[-10, 15]^n$ & $-$ & non-convex & multi-modal & $0.0001$ \\
			$30$ & \textit{Crosslegtable} & \cite{Gavana2021} & $2$ & $[-10, 15]^n$ & $-$ & non-convex & multi-modal & $-1.000$ \\
			$31$ & \textit{Cube}$^{\beta}$ & \cite{Gavana2021} & $2$ & $[-10, 10]^n$ & $-$ & convex & multi-modal & $0.000$ \\
			$32$ & \textit{Damavandi} & \cite{Gavana2021} & $2$ & $[0, 14]^n$ & $-$ & non-convex & multi-modal & $0.0000$ \\
			$33$ & \textit{Dejong5}$^{\beta}$ & \cite{Derek2013} & $2$ & $[-65.536, 65.536]^n$ & $-$ & non-convex & multi-modal & $0.9980$ \\
			$34$ & \textit{Dolan}$^{\beta}$ & \cite{Jamil2013} & $5$ & $[-100, 100]^n$ & $-$ & non-convex & multi-modal & $-529.8714$ \\
			$35$ & \textit{Drop\_wave}$^{\alpha}$ & \cite{Derek2013} & $2$ & $[-5.12, 5.12]^n$ & $[-4, 6]^n$ & non-convex & multi-modal & $-1.0000$ \\
			$36$ & \textit{Easom}$^{\alpha}$ & \cite{Hedar2005,Derek2013} & $2$ & $[-100, 100]^n$ & $[-100(i+1)^{-1}, 100i ]^i$ & non-convex & multi-modal & $-1.0000$ \\
			$37$ & \textit{Eggholder} & \cite{Derek2013} & $2$ & $[-512, 512]^n$ & $-$ & non-convex & multi-modal & $-959.6406$ \\
			$38$ & \textit{Giunta}$^{\beta}$ & \cite{Gavana2021} & $2$ & $[-1, 1]^n$ & $-$ & non-convex & multi-modal & $0.0644$ \\
			$39$ & \textit{Goldstein\_and\_Price}$^{\alpha}$ & \cite{Hedar2005,Dixon1978} & $2$ & $[-2, 2]^n$ & $[-1.1, 2.9]^n$ & non-convex & multi-modal & $3.0000$ \\
			$40$ & \textit{Hartman3} & \cite{Hedar2005,Derek2013} & $3$ & $[0, 1]^n$ & $-$ & non-convex & multi-modal & $-3.8627$ \\
			$41$ & \textit{Hartman4}$^{\beta}$ & \cite{Hedar2005,Derek2013} & $4$ & $[0, 1]^n$ & $-$ & non-convex & multi-modal & $-3.1344$ \\
			$42$ & \textit{Hartman6} & \cite{Hedar2005,Derek2013} & $6$ & $[0, 1]^n$ & $-$ & non-convex & multi-modal & $-3.3223$ \\
			$43$ & \textit{HelicalValley}$^{\beta}$ & \cite{Gavana2021} & $3$ & $[-10, 20]^n$ & $-$ & convex & multi-modal & $0.0000$ \\
			$44$ & \textit{HimmelBlau}$^{\beta}$ & \cite{Gavana2021} & $2$ & $[-5, 5]^n$ & $-$ & convex & multi-modal & $0.0000$ \\
			$45$ & \textit{Holder\_Table} & \cite{Derek2013} & $2$ & $[-10, 10]^n$ & $-$ & non-convex & multi-modal & $-19.2085$ \\
			$46$ & \textit{Hump} & \cite{Hedar2005,Derek2013} & $2$ & $[-5, 5]^n$ & $-$ & non-convex & multi-modal & $-1.0316$ \\
			$47$ & \textit{Langermann} & \cite{Derek2013} & $2$ & $[0, 10]^n$ & $-$ & non-convex & multi-modal & $-4.1558$ \\
			$48$ & \textit{Leon}$^{\beta}$ & \cite{Gavana2021} & $2$ & $[-1.2, 1.2]^n$ & $-$ & convex & multi-modal & $0.0000$ \\
			$49$ & \textit{Levi13}$^{\beta}$ & \cite{Gavana2021} & $2$ & $[-10, 10]^n$ & $-$ & non-convex & multi-modal & $0.0000$ \\
			$50$ & \textit{Matyas}$^{\alpha}$ & \cite{Hedar2005,Derek2013} & $2$ & $[-10, 10]^n$ & $[-5.5, 14.5]^n$ & convex & uni-modal & $0.0000$ \\
			$51$ & \textit{McCormick} & \cite{Derek2013} & $2$ & $[-1.5, 4] \times [-3,4]$ & $-$ & convex & multi-modal & $-1.9132$ \\
			$52$ & \textit{ModSchaffer1}$^{\alpha}$$^{\beta}$ & \cite{Mishra2006} & $2$ & $[-100, 100]^n$ & $[-100, 150]^n$ & non-convex & multi-modal & $0.0000$ \\
			$53$ & \textit{ModSchaffer2}$^{\alpha}$$^{\beta}$ & \cite{Mishra2006} & $2$ & $[-100, 100]^n$ & $[-100, 150]^n$ & non-convex & multi-modal & $0.0000$ \\
			$54$ & \textit{ModSchaffer3}$^{\alpha}$$^{\beta}$ & \cite{Mishra2006} & $2$ & $[-100, 100]^n$ & $[-100, 150]^n$ & non-convex & multi-modal & $0.0015$ \\
			$55$ & \textit{ModSchaffer4}$^{\alpha}$$^{\beta}$ & \cite{Mishra2006} & $2$ & $[-100, 100]^n$ & $[-100, 150]^n$ & non-convex & multi-modal & $0.2925$ \\
			$56$ & \textit{PenHolder}$^{\beta}$ & \cite{Hedar2005,Derek2013} & $2$ & $[-11, 11]^n$ & $-$ & non-convex & multi-modal & $-0.9635$ \\
			$57$ & \textit{Permdb4} & \cite{Hedar2005,Derek2013} & $4$ & $[-i, i]^i$ & $-$ & non-convex & multi-modal & $0.0000$ \\
			$58$ & \textit{Powell} & \cite{Hedar2005,Derek2013} & $4$ & $[-4, 5]^n$ & $-$ & convex & multi-modal & $0.0000$ \\
			$59$ & \textit{Power\_Sum}$^{\alpha}$ & \cite{Hedar2005,Derek2013} & $4$ & $[0, 4]^n$ & $[1, 4 + \sqrt[i]{2}]^i$ & convex & multi-modal & $0.0000$ \\
			$60$ & \textit{Shekel5} & \cite{Hedar2005,Derek2013} & $4$ & $[0, 10]^n$ & $-$ & non-convex & multi-modal & $-10.1531$ \\
			$61$ & \textit{Shekel7} & \cite{Hedar2005,Derek2013} & $4$ & $[0, 10]^n$ & $-$ & non-convex & multi-modal & $-10.4029$ \\
			$62$ & \textit{Shekel10} & \cite{Hedar2005,Derek2013} & $4$ & $[0, 10]^n$ & $-$ & non-convex & multi-modal & $-10.5364$ \\
			$63$ & \textit{Shubert} & \cite{Hedar2005,Derek2013} & $2$ & $[-10, 10]^n$ & $-$ & non-convex & multi-modal & $-186.7309$ \\
			$64$ & \textit{TestTubeHolder}$^{\beta}$ & \cite{Gavana2021} & $2$ & $[-10, 10]^n$ & $-$ & non-convex & multi-modal & $-10.8722$ \\
			$65$ & \textit{Trefethen} & \cite{Gavana2021} & $2$ & $[-2, 2]^n$ & $-$ & non-convex & multi-modal & $-3.3068$ \\
			$66$ & \textit{Wood}$^{\alpha}$$^{\beta}$ & \cite{Mishra2006} & $4$ & $[-100, 100]^n$ & $[-100, 150]^n$ & non-convex & multi-modal & $0.0000$ \\
			$67$ & \textit{Zettl}$^{\beta}$ & \cite{Gavana2021} & $2$ & $[-5, 5]^n$ & $-$ & convex & multi-modal & $-0.0037$ \\
			\midrule
			\multicolumn{8}{l}{$i$ -- indexes used for variable bounds ($1, ..., n$)} \\
			\multicolumn{8}{l}{$\alpha$ -- domain $D$ was perturbed} \\
			\multicolumn{8}{l}{$\beta$ -- New test problem in \directgolib} \\
			\bottomrule
	\end{tabular}}
	\label{tab:tests}
\end{table*}

\end{document}